\documentclass[12pt]{amsart}
\usepackage{amssymb}
\usepackage{amsmath}
\usepackage{mathabx}
\usepackage{amsthm}
\usepackage{amstext}
\usepackage{amsopn}
\usepackage{mathrsfs}
\usepackage{latexsym}
\usepackage{enumerate}
\usepackage{bm}
\usepackage{needspace}
\DeclareMathAlphabet{\mathpzc}{OT1}{pzc}{m}{it}
\allowdisplaybreaks
\newcommand{\subheading}[1]{\medskip\noindent{\it #1.}\par\smallskip}

\newtheorem{Theorem}{Theorem}[section]
\newtheorem{Definition}[Theorem]{Definition}
\newtheorem{Proposition}[Theorem]{Proposition}
\newtheorem{Lemma}[Theorem]{Lemma}
\newtheorem{Corollary}[Theorem]{Corollary}
\newtheorem{Remark}[Theorem]{Remark}
\newtheorem{Example}[Theorem]{Example}

\begin{document}
%%%%%%%%%%%%%%%%%%%%%%%%%%%%%%%%%%%%%%%%%%%%%%%%%%%%%%%%%%%%
\bibliographystyle{plain}
%%%%%%%%%%%%%%%%%%%%%%%%%%%%%%%%%%%%%%%%%%%%%%%%%%%%%%%%%%%%%%%%%
\title[Decomposition of free cumulants]
{Decomposition of free cumulants}
\author[R. Lenczewski]{Romuald Lenczewski}
\address{Romuald Lenczewski \newline
Katedra Matematyki, Wydzia\l{} Matematyki, Politechnika Wroc\l{}awska, \newline
Wybrze\.{z}e Wyspia\'{n}skiego 27, 50-370 Wroc{\l}aw, Poland}
\email{Romuald.Lenczewski@pwr.edu.pl}
%%%%%%%%%%%%%%%%%%%%%%%%%%%%%%%%%%%%%%%%%%%%%%%%%%%%%%%%%%%%%%%%
\begin{abstract}
Free cumulants are multilinear functionals defined 
in terms of the moment functional with the use of the family
$\{{\rm NC}(n):n\in \mathbb{N}\}$ of lattices of noncrossing partitions.
In the univariate case, they can be identified with the coefficients of the Voiculescu transform of 
a moment functional which
plays a role similar to that of the logarithm of the Fourier transform.
The associated linearization property 
is connected with free independence. In turn, the family $\{{\rm Int}(n):n\in \mathbb{N}\}$
of lattices of interval partitions gives rise to Boolean cumulants connected with Boolean independence.
In order to bridge the gap between these two families of lattices 
we introduce and study the family $\{\mathcal{NC}(w):w\in \mathpzc{M}\}$ of lattices 
of noncrossing partitions adapted to Motzkin paths $w\in \mathpzc{M}$. We define the 
associated (operator-valued) {\it Motzkin cumulants} and 
prove the corresponding M\"{o}bius inversion formula, which can be viewed as a lattice-level
refinement of the relation expressing free cumulants in terms of Boolean cumulants.
As an application, we obtain an additive decomposition of free cumulants in terms of scalar-valued counterparts of Motzkin cumulants.
\end{abstract}
\maketitle

%\tableofcontents

\section{Introduction}

Cumulants are related to moments of random variables by means of combinatorial recursions described with the use of 
certain partition lattices. In classical probability, these are the lattices of all partitions. The corresponding quantities, 
called {\it classical cumulants}, are related to classical independence by the linearization property of the logarithm of 
the Fourier transform.

In free probability, based on free independence or freeness, an analogous role is played by {\it free cumulants}, whose 
combinatorial definition, due to Speicher \cite{[Sp]}, uses the family 
$\{{\rm NC}(n):n\in \mathbb{N}\}$ of lattices of 
noncrossing partitions introduced by Kreweras \cite{[Kr]}.  
Free cumulants of a single random variable are numbers $(r_n)_{n\geq 1}$ 
defined implicitly in terms of its moments $(m_n)_{n\geq 1}$ by the moment-cumulant formula
\[
m_n=\sum_{\pi\in {\rm NC}(n)}r_{\pi},
\]
where $r_{\pi} =\prod_{V\in \pi}r_{|V|}$
is the product of free cumulants corresponding to blocks $V\in \pi$ and $|V|$ is the number of elements in $V$. 
This relation can be inverted via M\"obius inversion. Moreover, the cumulants $r_n$ are the coefficients 
of the Voiculescu R-transform \cite{[V1], [V2]} which plays a role analogous to that 
of the logarithm of the Fourier transform and is related to free independence by the associated linearization property.

In the multivariate case, the free cumulants $r_n(a_1, \ldots, a_n)$ of elements
$a_{1},\ldots, a_{n}$ of a unital algebra ${\mathcal A}$ equipped with
a normalized linear functional $\varphi$ are defined in terms of their mixed moments 
in a similar way, with the multivariate R-transform developed by Nica \cite{[N1]}.
In turn, multivariate {\it Boolean cumulants} $\beta_n(a_1, \ldots, a_n)$, which linearize the addition
of Boolean independent random variables, are defined with the use of the family 
$\{{\rm Int}(n):n\in \mathbb{N}\}$ of the lattices of interval partitions. We have
\begin{eqnarray*}
\varphi(a_1 \cdots a_n)&=&\sum_{\pi\in {\rm Int}(n)}\beta_{\pi}[a_1, \ldots, a_n],\\
\varphi(a_1 \cdots a_n)&=&\sum_{\pi\in {\rm NC}(n)}r_{\pi}[a_1, \ldots, a_n],
\end{eqnarray*}
where $r_{\pi}$ and $\beta_{\pi}$ factorize according to the block structure of $\pi$.
In particular, if $\hat{1}_{n}$ is the partition consisting of one block $\{1, \ldots, n\}$,
then
$r_{\hat{1}_{n}}[a_1, \ldots, a_n]= r_{n}(a_1, \ldots, a_n)$ and 
$\beta_{\hat{1}_{n}}[a_1, \ldots, a_n]=\beta_n(a_1, \ldots, a_n)$, 
and therefore the above moment-cumulant formulas serve as implicit definitions of the 
sequences $(r_{n})_{n\geq 1}$, $(\beta_{n})_{n\geq 1}$ of multilinear functionals, for which 
we use the same symbols as in the univariate case. See \cite{[N1],[NS],[Sp],[SW]} for further details.

There are combinatorial formulas that relate free cumulants to Boolean cumulants, obtained by 
Lehner \cite{[Leh1]} for one random variable and  
by Belinschi and Nica \cite{[BN]} in the multivariate case. They take the form 
\begin{eqnarray*}
\beta_{n}(a_1, \ldots, a_n)&=&\sum_{\pi\in {\rm NC}_{{\rm irr}}(n)}r_{\pi}[a_1, \ldots , a_n],\\
r_{n}(a_1, \ldots , a_n)&=&\sum_{\pi\in {\rm NC}_{{\rm irr}}(n)}(-1)^{|\pi|-1}\beta_{\pi}[a_1, \ldots, a_n],
\end{eqnarray*}
where the lattices ${\rm NC}_{{\rm irr}}(n)$ consist of irreducible 
noncrossing partitions. Our definition of new cumulants will resemble the first of the above equations.

Clearly, ${\rm Int}(n)\subseteq {\rm NC}(n)$ for all $n$, but there is a growing disparity between the 
cardinalities of ${\rm Int}(n)$ and ${\rm NC}(n)$ when $n$ increases since
\[
|{\rm NC}(n)|-|{\rm Int}(n)|
=\frac{1}{n+1}{2n \choose n} - 2^{n-1}
\rightarrow \infty
\]
as $n\rightarrow \infty$,
which reflects a structural gap between the corresponding families of lattices that, to our knowledge, 
has not been systematically studied. 
Our aim is to bridge this gap by introducing a suitably defined family of intermediate lattices and to identify the
associated cumulants which can be viewed as `partial free cumulants'.
For that purpose, we will use the {\it orthogonal replicas} obtained in the framework of
the unification of noncommutative independence by means of tensor independence \cite{[L3],[L4]}. 

We have recently shown in \cite{[L6]} that the mixed moments of free random 
variables can be decomposed into sums of the mixed moments of their orthogonal replicas.
Namely, if $({\mathcal A}_{1}, \varphi_1)$ and 
$({\mathcal A}_{2}, \varphi_2)$ are noncommutative probability spaces and 
$({\mathcal A}_{1}\star {\mathcal A}_{2}, \varphi_1\star \varphi_2)$ is their free product 
with identified units \cite{[Av], [V1]}, then 
the linear functional $\varphi$ on the free product without identification of units 
${\mathcal A}_{1}*{\mathcal A}_{2}$ given by 
\[
\varphi (a_1\cdots a_n)=\sum_{w=j_1\cdots j_n\in\mathpzc{M}_{n}}\Phi(a_1(j_1)\cdots a_n(j_n))
\]
satisfies the equation $\varphi=(\varphi_1\star \varphi_2)\circ \tau$,
where $\tau$ is the unit identification map and $\Phi$ is a suitably defined normalized linear functional.
Here, $a_1(j_1), \ldots, a_n(j_n)$ are the orthogonal replicas of 
variables $a_1\in \mathcal{A}_{i_1}, \ldots, a_n\in \mathcal{A}_{i_n}$, 
with {\it labels} $i_1, \ldots, i_n\in I=\{1,2\}$, where $i_1\neq \cdots \neq i_n$,
and {\it colors} $j_1, \ldots, j_n\in {\mathbb N}$ encoded by the word $w=j_1 \cdots j_n$. 
Their nonvanishing mixed moments are
associated with words called {\it reduced Motzkin words} which can be
identified with the lattice of {\it Motzkin paths} of length $n$ and height $h(w)=j_1=j_n=1$, 
denoted here by $\mathpzc{M}_{n}$. Their Motzkin subpaths, used later in block nesting, 
are of similar shapes, but they have heights $j\geq 1$. 

This decomposition formula motivates us to introduce the families of lattices of 
noncrossing and interval partitions adapted to Motzkin paths 
$\mathpzc{M}:=\bigcup_{n=1}^{\infty}\mathpzc{M}_n$, namely 
$\{\mathcal{NC}(w):w\in \mathpzc{M}\}$ and
$\{{\mathcal Int}(w):w\in \mathpzc{M}\}$, respectively. The blocks of $\pi\in \mathcal{NC}(w)$ 
are the usual blocks of the associated $\pi_0\in {\rm NC}(n)$, to which we assign appropriate 
Motzkin subpaths.
We will show that the disjoint union of these lattices
\[
\bigsqcup_{w\in \mathpzc{M}_{n}}\mathcal{NC}(w)
\]
with the partial order inherited from ${\rm NC}(n)\times \mathpzc{M}_{n}$ is a poset that 
encodes the decomposition of free cumulants of order $n$ into sums of
`partial free cumulants'. In turn, the lattices ${\mathcal Int}(w)$ correspond to 
the decomposition of mixed moments into sums of products of block Boolean cumulants adapted to $w$.

It is worth emphasizing that Motzkin paths play a constructive role in the combinatorics of freeness. Noncrossing partitions describe admissible block structures, but they do not specify how such structures arise. In contrast, a Motzkin path $w$
imposes blockwise constraints through the subwords assigned to individual blocks, together with compatibility 
conditions between neighboring blocks. This additional structure, encoded in the definition of 
${\mathcal NC}(w)$, leads to a finer combinatorial framework underlying the cumulant expansions considered in this paper.

In order to obtain a refinement of the relation between free and Boolean cumulants, we assume that we are given an 
operator-valued noncommutative probability space $(\mathcal{C}, \mathcal{B}, E)$.
Then, the formulas
\begin{eqnarray*}
E(a_1\cdots a_n)&=&\sum_{\pi\in {\mathcal Int}(w)}B_{\pi}[a_1, \ldots, a_n],\\
E(a_1\cdots a_n)&=&\sum_{\pi\in {\mathcal NC}(w)}K_{\pi}[a_1, \ldots, a_n]
\end{eqnarray*}
where $a_1, \ldots a_n\in \mathcal{C}$, may serve as
implicit definitions of multilinear functions given by
$B_{w}(a_1, \ldots, a_n)=B_{\hat{1}_{w}}[a_1, \ldots, a_n]$ and 
$K_w(a_1, \ldots, a_n)=K_{\hat{1}_{w}}[a_1, \ldots, a_n]$, where 
the symbol $\hat{1}_{w}=(\hat{1}_{n},w)$
stands for the partition consisting of one block $\{1, \ldots, n\}$ to which we 
assign the word $w$.
They will be called the {\it $w$-Boolean cumulants} and {\it Motzkin cumulants}, respectively.

Under some assumptions on the values of the $w$-Boolean cumulants, we obtain the following relations:
\begin{eqnarray*}
B_{w}(a_1, \ldots, a_n)&=&\sum_{\pi\in \mathcal{NC}_{{\rm irr}}(w)}K_{\pi}[a_1, \ldots , a_n],\\ 
K_{w}(a_1, \ldots, a_n)&=&\sum_{\pi\in \mathcal{NC}_{{\rm irr}}(w)}(-1)^{|\pi|-1}B_{\pi}[a_1, \ldots, a_n],
\end{eqnarray*}
where $\mathcal{NC}_{{\rm irr}}(w)$ is the sublattice of $\mathcal{NC}(w)$ consisting of irreducible 
partitions. The last equation is the M\"{o}bius inversion formula
which plays the role of a lattice refinement of the one expressing free cumulants in terms of Boolean cumulants.

We then apply these results to free probability and show that the mixed free
cumulants of free random variables can be decomposed:
\[
r_n(\tau(a_1), \ldots, \tau(a_n))=\sum_{w\in \mathpzc{M}_{n}}k_w(a_{1},\ldots, a_{n}),
\]
where the summands are {\it scalar-valued Motzkin cumulants}
\[
k_w(a_1, \ldots, a_n)=\zeta(K_w(a_{1}(j_1), \ldots , a_{n}(j_n))),
\]
where $a_1\in \mathcal{A}_{i_1}, \ldots, a_n\in \mathcal{A}_{i_n}$,
$w=j_1\cdots j_n\in \mathpzc{AM}$, and $\zeta:{\mathcal B}\rightarrow {\mathbb C}$ is a linear functional.
As before, the mapping $\tau$ is the unit identification map. Let us add that introducing $(\mathcal{C}, \mathcal{B}, E)$ in the 
context of orthogonal replicas is based on the decomposition of the unit 
in the unital algebra generated by these replicas, which leads to $\mathcal{B}$ and 
this enables us to apply the concept of Motzkin cumulants to free probability.

The scalar-valued Motzkin cumulants vanish if not all variables have the same labels, 
thus they satisfy a linearization property. Therefore, in view of the above decomposition of free cumulants,
they can be viewed as `partial free cumulants'. From a combinatorial perspective, they are naturally indexed by Motzkin paths and the associated lattices $\mathcal{NC}(w)$, which provides a refined framework interpolating between noncrossing and interval partitions. 
Since they are labeled by elements of $\mathpzc{M}_{n}$ and 
there exists a bijective correspondence between $\mathpzc{M}_{n}$ and the set ${\mathpzc T}_{n}^{(3)}$ 
of standard Young tableaux with $n$ cells and at most three rows \cite{[EFHH]},
we also obtain a connection between free cumulants and Standard Young Tableaux. This 
connection is realized by this bijection and is not of the fundamental type discovered by Biane \cite{[Bi]} 
and studied by other authors in the context of representations of the symmetric groups 
(see the paper of Feray and \'{S}niady \cite{[FS]} and its references).

Associating Motzkin words with mixed moments has its roots in the tensor product unification of independence 
\cite{[L1]}, where the $j$th color corresponds to the $j$th tensor factor and is closely 
related to the decompositions of the free additive convolution \cite{[L4]} and the free product of graphs \cite{[ALS]}.
In the cumulant context, it is linked to the combinatorics on words used in the construction of 
the ``Boolean-classical logarithm of the Fourier transform'' \cite{[L2]}.
Interestingly enough, monotone independent random variables are also built from orthogonal replicas: it suffices to 
take a replica of one color for one label and replicas of two colors for the other label (see \cite{[L6]} and Lemma 7.15 
in this paper). 
For combinatorial relations between classical, 
free, Boolean and monotone cumulants, see the paper of Arizmendi, Hasebe, Lehner and Vargas \cite{[AHLV]}.
For an operad approach to objects like convolutions and cumulants, see
the paper of Jekel and Liu \cite{[JL]}.

The paper is organized as follows. Section 2 recalls the relevant facts on
free and Boolean cumulants. In Section 3, we introduce the lattices
$\mathcal{NC}(w)$ associated with Motzkin words and study their basic
properties. Section 4 introduces $w$-Boolean cumulants and Motzkin cumulants,
and the corresponding M\"{o}bius inversion formula is proved in Section 5.
Sections 6--9 develop the framework of orthogonal replicas and
projection-valued moments, leading to the main vanishing results for Motzkin
cumulants. In Section 10, we derive the decomposition formula for free
cumulants of free random variables in terms of scalar-valued Motzkin
cumulants. Finally, Section 11 applies this decomposition to the free
convolution of distributions.

\section{Boolean and free cumulants}

Let us recall the definitions of free and Boolean cumulants and the combinatorial relations between them.

A noncommutative probability space is a pair $({\mathcal A},\varphi)$, where 
${\mathcal A}$ is a unital algebra with unit $1_{{\mathcal A}}$ and $\varphi$ is a normalized 
linear functional on ${\mathcal A}$ ($\varphi:{\mathcal A}\rightarrow {\mathbb C}$, $\varphi(1_{{\mathcal A}})=1$).
The {\it distribution} of $a\in \mathcal{A}$ with respect to $\varphi$ is the linear functional
\[
\mu_a:{\mathbb C}\langle X \rangle\rightarrow {\mathbb C}, \;\;\;
\mu_a(f)=\varphi(f(a)) \;\;{\rm for} \;\;f\in {\mathbb C}\langle a \rangle.
\]
where ${\mathbb C}\langle X \rangle$ is the unital algebra of polynomials in the indeterminate $X$.
Typical objects to study are the moments $m_n=\varphi(a^{n})$ and the associated cumulants
defined by means of some lattice of partitions.

More generally, if $\{a_{i}\}_{i\in I}$ is a family of random variables from 
${\mathcal A}$, where $I$ is an index set, then the {\it joint distribution} 
of $\{a_{i}\}_{i\in I}$ with respect to $\varphi$ is the linear functional 
\[
\mu:{\mathbb C}\langle \{X_{i}\}_{i\in I}\rangle \rightarrow {\mathbb C}, \;\;\;
\mu(f)=\varphi(f(\{a_{i}\}_{i\in I}))\;\;{\rm for}\;\; f\in {\mathbb C}\langle \{X_{i}\}_{i\in I}\rangle
\] 
where ${\mathbb C}\langle \{X_i\}_{i\in I}\rangle$ is the unital algebra of noncommutative polynomials
in the family of indeterminates $\{X_{i}\}_{i\in I}$. 
Typical objects to study are the mixed moments
$\varphi(a_{k_1} \cdots a_{k_n})$
and the associated mixed cumulants defined by means of some lattice of partitions.

\medskip

\noindent{\it Partitions.}

\begin{enumerate}
\item
A {\it partition} of the set $[n]:=\{1,2, \ldots, n\}$ is a collection 
$\pi=\{V_1,\ldots , V_k\}$ of non-empty subsets of $[n]$ such that $V_i\cap V_j=\emptyset$
if $i\neq j$ and $[n]=V_1\cup \ldots \cup V_k$. We adopt the convention 
that 
\[
V_1<\cdots <V_k,
\]
by which we understand that the smallest number in $V_j$ is smaller than the smallest number 
in $V_{j+1}$ for any $j=1, \ldots ,k-1$.
\item
We say that $\pi$ is an {\it interval partition} if each block $V_j$ is an interval, i.e. 
$V_j=\{k, k+1, \ldots, l\}$ for some $k\leq l$. The set of interval partitions of $[n]$ 
is a lattice and will be denoted by ${\rm Int}(n)$, however we will not explicitly use its lattice operations.
\item
We say that $\pi$ is {\it noncrossing} if there do 
not exist two distinct blocks $V_i,V_j$ and elements $k,l\in V_i$, $p,q\in V_j$, such that
$k<p<l<q$. The set of noncrossing partitions of $[n]$ is a lattice and is denoted by 
${\rm NC}(n)$. The partial order is given by reversed refinement, 
$\pi\preceq \rho$ if and only if $\pi$ is a (not necessarily proper) refinement of $\rho$.
\item
The {\it join} (least upper bound) of two partitions, $\pi, \rho\in {\rm NC}(n)$, will be denoted by $\pi\curlyvee \rho$.
\item
The partition of $[n]$ consisting of a single block is denoted by $\hat{1}_{n}$ (the greatest element of ${\rm NC}(n)$), whereas that consisting of $n$ singletons is denoted by $\hat{0}_{n}$ (the least element).
\end{enumerate}

\medskip
\noindent{\it Nesting.}

\begin{enumerate}
\item
Let $\pi\in {\rm NC}(n)$. A block $V_i$ is called an {\it inner} block of a block $V_j$ if there exist $p,q\in V_j$
such that $p<k<q$ for all $k\in V_i$. In that case, $V_j$ is called an {\it outer} block of $V_i$.
If $V_j$ does not have any outer blocks, it is called a {\it outermost block}.
If $V_i$ is not a outermost block, then its {\it nearest outer block} is the unique block $V_j$ such that 
$V_j$ is an outer block of $V_i$ and there is no $V_k$ such that $V_k$ is an outer block of $V_i$
and $V_j$ is an outer block of $V_k$. We write in this case 
$V_j=o(V_i)$.
%and we will also say that $V_j$ {\it covers} $V_i$.
\item
We say that a block $V_j$ has {\it depth} $d\in {\mathbb N}$
if there exists a sequence of blocks $(V_{k_1},\ldots , V_{k_{d-1}})$, such that 
\[
V_{k_{1}}=o(V_j),\;V_{k_{2}}=o(V_{k_1}), \ldots , V_{k_{d-1}}=o(V_{k_{d-2}})
\]
and $V_{k_{d-1}}$ is a outermost block. In that case we write $d(V_j)=d$.
We say that a noncrossing partition $\pi$ has depth $d$ if it contains a block of depth $d$ and no block of depth $d+1$. We then 
write $d(\pi)=d$.
\item
A partition $\pi\in {\rm NC}(n)$ is {\it irreducible} if it has exactly one outermost block.
The set of noncrossing irreducible partitions of $[n]$ is denoted by ${\rm NC}_{{\rm irr}}(n)$.
\end{enumerate}

\begin{Definition}
{\rm 
By the {\it free cumulants} $(r_n)_{n\in \mathbb{N}}$  and {\it Boolean cumulants} $(\beta_{n})$ of 
the normalized linear functional $\mu:{\mathbb C}\langle X \rangle \rightarrow {\mathbb C}$ we mean the numbers defined 
in terms of its moments $(m_n)_{n\in \mathbb{N}}$ by the recursive formulas
\begin{eqnarray*}
m_{n}&=&\sum_{\pi\in {\rm NC}(n)}r_{\pi}\\
m_{n}&=&\sum_{\pi\in {\rm Int}(n)}\beta_{\pi},
\end{eqnarray*}
where 
$r_{\pi}=\prod_{V\in \pi}r_{|V|}$
for any $\pi\in {\rm NC}(n)$ and
$\beta_{\pi}=\prod_{V\in\pi}\beta_{|V|}$
for any $\pi\in {\rm Int}(n)$, respectively, with $|V|$ being the cardinality of $V$. 
}
\end{Definition}

Both Boolean and free cumulants can be expressed in terms of moments via M\"{o}bius inversion.
Combinatorial relations between Boolean and free cumulants of one variable were obtained by Lehner \cite{[Leh1]}, who
used the family of {\it irreducible} noncrossing partitions of $[n]$, denoted by ${\rm NC}_{{\rm irr}}(n)$,
which is the subset of ${\rm NC}(n)$ consisting of those partitions in which $1$ and $n$ belong to 
the same block. 
\begin{Theorem} \cite{[Leh1]}
The following relations hold:
\begin{eqnarray*}
\beta_{n}&=&\sum_{\pi\in {\rm NC}_{{\rm irr}}(n)}r_{\pi}\\
r_{n}&=&\sum_{\pi\in {\rm NC}_{{\rm irr}}(n)}(-1)^{|\pi|-1}\beta_{\pi},
\end{eqnarray*}
where $|\pi|$ is the number of blocks in $\pi$.
\end{Theorem}

\begin{Example}
{\rm For illustration, for the lowest-order cumulants, we have
\[
r_1=\beta_1,\;\;
r_2=\beta_2,\;\;
r_3=\beta_3-\beta_2\beta_1,\;\;
r_4=\beta_4-2\beta_3\beta_1-\beta_2^{2}+\beta_2\beta_1^2.
\]
After expressing the Boolean cumulants in terms of moments, we obtain
\begin{eqnarray*}
r_1&=&m_1,\\
r_2&=&m_2-m_1^2,\\
r_3&=&m_3-3m_2m_1+2m_1^3,\\
r_4&=&m_4-4m_3m_1 -2m_2^2 +10m_2m_1^2-5m_1^4,
\end{eqnarray*}
where the coefficients standing by the $m_1^n$ are Catalan numbers multiplied by $(-1)^{n-1}$.
}
\end{Example}

While cumulants of a single variable serve as a reference point, in our framework the multivariate ones are of primary interest, since the variables under consideration are typically constructed from more elementary components.
The multivariate free and Boolean cumulants \cite{[NS],[Sp]} are defined by means of
partitioned multilinear functionals which are similar to
$r_{\pi}, \beta_{\pi}$ defined above. In our notation, we use the same letters as in the univariate case.

\begin{Definition}
{\rm 
Let ${\mathcal A}$ be an algebra and let $(f_n)_{n\in {\mathbb N}}$ be a sequence of 
$n$-linear functio\-nals $f_n:{\mathcal A}\times \cdots \times {\mathcal A}\rightarrow {\mathbb C}$
with values $f_n(a_1, \ldots , a_n)$. We introduce the associated families of multilinear functionals
\begin{eqnarray*}
f_{\pi}[a_1, \ldots , a_n]&:=&\prod_{V\in \pi}f_V[a_1, \ldots , a_n],
\end{eqnarray*}
where
$
f_V[a_1, \ldots, a_n]:=f_k(a_{i_1}, \ldots , a_{i_k})
$
for any set $V=\{i_1< \cdots  <i_k\}\subseteq [n]$ and any $\pi\in {\rm NC}(n)$, where
$a_1, \ldots , a_n\in \mathcal{A}$. 
In particular, if $\pi=\hat{1}_{n}$, then 
$f_{\pi}\equiv f_V\equiv f_{n}$. 
}
\end{Definition}

\begin{Definition}
{\rm 
The {\it multivariate free cumulants} and {\it multivariate Boolean cumulants}
are multilinear functionals $(r_{n})_{n\in {\mathbb N}}$ 
and $(\beta_{n})_{n\in {\mathbb N}}$, respectively, defined by the recursive formulas
\begin{eqnarray*}
\varphi(a_1\cdots a_n)&=&\sum_{\pi\in {\rm NC}(n)}r_{\pi}[a_1, \ldots , a_n]\\
\varphi(a_1\cdots a_n)&=&\;\;\sum_{\pi\in {\rm Int}(n)}\;\beta_{\pi}[a_1, \ldots , a_n]
\end{eqnarray*}
respectively, where $a_1, \ldots , a_n\in \mathcal{A}$ and
$r_{\pi}$ and $\beta_{\pi}$ coincide with $r_n$ and $\beta_n$, respectively, 
when $\pi=\hat{1}_{n}$ (one block). 
}
\end{Definition}

The relations between multivariate free cumulants and multivariate Boolean cumulants are similar to those
in the univariate case, as shown by Belinschi and Nica \cite{[BN]}.

\begin{Theorem}\cite{[BN]}
The following relations hold:
\begin{eqnarray*}
\beta_{n}(a_1, \ldots, a_n)&=&\sum_{\pi\in {\rm NC}_{{\rm irr}}(n)}r_{\pi}[a_1, \ldots , a_n]\\
r_{n}(a_1, \ldots , a_n)&=&\sum_{\pi\in {\rm NC}_{{\rm irr}}(n)}(-1)^{|\pi|-1}\beta_{\pi}[a_1, \ldots, a_n],
\end{eqnarray*}
for any $a_1, \ldots, a_n\in {\mathcal A}$ and any $n\in \mathbb{N}$.
\end{Theorem}

\begin{Example}
{\rm For the lowest-order cumulants, we have relations
\begin{eqnarray*}
r_1(a_1)&=&\beta_1(a_1),\\
r_2(a_1,a_2)&=&\beta_2(a_1,a_2),\\
r_3(a_1,a_2,a_3)&=&\beta_3(a_1,a_2,a_3)-\beta_2(a_1,a_3)\beta_1(a_2),\\
r_4(a_1,a_2,a_3,a_4)&=&\beta_4(a_1,a_2,a_3,a_4)-\beta_3(a_1,a_2,a_4)\beta_1(a_3)\\
&&
\!\!\!\!\!\!-\;\beta_3(a_1,a_3,a_4)\beta_1(a_2)-\beta_2(a_1,a_4)\beta_2(a_2,a_3)\\
&&
\!\!\!\!\!\!+\;\beta_2(a_1,a_4)\beta_1(a_2)\beta_1(a_3),
\end{eqnarray*}
for any $a_1,a_2,a_3,a_4\in \mathcal{A}$, where the corresponding irreducible noncrossing partitions can be easily
identified from the indices of the variables. 
It is not hard to express these free cumulants in terms of mixed moments as in Example 2.3. 
}
\end{Example}

Of course, the leading terms in the formulas for the free cumulants in Theorem 2.6, and thus in Example 2.7,
are Boolean cumulants, but we would like to write the right-hand side of each equation as a sum of cumulants
and not only as a linear combination of products of cumulants.

\section{Lattices $\mathcal{NC}(w)$}

We recall notions conerning Motzkin paths and words (see \cite{[L6]}), and introduce partitions adapted
to them which is needed for our construction.

\medskip
\noindent{\it Motzkin paths and words.}

\begin{enumerate}
\item
A {\it Motzkin path} of length $n$ is a lattice path in the integer plane from 
$(0,1)$ to $(n,1)$ that never passes below the line $y=1$, and is composed exclusively of steps 
of three types: an up step $(1,1)$, a down step $(1,-1)$ and a level step $(1,0)$; 
equivalently, it is a sequence of such steps whose cumulative vertical displacement relative 
to the starting level is always nonnegative and returns to zero at the endpoint.
Our definition, named after Motzkin \cite{[Mo]}, slightly differs from that used by most authors, 
in which the paths start at $(0,0)$, end at $(n,0)$ and never pass below the $x$-axis.
\item
Using the heights of the path, we can assign a word corresponding to each Motzkin path.
By a {\it reduced Motzkin word} of length $n\in \mathbb{N}$ we mean a word over the alphabet 
${\mathbb N}$ of the form
\[
w=j_1j_2\cdots j_n
\]
where $j_1, \ldots, j_n\in \mathbb{N}$, $j_1=j_n=1$ and $j_{i}-j_{i-1}\in \{-1,0,1\}$ for $i=2, \ldots, n$.
The set of reduced Motzkin words of length $n$ will be denoted by $\mathpzc{M}_{n}$. Then
\[
\mathpzc{M}=\bigcup_{n=1}^{\infty}\mathpzc{M}_{n}
\] 
will stand for the set of all nonempty reduced Motzkin words. 
The set of Motzkin paths of length $n-1$ will also be denoted by $\mathpzc{M}_{n}$ via the natural identification with reduced Motzkin words.
\item
It is known that each $\mathpzc{M}_{n}$ is a lattice with the natural 
order obtained by identifying a word $j_1\cdots j_n$ with its height function on $\{0, \ldots, n-1\}$; namely
for $f,g\in \mathpzc{M}_{n}$, we write $f\leq g$ if and only if $f(x)\leq g(x)$ 
for all $x$.
\item
By a {\it Motzkin word} of length $n\in \mathbb{N}$ we shall understand 
a word of the above form, where $j_1=j_n=j$ for some $j\in {\mathbb N}$ 
and such that
\[
w':=j_1'j_2' \cdots j_n'\in \mathpzc{M}
\]
where $j_k'=j_k-j+1$ for any $k$. We will write in this case $h(w)=j$ and say that the {\it height} of $w$ is $j$. 
More generally, the height of any word $v=j_1\ldots j_n$ is $h(v)=\max\{j_1,\ldots, j_n\}$.
In particular, $h(j^{k})=j$ for $k\in \mathbb{N}$ and in that special case we will also refer to $j$ as the {\it color} of $j^k$. 
The set of all Motzkin words will be denoted by $\mathpzc{AM}$
and $\mathpzc{AM}_{n}$ will be its subset consisting of words of length $n$.
\end{enumerate}

\needspace{6\baselineskip}
\medskip
\noindent{\it Partitions decorated by words.}

\begin{enumerate}
\item
Let $\pi_{0}=\{V_{1},\ldots , V_{k}\}\in {\rm NC}(n)$ and let
$w=j_1\cdots j_n\in \mathpzc{AM}$, where $j_1,\ldots,j_n\in \mathbb{N}$.
For each block $V=\{i_1<\cdots<i_p\}\in \pi_0$, consider the subword
\[
v:=j_{i_1}\cdots j_{i_p}.
\]
We restrict attention to those pairs $(\pi_0,w)$ for which each such subword $v$ belongs to $\mathpzc{AM}$. 
In that case, we associate to $(\pi_0,w)$ the collection of pairs
\[
\{(V_1,v_1),\ldots,(V_k,v_k)\}.
\]
\item
For fixed $\pi$, it is often
convenient to identify $\pi$ with the set
$\{v_1,\ldots , v_k\}$ and we will refer to both $(V_j,v_j)$ and $v_j$ as {\it blocks} of $\pi$, with the understanding that 
$v_j$ encodes $V_j$. 
\item
The {\it depth} of $v_j$ is $d(v_j):=d(V_j)$
and $v_j$ is the {\it nearest outer block} of $v_i$, denoted by $v_j=o(v_i)$,
whenever $V_j=o(V_i)$.
\item
If $v_j=i_1\cdots i_p=o(v_i)$, where
$v_i=k_1\cdots k_r$, then there exists $q$ such
that $i_q<k_1<k_r<i_{q+1}$. In this case we define the {\it bridge} over $v_i$ by
$b(v_i):=i_qi_{q+1}$.
\end{enumerate}

We now introduce the notion of a noncrossing partition adapted to $w\in \mathpzc{AM}$.

\begin{Definition}
{\rm
Let $w\in \mathpzc{AM}$ be of the form $w=j_1\cdots j_n$, where $j_1, \ldots , j_n\in {\mathbb N}$.
A {\it noncrossing partition adapted to} $w$ is a pair
\[
\pi=(\pi_0,w)\equiv \{(V_1,v_1), \ldots, (V_k,v_k)\},
\]
where $\pi_{0}=\{V_{1},\ldots , V_{k}\}\in {\rm NC}(n)$ and each $v_j$ is the subword of $w$ corresponding to the block $V_j$ and belongs to $\mathpzc{AM}$, such that:
\begin{enumerate}
\item
for any inner block $v_j$, 
\[
d(v_j)\leq h(v_j) \quad \text{and} \quad h(v_j)\geq h(b(v_j)),
\]
\item
neighboring blocks of the same depth have equal heights.
\end{enumerate}
The set of such partitions is denoted by
$\mathcal{NC}(w)$, and its subset consisting of those for which $\pi_0$ is irreducible is denoted by $\mathcal{NC}_{{\rm irr}}(w)$.
}
\end{Definition}

\begin{Example}
{\rm 
We first consider two examples of noncrossing partitions $\pi$ adapted to the Motzkin words $w$ constructed from the colors 
assigned to the legs of $\pi_0$. 
\\
\unitlength=1mm
\special{em.linewidth 0.5pt}
\linethickness{0.5pt}
\begin{picture}(140.00,20.00)(22.00,103.00)

%%%%%%%%%%%%% P2 %%%%%%%%%%%%%%%%
\put(56.00,117.00){\line(1,0){30.00}}   
\put(59.00,114.50){\line(1,0){24.00}}  
\put(62.00,112.50){\line(1,0){6.00}}
\put(74.00,114.50){\line(1,0){9.00}}
\put(56.00,110.00){\line(0,1){7.00}}
\put(59.00,110.00){\line(0,1){4.50}}
\put(62.00,110.00){\line(0,1){2.50}}
\put(65.00,110.00){\line(0,1){2.50}}
\put(68.00,110.00){\line(0,1){2.50}}
\put(71.00,110.00){\line(0,1){2.50}}
\put(74.00,110.00){\line(0,1){4.50}}
\put(77.00,110.00){\line(0,1){4.50}}
\put(80.00,110.00){\line(0,1){4.50}}
\put(83.00,110.00){\line(0,1){4.50}}
\put(86.00,110.00){\line(0,1){7.00}}

\put(55.00,106.00){$\scriptstyle{1}$}
\put(58.00,106.00){$\scriptstyle{2}$}
\put(61.00,106.00){$\scriptstyle{3}$}
\put(64.00,106.00){$\scriptstyle{4}$}
\put(67.00,106.00){$\scriptstyle{3}$}
\put(70.00,106.00){$\scriptstyle{3}$}
\put(73.00,106.00){$\scriptstyle{3}$}
\put(76.00,106.00){$\scriptstyle{3}$}
\put(79.00,106.00){$\scriptstyle{3}$}
\put(82.00,106.00){$\scriptstyle{2}$}
\put(85.00,106.00){$\scriptstyle{1}$}

%%%%%%%%%%%%% P3 %%%%%%%%%%%%%%%%
\put(106.00,117.00){\line(1,0){21.00}}
\put(130.00,117.00){\line(1,0){6.00}}
\put(109.00,114.50){\line(1,0){12.00}}
\put(106.00,110.00){\line(0,1){7.00}}
\put(109.00,110.00){\line(0,1){4.50}}
\put(112.00,110.00){\line(0,1){4.50}}
\put(115.00,110.00){\line(0,1){2.50}}
\put(118.00,110.00){\line(0,1){4.50}}
\put(121.00,110.00){\line(0,1){4.50}}
\put(124.00,110.00){\line(0,1){7.00}}
\put(127.00,110.00){\line(0,1){7.00}}
\put(130.00,110.00){\line(0,1){7.00}}
\put(133.00,110.00){\line(0,1){7.00}}
\put(136.00,110.00){\line(0,1){7.00}}

\put(105.00,106.00){$\scriptstyle{1}$}
\put(108.00,106.00){$\scriptstyle{2}$}
\put(111.00,106.00){$\scriptstyle{3}$}
\put(114.00,106.00){$\scriptstyle{4}$}
\put(117.00,106.00){$\scriptstyle{3}$}
\put(120.00,106.00){$\scriptstyle{2}$}
\put(123.00,106.00){$\scriptstyle{2}$}
\put(126.00,106.00){$\scriptstyle{1}$}
\put(129.00,106.00){$\scriptstyle{1}$}
\put(132.00,106.00){$\scriptstyle{2}$}
\put(135.00,106.00){$\scriptstyle{1}$}
\end{picture}
\\
We now analyze the first partition in more detail: it is irreducible and consists of subwords
$v_1=1^2$, $v_2=23^32$, $v_3=343$, $v_4=3$,
associated with the blocks of $\pi_0$:
\[
V_1=\{1,11\}, \;V_2=\{2,7,8,9,10\}, \;V_3=\{3,4,5\},\;V_4=\{6\},
\]
respectively. They satisfy condition (1) of Definition 3.1 since 
$d(v_j)=h(v_j)=j$ for each $j$ and $h(b(v_2))=1<2=h(v_2)$ as well as
$h(b(v_3))=h(b(v_4))=3=h(v_3)=h(v_4)$, so all bridges have appropriate heights.
Clearly, condition (2) is satisfied. The second partition is also
adapted to the corresponding word, but it is not irreducible.
}
\end{Example}

\begin{Example}
{\rm 
The partitions shown below are not adapted 
to the corresponding Motzkin words.

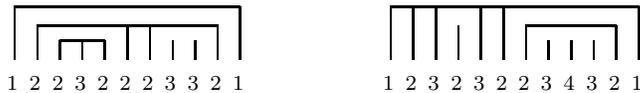
\unitlength=1mm
\special{em.linewidth 0.5pt}
\linethickness{0.5pt}
\begin{picture}(140.00,20.00)(22.00,103.00)

%%%%%%%%%%%%% P2 %%%%%%%%%%%%%%%% 

\put(56.00,117.00){\line(1,0){30.00}}
\put(59.00,114.50){\line(1,0){24.00}}
\put(62.00,112.50){\line(1,0){6.00}}

\put(56.00,110.00){\line(0,1){7.00}}
\put(59.00,110.00){\line(0,1){4.50}}
\put(62.00,110.00){\line(0,1){2.50}}
\put(65.00,110.00){\line(0,1){2.50}}
\put(68.00,110.00){\line(0,1){2.50}}
\put(71.00,110.00){\line(0,1){4.50}}
\put(74.00,110.00){\line(0,1){4.50}}
\put(77.00,110.00){\line(0,1){2.50}}
\put(80.00,110.00){\line(0,1){2.50}}
\put(83.00,110.00){\line(0,1){4.50}}
\put(86.00,110.00){\line(0,1){7.00}}

\put(55.00,106.00){$\scriptstyle{1}$}
\put(58.00,106.00){$\scriptstyle{2}$}
\put(61.00,106.00){$\scriptstyle{2}$}
\put(64.00,106.00){$\scriptstyle{3}$}
\put(67.00,106.00){$\scriptstyle{2}$}
\put(70.00,106.00){$\scriptstyle{2}$}
\put(73.00,106.00){$\scriptstyle{2}$}
\put(76.00,106.00){$\scriptstyle{3}$}
\put(79.00,106.00){$\scriptstyle{3}$}
\put(82.00,106.00){$\scriptstyle{2}$}
\put(85.00,106.00){$\scriptstyle{1}$}

%%%%%%%%%%%%% P3 %%%%%%%%%%%%%%%% 

\put(106.00,117.00){\line(1,0){33.00}}
\put(124.00,114.50){\line(1,0){12.00}}
\put(106.00,110.00){\line(0,1){7.00}}
\put(109.00,110.00){\line(0,1){7.00}}
\put(112.00,110.00){\line(0,1){7.00}}
\put(115.00,110.00){\line(0,1){4.50}}
\put(118.00,110.00){\line(0,1){7.00}}
\put(121.00,110.00){\line(0,1){7.00}}
\put(124.00,110.00){\line(0,1){4.50}}
\put(127.00,110.00){\line(0,1){2.50}}
\put(130.00,110.00){\line(0,1){2.50}}
\put(133.00,110.00){\line(0,1){2.50}}
\put(136.00,110.00){\line(0,1){4.50}}
\put(139.00,110.00){\line(0,1){7.00}}

\put(105.00,106.00){$\scriptstyle{1}$}
\put(108.00,106.00){$\scriptstyle{2}$}
\put(111.00,106.00){$\scriptstyle{3}$}
\put(114.00,106.00){$\scriptstyle{2}$}
\put(117.00,106.00){$\scriptstyle{3}$}
\put(120.00,106.00){$\scriptstyle{2}$}
\put(123.00,106.00){$\scriptstyle{2}$}
\put(126.00,106.00){$\scriptstyle{3}$}
\put(129.00,106.00){$\scriptstyle{4}$}
\put(132.00,106.00){$\scriptstyle{3}$}
\put(135.00,106.00){$\scriptstyle{2}$}
\put(138.00,106.00){$\scriptstyle{1}$}

\end{picture}
\\
\noindent
Indeed, the first partition contains a block $v_3=232$ with $d(v_3)=3>2=h(v_3)$, 
and thus condition (1) of Definition 3.1 does not hold. 
In the second partition, condition (1) fails for the block 
$v_2=2$ since the corresponding bridge 
$b(v_2)=3^2$ has height $h(b(v_2))=3>2=h(v_2)$; moreover, the third block $v_3=2^2$ 
is the nearest outer block of three singletons which are not of the same height, so condition (2) does not hold.
}
\end{Example}

\begin{Definition}
{\rm
For a given $w\in \mathpzc{AM}$, the set of all blocks of 
partitions $\pi\in \mathcal{NC}(w)$ will be denoted by
\[
\mathpzc{B}(w)=\bigcup_{\pi\in \mathcal{NC}(w)}\mathpzc{B}_{\pi}
\]
where $\mathpzc{B}_{\pi}$ is the set of blocks $(V_j,v_j)$ of $\pi$. If $\pi\neq \pi'$, the sets 
$\mathpzc{B}_{\pi}$ and $\mathpzc{B}_{\pi'}$ may not be disjoint.
}
\end{Definition}

\begin{Example}
{\rm
For instance, if $w=12^21$, then
\begin{eqnarray*}
\mathpzc{B}(w)&=&\{(\{1,2,3,4\},12^21), (\{1,2,4\},121), (\{1,3,4\}, 121), \\
&&(\{1,4\}, 1^2), (\{2,3\},2^2), (\{2\},2), (\{3\},2)\}
\end{eqnarray*}
gives all blocks of partitions $\pi\in \mathcal{NC}(w)$.
}
\end{Example}

\medskip
\noindent{\it Order structure.}

We can equip each family $\mathcal{NC}(w)$ with a partial order.
Since ${\rm NC}(n)$ is a lattice with the reversed refinement order -- namely, $\pi_0\preceq \pi_0'$
whenever $\pi_0$ is a refinement of $\pi_0'$ -- it is natural
to use this order to define a partial order on $\mathcal{NC}(w)$. 
More precisely, if $\pi=(\pi_0,w), \pi'=(\pi'_0,w)$ are elements of $\mathcal{NC}(w)$, we set 
\[
\pi\preceq \pi'
\;\;
{\rm whenever}\;\;
\pi_0\preceq \pi_{0}',
\]
where we use the same symbol $\preceq$ for the partial orders on $\mathcal{NC}(w)$ and ${\rm NC}(n)$. Moreover, 
this correspondence is compatible with the lattice structure in the sense that 
\[
\pi\vee \pi'=(\pi_0\vee\pi'_{0},w),
\]
where, again, the same symbol $\vee$ denotes the join operations in both $\mathcal{NC}(w)$ and ${\rm NC}(n)$.
To justify this, one needs to verify that the partition $\pi_0\vee \pi_0'$ is adapted to $w$, which is not immediate.
Consequently, $\mathcal{NC}(w)$ is (lattice)-isomorphic to a sublattice of ${\rm NC}(n)$. 

\begin{Example}
{\rm The cardinality of $\mathcal{NC}(w)$ depends on the word $w$. 
If $w$ is constant (equivalently, corresponds to a constant path), then $\mathcal{NC}(w)$ 
consists only of interval partitions, and in this case $\mathcal{NC}_{{\rm irr}}(w)=\{\hat{1}_{w}\}$. 
The lattice of all noncrossing partitions adapted to
$w=1^22^21\in \mathpzc{M}_{5}$ is shown in Fig.~1, where smaller and larger circles represent
$1$ and $2$, respectively, and the partial order is the reversed refinement order inherited from ${\rm NC}(5)$.
In the sequel, we will also use circles of different sizes (instead of numbers) to 
represent 
colors.
It is easy to see that $\mathcal{NC}(w)$ can be 
obtained from irreducible noncrossing partitions adapted to subwords of $w$. Indeed, 
consider the two decompositions of $w$ into reduced Motzkin words: the trivial one $w=w$ and the non-trivial 
one $w=w_1w_2$, where $w_1=1$ and $w_2=12^21$. This leads to a bijection
\[
\mathcal{NC}(w)\cong\mathcal{NC}_{{\rm irr}}(w)\cup (\{\hat{1}_{w_1}\}\times \mathcal{NC}_{{\rm irr}}(w_2)),
\]
where of course $\mathcal{NC}_{{\rm irr}}(w_1)=\{\hat{1}_{w_{1}}\}$. 
Observe that $\mathcal{NC}_{{\rm irr}}(w)\cong \mathcal{NC}_{{\rm irr}}(12^21)$ which is the diamond lattice. Hence, $\mathcal{NC}(w)$
consists of two diamond lattices. 
This type of decomposition -- at the level of lattices -- extends to arbitrary $w\in \mathpzc{M}$ via irreducible components, see 
Proposition 3.8 below.
A complete list of the families $\mathcal{NC}_{{\rm irr}}(w)$ for all $w\in \mathpzc{M}_5$ is given in Fig.~2.
}
\end{Example}

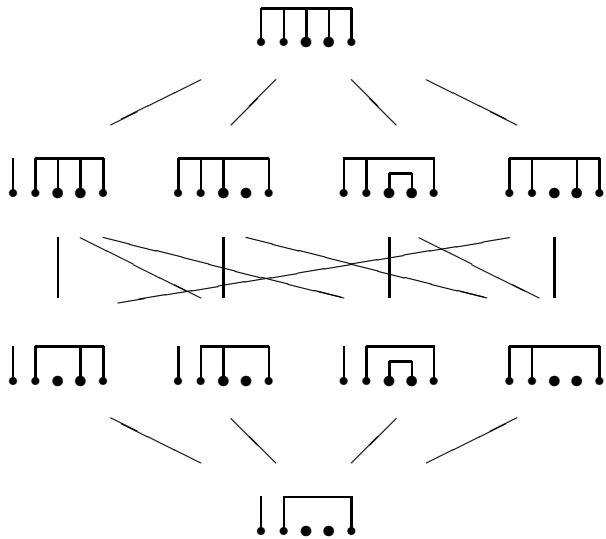
\begin{figure}
\unitlength=1mm
\special{em.linewidth 0.5pt}
\linethickness{0.5pt}
\begin{picture}(140.00,90.00)(-45.00,23.00)

%%%%%%%%%%%%% 1 %%%%%%%%%%%%%%%% 

\put(18.00,30.00){\circle*{1.00}}
\put(21.00,30.00){\circle*{1.00}}
\put(24.00,30.00){\circle*{1.25}}
\put(27.00,30.00){\circle*{1.25}}
\put(30.00,30.00){\circle*{1.00}}
\put(18.00,30.00){\line(0,1){4,5}}
\put(21.00,30.00){\line(0,1){4,5}}
\put(30.00,30.00){\line(0,1){4,5}}
\put(21.00,34.50){\line(1,0){9}}
%%%%%%%%%%%%%%%% links %%%%%%%%%%%%%%%%%%%%%%%

\put(20.00,90.00){\line(-1,-1){6.00}}
\put(30.00,90.00){\line(1,-1){6.00}}
\put(10.00,90.00){\line(-2,-1){12.00}}
\put(40.00,90.00){\line(2,-1){12.00}}

\put(-03.00,69.00){\line(4,-1){32.00}}
\put(-06.00,69.00){\line(2,-1){16.00}}
\put(16.00,69.00){\line(4,-1){32.00}}
\put(39.00,69.00){\line(2,-1){16.00}}
\put(51.00,69.00){\line(-6,-1){52.00}}

\put(-09.00,69.00){\line(0,-1){8.00}}
\put(13.00,69.00){\line(0,-1){8.00}}
\put(35.00,69.00){\line(0,-1){8.00}}
\put(57.00,69.00){\line(0,-1){8.00}}

%%%%%%%%%%%%%%%%% 2 %%%%%%%%%%%%%%%%%

\put(-15.00,50.00){\circle*{1.00}}
\put(-12.00,50.00){\circle*{1.00}}
\put(-09.00,50.00){\circle*{1.25}}
\put(-06.00,50.00){\circle*{1.25}}
\put(-03.00,50.00){\circle*{1.00}}
\put(-15.00,50.00){\line(0,1){4,5}}
\put(-12.00,50.00){\line(0,1){4,5}}
\put(-06.00,50.00){\line(0,1){4,5}}
\put(-03.00,50.00){\line(0,1){4,5}}
\put(-12.00,54.50){\line(1,0){9}}

%%%%%%%%%%%%% 3 %%%%%%%%%%%%%%%% 

\put(07.00,50.00){\circle*{1.00}}
\put(10.00,50.00){\circle*{1.00}}
\put(13.00,50.00){\circle*{1.25}}
\put(16.00,50.00){\circle*{1.25}}
\put(19.00,50.00){\circle*{1.00}}
\put(07.00,50.00){\line(0,1){4,5}}
\put(10.00,50.00){\line(0,1){4,5}}
\put(13.00,50.00){\line(0,1){4,5}}
\put(19.00,50.00){\line(0,1){4,5}}
\put(10.00,54.50){\line(1,0){9}}

%%%%%%%%%%%%% 4 %%%%%%%%%%%%%%%% 

\put(29.00,50.00){\circle*{1.00}}
\put(32.00,50.00){\circle*{1.00}}
\put(35.00,50.00){\circle*{1.25}}
\put(38.00,50.00){\circle*{1.25}}
\put(41.00,50.00){\circle*{1.00}}
\put(29.00,50.00){\line(0,1){4,5}}
\put(32.00,50.00){\line(0,1){4,5}}
\put(35.00,50.00){\line(0,1){2,5}}
\put(38.00,50.00){\line(0,1){2,5}}
\put(41.00,50.00){\line(0,1){4,5}}

\put(35.00,52.50){\line(1,0){3}}
\put(32.00,54.50){\line(1,0){9}}

%%%%%%%%%%%%%%%% 5 %%%%%%%%%%%%%%%%%%%%

\put(51.00,50.00){\circle*{1.00}}
\put(54.00,50.00){\circle*{1.00}}
\put(57.00,50.00){\circle*{1.25}}
\put(60.00,50.00){\circle*{1.25}}
\put(63.00,50.00){\circle*{1.00}}
\put(51.00,50.00){\line(0,1){4,5}}
\put(54.00,50.00){\line(0,1){4,5}}
\put(63.00,50.00){\line(0,1){4,5}}
\put(51.00,54.50){\line(1,0){12}}

%%%%%%%%%%%%% 6 %%%%%%%%%%%%%%%% 

\put(-15.00,75.00){\circle*{1.00}}
\put(-12.00,75.00){\circle*{1.00}}
\put(-09.00,75.00){\circle*{1.25}}
\put(-06.00,75.00){\circle*{1.25}}
\put(-03.00,75.00){\circle*{1.00}}
\put(-15.00,75.00){\line(0,1){4,5}}
\put(-12.00,75.00){\line(0,1){4,5}}
\put(-09.00,75.00){\line(0,1){4,5}}
\put(-06.00,75.00){\line(0,1){4,5}}
\put(-03.00,75.00){\line(0,1){4,5}}
\put(-12.00,79.50){\line(1,0){9}}

%%%%%%%%%%%%% 7 %%%%%%%%%%%%%%%% 

\put(07.00,75.00){\circle*{1.00}}
\put(10.00,75.00){\circle*{1.00}}
\put(13.00,75.00){\circle*{1.25}}
\put(16.00,75.00){\circle*{1.25}}
\put(19.00,75.00){\circle*{1.00}}
\put(07.00,75.00){\line(0,1){4,5}}
\put(10.00,75.00){\line(0,1){4,5}}
\put(13.00,75.00){\line(0,1){4,5}}
\put(19.00,75.00){\line(0,1){4,5}}
\put(07.00,79.50){\line(1,0){12}}

%%%%%%%%%%%% links %%%%%%%%%%%%%%%%%%%%

\put(10.00,39.00){\line(-2,1){12.00}}
\put(20.00,39.00){\line(-1,1){6.00}}
\put(30.00,39.00){\line(1,1){6.00}}
\put(40.00,39.00){\line(2,1){12.00}}

%%%%%%%%%%%%% 8 %%%%%%%%%%%%%%%% 

\put(29.00,75.00){\circle*{1.00}}
\put(32.00,75.00){\circle*{1.00}}
\put(35.00,75.00){\circle*{1.25}}
\put(38.00,75.00){\circle*{1.25}}
\put(41.00,75.00){\circle*{1.00}}
\put(29.00,75.00){\line(0,1){4,5}}
\put(32.00,75.00){\line(0,1){4,5}}
\put(35.00,75.00){\line(0,1){2,5}}
\put(38.00,75.00){\line(0,1){2,5}}
\put(41.00,75.00){\line(0,1){4,5}}
\put(35.00,77.50){\line(1,0){3}}
\put(29.00,79.50){\line(1,0){12}}

%%%%%%%%%%%%%%%% 9 %%%%%%%%%%%%%%%%%%%%

\put(51.00,75.00){\circle*{1.00}}
\put(54.00,75.00){\circle*{1.00}}
\put(57.00,75.00){\circle*{1.25}}
\put(60.00,75.00){\circle*{1.25}}
\put(63.00,75.00){\circle*{1.00}}
\put(51.00,75.00){\line(0,1){4,5}}
\put(54.00,75.00){\line(0,1){4,5}}
\put(60.00,75.00){\line(0,1){4,5}}
\put(63.00,75.00){\line(0,1){4,5}}
\put(51.00,79.50){\line(1,0){12}}

%%%%%%%%%%%%% 10 %%%%%%%%%%%%%%%% 
\put(18.00,95.00){\circle*{1.00}}
\put(21.00,95.00){\circle*{1.00}}
\put(24.00,95.00){\circle*{1.25}}
\put(27.00,95.00){\circle*{1.25}}
\put(30.00,95.00){\circle*{1.00}}
\put(18.00,95.00){\line(0,1){4,5}}
\put(21.00,95.00){\line(0,1){4,5}}
\put(24.00,95.00){\line(0,1){4,5}}
\put(27.00,95.00){\line(0,1){4,5}}
\put(30.00,95.00){\line(0,1){4,5}}
\put(18.00,99.50){\line(1,0){12}}

\end{picture}
\caption{The Hasse diagram of the lattice $\mathcal{NC}(1^22^21)$.}
\end{figure}

Before proving that each $\mathcal{NC}(w)$ is a lattice, 
we describe their collective structure for words $w$ of fixed length, in the category of posets.

\begin{Definition}
{\rm 
Let ${\rm NC}(n)\times \mathpzc{M}_{n}$ be the Cartesian product of posets, equipped with the product partial order
\[
(\pi_0,w)\preceq (\pi'_0,w')\quad \text{whenever}\quad \pi_0\preceq \pi_0'\;\;{\rm and}\;\;w\preceq w',
\]
where $\pi_0,\pi_0' \in {\rm NC}(n)$ and $w,w'\in \mathpzc{M}_n$.
Define
\[
\mathcal{NC}_{n}:= 
\bigsqcup_{w\in \mathpzc{M}_{n}}\mathcal{NC}(w)
\]
to be the induced subposet of ${\rm NC}(n)\times \mathpzc{M}_{n}$ and let 
$\mathcal{NC}_{n, {\rm irr}}$ be the subposet of $\mathcal{NC}_{n}$ consisting of irreducible partitions.
}
\end{Definition} 

\begin{figure}
\unitlength=1mm
\special{em:linewidth 1pt}
\linethickness{0.5pt}
\begin{picture}(180.00,120.00)(-15.00,-30.00)

\put(-5.00, 82.00){$w$}
\put(20.00, 82.00){$\mathcal{NC}_{{\rm irr}}(w)$}

%%%%%%%%%%%%%%%path%%%%%%%%%%%%%%%%%%%%
%\put(-5.00,70.00){$1^5$}

\put(-5.00,70.00){\line(1,0){4.00}}
\put(-1.00,70.00){\line(1,0){4.00}}
\put(3.00,70.00){\line(1,0){4.00}}
\put(7.00,70.00){\line(1,0){4.00}}
\put(-5.00,70.00){\circle*{1.00}}
\put(-1.00,70.00){\circle*{1.00}}
\put(3.00,70.00){\circle*{1.00}}
\put(7.00,70.00){\circle*{1.00}}
\put(11.00,70.00){\circle*{1.00}}

\put(20.00,70.00){\circle*{1.00}}
\put(23.00,70.00){\circle*{1.00}}
\put(26.00,70.00){\circle*{1.00}}
\put(29.00,70.00){\circle*{1.00}}
\put(32.00,70.00){\circle*{1.00}}

\put(20.00,70.00){\line(0,1){4,5}}
\put(23.00,70.00){\line(0,1){4,5}}
\put(26.00,70.00){\line(0,1){4,5}}
\put(29.00,70.00){\line(0,1){4,5}}
\put(32.00,70.00){\line(0,1){4,5}}
\put(20.00,74.50){\line(1,0){12}}

%%%%%%%%%%%%%path%%%%%%%%%%%%%%%%%%%%%%%
%\put(-5.00,60.00){$1^321$}

\put(-5.00,60.00){\line(1,0){4.00}}
\put(-1.00,60.00){\line(1,0){4.00}}
\put(3.00,60.00){\line(1,1){4.00}}
\put(7.00,64.00){\line(1,-1){4.00}}
\put(-5.00,60.00){\circle*{1.00}}
\put(-1.00,60.00){\circle*{1.00}}
\put(3.00,60.00){\circle*{1.00}}
\put(7.00,64.00){\circle*{1.00}}
\put(11.00,60.00){\circle*{1.00}}

\put(20.00,60.00){\circle*{1.00}}
\put(23.00,60.00){\circle*{1.00}}
\put(26.00,60.00){\circle*{1.00}}
\put(29.00,60.00){\circle*{1.25}}
\put(32.00,60.00){\circle*{1.00}}
\put(20.00,60.00){\line(0,1){4,5}}
\put(23.00,60.00){\line(0,1){4,5}}
\put(26.00,60.00){\line(0,1){4,5}}
\put(29.00,60.00){\line(0,1){4,5}}
\put(32.00,60.00){\line(0,1){4,5}}
\put(20.00,64.50){\line(1,0){12}}
\put(38.00,60.00){\circle*{1.00}}
\put(41.00,60.00){\circle*{1.00}}
\put(44.00,60.00){\circle*{1.00}}
\put(47.00,60.00){\circle*{1.25}}
\put(50.00,60.00){\circle*{1.00}}
\put(38.00,60.00){\line(0,1){4,5}}
\put(41.00,60.00){\line(0,1){4,5}}
\put(44.00,60.00){\line(0,1){4,5}}
\put(50.00,60.00){\line(0,1){4,5}}
\put(38.00,64.50){\line(1,0){12}}

%%%%%%%%%%%%%%%path %%%%%%%%%%%%%%%%
%\put(-5.00,50.00){$1^221^2$}

\put(-5.00,50.00){\line(1,0){4.00}}
\put(-1.00,50.00){\line(1,1){4.00}}
\put(3.00,54.00){\line(1,-1){4.00}}
\put(7.00,50.00){\line(1,0){4.00}}
\put(-5.00,50.00){\circle*{1.00}}
\put(-1.00,50.00){\circle*{1.00}}
\put(3.00,54.00){\circle*{1.00}}
\put(7.00,50.00){\circle*{1.00}}
\put(11.00,50.00){\circle*{1.00}}

\put(20.00,50.00){\circle*{1.00}}
\put(23.00,50.00){\circle*{1.00}}
\put(26.00,50.00){\circle*{1.25}}
\put(29.00,50.00){\circle*{1.00}}
\put(32.00,50.00){\circle*{1.00}}
\put(20.00,50.00){\line(0,1){4,5}}
\put(23.00,50.00){\line(0,1){4,5}}
\put(26.00,50.00){\line(0,1){4,5}}
\put(29.00,50.00){\line(0,1){4,5}}
\put(32.00,50.00){\line(0,1){4,5}}
\put(20.00,54.50){\line(1,0){12}}
\put(38.00,50.00){\circle*{1.00}}
\put(41.00,50.00){\circle*{1.00}}
\put(44.00,50.00){\circle*{1.25}}
\put(47.00,50.00){\circle*{1.00}}
\put(50.00,50.00){\circle*{1.00}}
\put(38.00,50.00){\line(0,1){4,5}}
\put(41.00,50.00){\line(0,1){4,5}}
\put(47.00,50.00){\line(0,1){4,5}}
\put(50.00,50.00){\line(0,1){4,5}}
\put(38.00,54.50){\line(1,0){12}}

%%%%%%%%%%%%%%%%%%%%%%%%%%%%%%%%%%%

%\put(-5.00,40.00){$121^3$}

\put(-5.00,40.00){\line(1,1){4.00}}
\put(-1.00,44.00){\line(1,-1){4.00}}
\put(3.00,40.00){\line(1,0){4.00}}
\put(7.00,40.00){\line(1,0){4.00}}
\put(-5.00,40.00){\circle*{1.00}}
\put(-1.00,44.00){\circle*{1.00}}
\put(3.00,40.00){\circle*{1.00}}
\put(7.00,40.00){\circle*{1.00}}
\put(11.00,40.00){\circle*{1.00}}

\put(20.00,40.00){\circle*{1.00}}
\put(23.00,40.00){\circle*{1.25}}
\put(26.00,40.00){\circle*{1.00}}
\put(29.00,40.00){\circle*{1.00}}
\put(32.00,40.00){\circle*{1.00}}
\put(20.00,40.00){\line(0,1){4,5}}
\put(23.00,40.00){\line(0,1){4,5}}
\put(26.00,40.00){\line(0,1){4,5}}
\put(29.00,40.00){\line(0,1){4,5}}
\put(32.00,40.00){\line(0,1){4,5}}
\put(20.00,44.50){\line(1,0){12}}
\put(38.00,40.00){\circle*{1.00}}
\put(41.00,40.00){\circle*{1.25}}
\put(44.00,40.00){\circle*{1.00}}
\put(47.00,40.00){\circle*{1.00}}
\put(50.00,40.00){\circle*{1.00}}
\put(38.00,40.00){\line(0,1){4,5}}
\put(44.00,40.00){\line(0,1){4,5}}
\put(47.00,40.00){\line(0,1){4,5}}
\put(50.00,40.00){\line(0,1){4,5}}
\put(38.00,44.50){\line(1,0){12}}

%%%%%%%%%%%%%%%%%path %%%%%%%%%%%%%%%%%%%
%\put(-5.00,30.00){$12121$}

\put(-5.00,30.00){\line(1,1){4.00}}
\put(-1.00,34.00){\line(1,-1){4.00}}
\put(3.00,30.00){\line(1,1){4.00}}
\put(7.00,34.00){\line(1,-1){4.00}}
\put(-5.00,30.00){\circle*{1.00}}
\put(-1.00,34.00){\circle*{1.00}}
\put(3.00,30.00){\circle*{1.00}}
\put(7.00,34.00){\circle*{1.00}}
\put(11.00,30.00){\circle*{1.00}}

\put(20.00,30.00){\circle*{1.00}}
\put(23.00,30.00){\circle*{1.25}}
\put(26.00,30.00){\circle*{1.00}}
\put(29.00,30.00){\circle*{1.25}}
\put(32.00,30.00){\circle*{1.00}}
\put(20.00,30.00){\line(0,1){4,5}}
\put(23.00,30.00){\line(0,1){4,5}}
\put(26.00,30.00){\line(0,1){4,5}}
\put(29.00,30.00){\line(0,1){4,5}}
\put(32.00,30.00){\line(0,1){4,5}}
\put(20.00,34.50){\line(1,0){12}}
\put(38.00,30.00){\circle*{1.00}}
\put(41.00,30.00){\circle*{1.25}}
\put(44.00,30.00){\circle*{1.00}}
\put(47.00,30.00){\circle*{1.25}}
\put(50.00,30.00){\circle*{1.00}}
\put(38.00,30.00){\line(0,1){4,5}}
\put(41.00,30.00){\line(0,1){4,5}}
\put(44.00,30.00){\line(0,1){4,5}}
\put(50.00,30.00){\line(0,1){4,5}}
\put(38.00,34.50){\line(1,0){12}}
\put(56.00,30.00){\circle*{1.00}}
\put(59.00,30.00){\circle*{1.25}}
\put(62.00,30.00){\circle*{1.00}}
\put(65.00,30.00){\circle*{1.25}}
\put(68.00,30.00){\circle*{1.00}}
\put(56.00,30.00){\line(0,1){4,5}}
\put(62.00,30.00){\line(0,1){4,5}}
\put(65.00,30.00){\line(0,1){4,5}}
\put(68.00,30.00){\line(0,1){4,5}}
\put(56.00,34.50){\line(1,0){12}}
\put(74.00,30.00){\circle*{1.00}}
\put(77.00,30.00){\circle*{1.25}}
\put(80.00,30.00){\circle*{1.00}}
\put(83.00,30.00){\circle*{1.25}}
\put(86.00,30.00){\circle*{1.00}}
\put(74.00,30.00){\line(0,1){4,5}}
\put(80.00,30.00){\line(0,1){4,5}}
\put(86.00,30.00){\line(0,1){4,5}}
\put(74.00,34.50){\line(1,0){12}}

%%%%%%%%%%%%%%%path %%%%%%%%%%%%%%%

%\put(-5.00,20.00){$1^22^21$}

\put(-5.00,20.00){\line(1,0){4.00}}
\put(-1.00,20.00){\line(1,1){4.00}}
\put(3.00,24.00){\line(1,0){4.00}}
\put(7.00,24.00){\line(1,-1){4.00}}
\put(-5.00,20.00){\circle*{1.00}}
\put(-1.00,20.00){\circle*{1.00}}
\put(3.00,24.00){\circle*{1.00}}
\put(7.00,24.00){\circle*{1.00}}
\put(11.00,20.00){\circle*{1.00}}

\put(20.00,20.00){\circle*{1.00}}
\put(23.00,20.00){\circle*{1.00}}
\put(26.00,20.00){\circle*{1.25}}
\put(29.00,20.00){\circle*{1.25}}
\put(32.00,20.00){\circle*{1.00}}
\put(20.00,20.00){\line(0,1){4,5}}
\put(23.00,20.00){\line(0,1){4,5}}
\put(26.00,20.00){\line(0,1){4,5}}
\put(29.00,20.00){\line(0,1){4,5}}
\put(32.00,20.00){\line(0,1){4,5}}
\put(20.00,24.50){\line(1,0){12}}
\put(38.00,20.00){\circle*{1.00}}
\put(41.00,20.00){\circle*{1.00}}
\put(44.00,20.00){\circle*{1.25}}
\put(47.00,20.00){\circle*{1.25}}
\put(50.00,20.00){\circle*{1.00}}
\put(38.00,20.00){\line(0,1){4,5}}
\put(41.00,20.00){\line(0,1){4,5}}
\put(47.00,20.00){\line(0,1){4,5}}
\put(50.00,20.00){\line(0,1){4,5}}
\put(38.00,24.50){\line(1,0){12}}
\put(56.00,20.00){\circle*{1.00}}
\put(59.00,20.00){\circle*{1.00}}
\put(62.00,20.00){\circle*{1.25}}
\put(65.00,20.00){\circle*{1.25}}
\put(68.00,20.00){\circle*{1.00}}
\put(56.00,20.00){\line(0,1){4,5}}
\put(59.00,20.00){\line(0,1){4,5}}
\put(62.00,20.00){\line(0,1){4,5}}
\put(68.00,20.00){\line(0,1){4,5}}
\put(56.00,24.50){\line(1,0){12}}
\put(74.00,20.00){\circle*{1.00}}
\put(77.00,20.00){\circle*{1.00}}
\put(80.00,20.00){\circle*{1.25}}
\put(83.00,20.00){\circle*{1.25}}
\put(86.00,20.00){\circle*{1.00}}
\put(74.00,20.00){\line(0,1){4,5}}
\put(77.00,20.00){\line(0,1){4,5}}
\put(80.00,20.00){\line(0,1){2,5}}
\put(83.00,20.00){\line(0,1){2,5}}
\put(86.00,20.00){\line(0,1){4,5}}
\put(74.00,24.50){\line(1,0){12}}
\put(80.00,22.50){\line(1,0){3}}
\put(92.00,20.00){\circle*{1.00}}
\put(95.00,20.00){\circle*{1.00}}
\put(98.00,20.00){\circle*{1.25}}
\put(101.00,20.00){\circle*{1.25}}
\put(104.00,20.00){\circle*{1.00}}
\put(92.00,20.00){\line(0,1){4,5}}
\put(95.00,20.00){\line(0,1){4,5}}
\put(104.00,20.00){\line(0,1){4,5}}
\put(92.00,24.50){\line(1,0){12}}

%%%%%%%%%%%%%%path %%%%%%%%%%%%%%%%%%%%%
%\put(-5.00,10.00){$12^21^2$}

\put(-5.00,10.00){\line(1,1){4.00}}
\put(-1.00,14.00){\line(1,0){4.00}}
\put(3.00,14.00){\line(1,-1){4.00}}
\put(7.00,10.00){\line(1,0){4.00}}
\put(-5.00,10.00){\circle*{1.00}}
\put(-1.00,14.00){\circle*{1.00}}
\put(3.00,14.00){\circle*{1.00}}
\put(7.00,10.00){\circle*{1.00}}
\put(11.00,10.00){\circle*{1.00}}

\put(20.00,10.00){\circle*{1.00}}
\put(23.00,10.00){\circle*{1.25}}
\put(26.00,10.00){\circle*{1.25}}
\put(29.00,10.00){\circle*{1.00}}
\put(32.00,10.00){\circle*{1.00}}
\put(20.00,10.00){\line(0,1){4,5}}
\put(23.00,10.00){\line(0,1){4,5}}
\put(26.00,10.00){\line(0,1){4,5}}
\put(29.00,10.00){\line(0,1){4,5}}
\put(32.00,10.00){\line(0,1){4,5}}
\put(20.00,14.50){\line(1,0){12}}
\put(38.00,10.00){\circle*{1.00}}
\put(41.00,10.00){\circle*{1.25}}
\put(44.00,10.00){\circle*{1.25}}
\put(47.00,10.00){\circle*{1.00}}
\put(50.00,10.00){\circle*{1.00}}
\put(38.00,10.00){\line(0,1){4,5}}
\put(41.00,10.00){\line(0,1){4,5}}
\put(47.00,10.00){\line(0,1){4,5}}
\put(50.00,10.00){\line(0,1){4,5}}
\put(38.00,14.50){\line(1,0){12}}
\put(56.00,10.00){\circle*{1.00}}
\put(59.00,10.00){\circle*{1.25}}
\put(62.00,10.00){\circle*{1.25}}
\put(65.00,10.00){\circle*{1.00}}
\put(68.00,10.00){\circle*{1.00}}
\put(56.00,10.00){\line(0,1){4,5}}
\put(62.00,10.00){\line(0,1){4,5}}
\put(65.00,10.00){\line(0,1){4,5}}
\put(68.00,10.00){\line(0,1){4,5}}
\put(56.00,14.50){\line(1,0){12}}
\put(74.00,10.00){\circle*{1.00}}
\put(77.00,10.00){\circle*{1.25}}
\put(80.00,10.00){\circle*{1.25}}
\put(83.00,10.00){\circle*{1.00}}
\put(86.00,10.00){\circle*{1.00}}
\put(74.00,10.00){\line(0,1){4,5}}
\put(77.00,10.00){\line(0,1){2,5}}
\put(80.00,10.00){\line(0,1){2,5}}
\put(83.00,10.00){\line(0,1){4,5}}
\put(86.00,10.00){\line(0,1){4,5}}
\put(74.00,14.50){\line(1,0){12}}
\put(77.00,12.50){\line(1,0){3}}
\put(92.00,10.00){\circle*{1.00}}
\put(95.00,10.00){\circle*{1.25}}
\put(98.00,10.00){\circle*{1.25}}
\put(101.00,10.00){\circle*{1.00}}
\put(104.00,10.00){\circle*{1.00}}
\put(92.00,10.00){\line(0,1){4,5}}
\put(101.00,10.00){\line(0,1){4,5}}
\put(104.00,10.00){\line(0,1){4,5}}
\put(92.00,14.50){\line(1,0){12}}

%%%%%%%%%%%%%%%%%%%%%%%%%%%%%%%%%%%%%
\put(-5.00,0.00){\line(1,1){4.00}}
\put(-1.00,4.00){\line(1,0){4.00}}
\put(3.00,4.00){\line(1,0){4.00}}
\put(7.00,4.00){\line(1,-1){4.00}}
\put(-5.00,0.00){\circle*{1.00}}
\put(-1.00,4.00){\circle*{1.00}}
\put(3.00,4.00){\circle*{1.00}}
\put(7.00,4.00){\circle*{1.00}}
\put(11.00,0.00){\circle*{1.00}}

%\put(-5.00,0.00){$12^31$}
\put(20.00,0.00){\circle*{1.00}}
\put(23.00,0.00){\circle*{1.25}}
\put(26.00,0.00){\circle*{1.25}}
\put(29.00,0.00){\circle*{1.25}}
\put(32.00,0.00){\circle*{1.00}}
\put(20.00,0.00){\line(0,1){4,5}}
\put(23.00,0.00){\line(0,1){4,5}}
\put(26.00,0.00){\line(0,1){4,5}}
\put(29.00,0.00){\line(0,1){4,5}}
\put(32.00,0.00){\line(0,1){4,5}}
\put(20.00,04.50){\line(1,0){12}}

\put(38.00,00.00){\circle*{1.00}}
\put(41.00,00.00){\circle*{1.25}}
\put(44.00,00.00){\circle*{1.25}}
\put(47.00,00.00){\circle*{1.25}}
\put(50.00,00.00){\circle*{1.00}}
\put(38.00,00.00){\line(0,1){4,5}}
\put(41.00,00.00){\line(0,1){4,5}}
\put(44.00,00.00){\line(0,1){4,5}}
\put(50.00,00.00){\line(0,1){4,5}}
\put(38.00,04.50){\line(1,0){12}}

\put(56.00,00.00){\circle*{1.00}}
\put(59.00,00.00){\circle*{1.25}}
\put(62.00,00.00){\circle*{1.25}}
\put(65.00,00.00){\circle*{1.25}}
\put(68.00,00.00){\circle*{1.00}}
\put(56.00,00.00){\line(0,1){4,5}}
\put(62.00,00.00){\line(0,1){4,5}}
\put(65.00,00.00){\line(0,1){4,5}}
\put(68.00,00.00){\line(0,1){4,5}}
\put(56.00,04.50){\line(1,0){12}}

\put(74.00,00.00){\circle*{1.00}}
\put(77.00,00.00){\circle*{1.25}}
\put(80.00,00.00){\circle*{1.25}}
\put(83.00,00.00){\circle*{1.25}}
\put(86.00,00.00){\circle*{1.00}}
\put(74.00,00.00){\line(0,1){4,5}}
\put(77.00,00.00){\line(0,1){4,5}}
\put(83.00,00.00){\line(0,1){4,5}}
\put(86.00,00.00){\line(0,1){4,5}}
\put(74.00,04.50){\line(1,0){12}}

\put(92.00,00.00){\circle*{1.00}}
\put(95.00,00.00){\circle*{1.25}}
\put(98.00,00.00){\circle*{1.25}}
\put(101.00,00.00){\circle*{1.25}}
\put(104.00,00.00){\circle*{1.00}}
\put(92.00,00.00){\line(0,1){4,5}}
\put(95.00,00.00){\line(0,1){4,5}}
\put(98.00,00.00){\line(0,1){2,5}}
\put(101.00,00.00){\line(0,1){2,5}}
\put(104.00,00.00){\line(0,1){4,5}}
\put(92.00,04.50){\line(1,0){12}}
\put(98.00,02.50){\line(1,0){3}}

\put(110.00,00.00){\circle*{1.00}}
\put(113.00,00.00){\circle*{1.25}}
\put(116.00,00.00){\circle*{1.25}}
\put(119.00,00.00){\circle*{1.25}}
\put(122.00,00.00){\circle*{1.00}}
\put(110.00,00.00){\line(0,1){4,5}}
\put(113.00,00.00){\line(0,1){2,5}}
\put(116.00,00.00){\line(0,1){2,5}}
\put(119.00,00.00){\line(0,1){4,5}}
\put(122.00,00.00){\line(0,1){4,5}}
\put(110.00,04.50){\line(1,0){12}}
\put(113.00,02.50){\line(1,0){3}}

\put(20.00,-10.00){\circle*{1.00}}
\put(23.00,-10.00){\circle*{1.25}}
\put(26.00,-10.00){\circle*{1.25}}
\put(29.00,-10.00){\circle*{1.25}}
\put(32.00,-10.00){\circle*{1.00}}
\put(20.00,-10.00){\line(0,1){4,5}}
\put(23.00,-10.00){\line(0,1){4,5}}
\put(32.00,-10.00){\line(0,1){4,5}}
\put(20.00,-5.50){\line(1,0){12}}

\put(38.00,-10.00){\circle*{1.00}}
\put(41.00,-10.00){\circle*{1.25}}
\put(44.00,-10.00){\circle*{1.25}}
\put(47.00,-10.00){\circle*{1.25}}
\put(50.00,-10.00){\circle*{1.00}}
\put(38.00,-10.00){\line(0,1){4,5}}
\put(44.00,-10.00){\line(0,1){4,5}}
\put(50.00,-10.00){\line(0,1){4,5}}
\put(38.00,-5.50){\line(1,0){12}}

\put(56.00,-10.00){\circle*{1.00}}
\put(59.00,-10.00){\circle*{1.25}}
\put(62.00,-10.00){\circle*{1.25}}
\put(65.00,-10.00){\circle*{1.25}}
\put(68.00,-10.00){\circle*{1.00}}
\put(56.00,-10.00){\line(0,1){4,5}}
\put(65.00,-10.00){\line(0,1){4,5}}
\put(68.00,-10.00){\line(0,1){4,5}}
\put(56.00,-5.50){\line(1,0){12}}

\put(74.00,-10.00){\circle*{1.00}}
\put(77.00,-10.00){\circle*{1.25}}
\put(80.00,-10.00){\circle*{1.25}}
\put(83.00,-10.00){\circle*{1.25}}
\put(86.00,-10.00){\circle*{1.00}}
\put(74.00,-10.00){\line(0,1){4,5}}
\put(77.00,-10.00){\line(0,1){2,5}}
\put(80.00,-10.00){\line(0,1){2,5}}
\put(83.00,-10.00){\line(0,1){2,5}}
\put(86.00,-10.00){\line(0,1){4,5}}
\put(74.00,-5.50){\line(1,0){12}}
\put(77.00,-7.50){\line(1,0){6}}

\put(92.00,-10.00){\circle*{1.00}}
\put(95.00,-10.00){\circle*{1.25}}
\put(98.00,-10.00){\circle*{1.25}}
\put(101.00,-10.00){\circle*{1.25}}
\put(104.00,-10.00){\circle*{1.00}}
\put(92.00,-10.00){\line(0,1){4,5}}
\put(95.00,-10.00){\line(0,1){2,5}}
\put(98.00,-10.00){\line(0,1){2,5}}
\put(104.00,-10.00){\line(0,1){4,5}}
\put(92.00,-5.50){\line(1,0){12}}
\put(95.00,-7.50){\line(1,0){3}}

\put(110.00,-10.00){\circle*{1.00}}
\put(113.00,-10.00){\circle*{1.25}}
\put(116.00,-10.00){\circle*{1.25}}
\put(119.00,-10.00){\circle*{1.25}}
\put(122.00,-10.00){\circle*{1.00}}
\put(110.00,-10.00){\line(0,1){4,5}}
\put(116.00,-10.00){\line(0,1){2,5}}
\put(119.00,-10.00){\line(0,1){2,5}}
\put(122.00,-10.00){\line(0,1){4,5}}
\put(110.00,-5.50){\line(1,0){12}}
\put(116.00,-7.50){\line(1,0){3}}

\put(128.00,-10.00){\circle*{1.00}}
\put(131.00,-10.00){\circle*{1.25}}
\put(134.00,-10.00){\circle*{1.25}}
\put(137.00,-10.00){\circle*{1.25}}
\put(140.00,-10.00){\circle*{1.00}}
\put(128.00,-10.00){\line(0,1){4,5}}
\put(140.00,-10.00){\line(0,1){4,5}}
\put(128.00,-5.50){\line(1,0){12}}

%\put(-5.00,-20.00){$12321$}

%%%%%%%%%%%%% path1 %%%%%%%%%%%%%%%% 
\put(-5.00,-22.00){\line(1,1){4.00}}
\put(-1.00,-18.00){\line(1,1){4.00}}
\put(3.00,-14.00){\line(1,-1){4.00}}
\put(7.00,-18.00){\line(1,-1){4.00}}
\put(-5.00,-22.00){\circle*{1.00}}
\put(-1.00,-18.00){\circle*{1.00}}
\put(3.00,-14.00){\circle*{1.00}}
\put(7.00,-18.00){\circle*{1.00}}
\put(11.00,-22.00){\circle*{1.00}}
%%%%%%%%%%%%%%%%%%%%%%%%%%%%%%%%%%%%%

\put(20.00,-22.00){\circle*{1.00}}
\put(23.00,-22.00){\circle*{1.25}}
\put(26.00,-22.00){\circle*{1.60}}
\put(29.00,-22.00){\circle*{1.25}}
\put(32.00,-22.00){\circle*{1.00}}
\put(20.00,-22.00){\line(0,1){4,5}}
\put(23.00,-22.00){\line(0,1){4,5}}
\put(26.00,-22.00){\line(0,1){4,5}}
\put(29.00,-22.00){\line(0,1){4,5}}
\put(32.00,-22.00){\line(0,1){4,5}}
\put(20.00,-17.50){\line(1,0){12}}

\put(38.00,-22.00){\circle*{1.00}}
\put(41.00,-22.00){\circle*{1.25}}
\put(44.00,-22.00){\circle*{1.60}}
\put(47.00,-22.00){\circle*{1.25}}
\put(50.00,-22.00){\circle*{1.00}}
\put(38.00,-22.00){\line(0,1){4,5}}
\put(41.00,-22.00){\line(0,1){2,5}}
\put(44.00,-22.00){\line(0,1){2,5}}
\put(47.00,-22.00){\line(0,1){2,5}}
\put(50.00,-22.00){\line(0,1){4,5}}
\put(38.00,-17.50){\line(1,0){12}}
\put(41.00,-19.50){\line(1,0){6}}

\put(56.00,-22.00){\circle*{1.00}}
\put(59.00,-22.00){\circle*{1.25}}
\put(62.00,-22.00){\circle*{1.60}}
\put(65.00,-22.00){\circle*{1.25}}
\put(68.00,-22.00){\circle*{1.00}}
\put(56.00,-22.00){\line(0,1){4,5}}
\put(59.00,-22.00){\line(0,1){4,5}}
%\put(62.00,-20.00){\line(0,1){4,5}}
\put(65.00,-22.00){\line(0,1){4,5}}
\put(68.00,-22.00){\line(0,1){4,5}}
\put(56.00,-17.50){\line(1,0){12}}

\put(74.00,-22.00){\circle*{1.00}}
\put(77.00,-22.00){\circle*{1.25}}
\put(80.00,-22.00){\circle*{1.60}}
\put(83.00,-22.00){\circle*{1.25}}
\put(86.00,-22.00){\circle*{1.00}}
\put(74.00,-22.00){\line(0,1){4,5}}
\put(77.00,-22.00){\line(0,1){2,5}}
\put(83.00,-22.00){\line(0,1){2,5}}
\put(86.00,-22.00){\line(0,1){4,5}}
\put(74.00,-17.50){\line(1,0){12}}
\put(77.00,-19.50){\line(1,0){6}}

\end{picture}
\caption{Families $\mathcal{NC}_{{\rm irr}}(w)$ for $w\in \mathpzc{M}_5$.}
\end{figure}

\begin{Proposition}
For any $n\in \mathbb{N}$, there is a decomposition
\[
\mathcal{NC}_{n}\cong 
\bigsqcup_{\pi_0\in {\rm NC}(n)}\left(\pi_{0}\times \mathpzc{M}(\pi_0)\right),
\]
where each component is endowed with the product order inherited from ${\rm NC}(n)\times \mathpzc{M}_{n}$, 
where $\pi_0\times \mathpzc{M}(\pi_0)$ denotes the set of pairs $(\pi_0, w)$ with $w\in \mathpzc{M}(\pi_0)$.
Here
\[
\mathpzc{M}(\pi_0):=\{w\in \mathpzc{M}_{n}: (\pi_0,w)\in \mathcal{NC}(w)\}.
\]
In particular, $\mathpzc{M}(\hat{1}_n)=\mathpzc{M}_n$ for any $n\in \mathbb{N}$. 
An analogous decomposition holds for $\mathcal{NC}_{n, {\rm irr}}$.
\end{Proposition}

{\it Proof.}
These decompositions follow directly from Definition 3.7. \hfill $\boxvoid$\\

\begin{figure}
\unitlength=1mm
\special{em.linewidth 0.5pt}
\linethickness{0.5pt}
\begin{picture}(140.00,120.00)(-45.00,0.00)

%%%%%%%%%%%%%%%%% minidiamond %%%%%%%%%%%%%%%%%%%

\put(-40.00,38.00){${\rm NC}_{{\rm irr}}(4)$}
\put(-35.00,34.00){\circle*{2.00}}
\put(-35.00,22.00){\circle*{2.00}}
\put(-35.00,10.00){\circle*{2.00}}
\put(-29.00,22.00){\circle*{2.00}}
\put(-41.00,22.00){\circle*{2.00}}

\put(-33.50,33.00){$\scriptstyle{A}$}
\put(-33.50,21.00){$\scriptstyle{C}$}
\put(-33.50,9.00){$\scriptstyle{E}$}
\put(-28.00,21.00){$\scriptstyle{D}$}
\put(-40.00,21.00){$\scriptstyle{B}$}

\put(-35.00,10.00){\line(0,1){24}}
\put(-35.00,10.00){\line(1,2){6}}
\put(-35.00,10.00){\line(-1,2){6}}
\put(-41.00,22.00){\line(1,2){6}}
\put(-29.00,22.00){\line(-1,2){6}}

%%%%%%%%%%%%%%%%%% minisquare %%%%%%%%%%%%%%%%%

\put(-9.00,38.00){$\mathpzc{M}_{4}$}
\put(-5.00,34.00){\line(-1,-1){12}}
\put(-5.00,34.00){\line(1,-1){12}}
\put(-5.00,10.00){\line(-1,1){12}}
\put(-5.00,10.00){\line(1,1){12}}

\put(-5.00,34.00){\circle*{2.00}}
\put(-17.00,22.00){\circle*{2.00}}
\put(-5.00,10.00){\circle*{2.00}}
\put(7.00,22.00){\circle*{2.00}}

\put(-3.00,34.00){$\scriptscriptstyle{1221}$}
\put(-15.00,22.00){$\scriptscriptstyle{1211}$}
\put(-3.00,10.00){$\scriptscriptstyle{1111}$}
\put(9.00,22.00){$\scriptscriptstyle{1121}$}

%%%%%%%%%%%%% diamond 1 %%%%%%%%%%%%%%%% 

\put(38.00, 106.00){${\rm NC}_{{\rm irr}}(4)\times \mathpzc{M}_{4}$}

\put(50.00,100.00){\circle*{2.00}}
\put(50.00,80.00){\circle*{2.00}}
\put(50.00,60.00){\circle*{2.00}}
\put(60.00,80.00){\circle*{2.00}}
\put(40.00,80.00){\circle*{2.00}}

\put(50.00,100.00){\circle{6.00}}
\put(50.00,80.00){\circle{6.00}}
\put(50.00,60.00){\circle{6.00}}
\put(60.00,80.00){\circle{6.00}}
\put(40.00,80.00){\circle{6.00}}

\put(50.00,60.00){\line(0,1){40}}
\put(50.00,60.00){\line(1,2){10}}
\put(50.00,60.00){\line(-1,2){10}}
\put(40.00,80.00){\line(1,2){10}}
\put(60.00,80.00){\line(-1,2){10}}

%%%%%%%%%%%% lines connecting diamond 1 %%%%%%%%%%%%%%%%%%

\put(50.00,100.00){\line(-2,-1){40}}
\put(50.00,100.00){\line(1,-1){30}}
\put(10.00,80.00){\line(1,-1){30}}
\put(40.00,50.00){\line(2,1){40}}

%%%%%%%%%%% diamond 2 %%%%%%%%%%%%%
\put(10.00,80.00){\circle*{2.00}}
\put(10.00,60.00){\circle*{2.00}}
\put(10.00,40.00){\circle*{2.00}}
\put(20.00,60.00){\circle*{2.00}}
\put(0.00,60.00){\circle*{2.00}}

\put(10.00,80.00){\circle{6.00}}
\put(0.00,60.00){\circle{6.00}}

\put(10.00,40.00){\line(0,1){40}}
\put(10.00,40.00){\line(1,2){10}}

\put(10.00,40.00){\line(-1,2){10}}
\put(0.00,60.00){\line(1,2){10}}
\put(20.00,60.00){\line(-1,2){10}}

%%%%%%%%%%% diamond 3 %%%%%%%%%%%%%%

\put(40.00,50.00){\circle*{2.00}}
\put(40.00,30.00){\circle*{2.00}}
\put(40.00,10.00){\circle*{2.00}}
\put(50.00,30.00){\circle*{2.00}}
\put(30.00,30.00){\circle*{2.00}}

\put(40.00,50.00){\circle{6.00}}

\put(40.00,10.00){\line(0,1){40}}
\put(40.00,10.00){\line(1,2){10}}
\put(40.00,10.00){\line(-1,2){10}}
\put(30.00,30.00){\line(1,2){10}}
\put(50.00,30.00){\line(-1,2){10}}

%%%%%%%%%%%% damond 4 %%%%%%%%%%%%%%

\put(80.00,70.00){\circle*{2.00}}
\put(80.00,50.00){\circle*{2.00}}
\put(80.00,30.00){\circle*{2.00}}
\put(90.00,50.00){\circle*{2.00}}
\put(70.00,50.00){\circle*{2.00}}

\put(80.00,70.00){\circle{6.00}}
\put(90.00,50.00){\circle{6.00}}

\put(80.00,30.00){\line(0,1){40}}
\put(80.00,30.00){\line(1,2){10}}
\put(80.00,30.00){\line(-1,2){10}}
\put(70.00,50.00){\line(1,2){10}}
\put(90.00,50.00){\line(-1,2){10}}

%%%%%%%%%%% lines of type 1 %%%%%%%%%%%%
\put(50.00,80.00){\line(-2,-1){40}}
\put(50.00,60.00){\line(-2,-1){40}}
\put(60.00,80.00){\line(-2,-1){40}}
\put(40.00,80.00){\line(-2,-1){40}}

\put(80.00,50.00){\line(-2,-1){40}}
\put(80.00,30.00){\line(-2,-1){40}}
\put(90.00,50.00){\line(-2,-1){40}}
\put(70.00,50.00){\line(-2,-1){40}}

%%%%%%%%%%%% lines of type 2 %%%%%%%%%

\put(50.00,80.00){\line(1,-1){30}}
\put(50.00,60.00){\line(1,-1){30}}
\put(60.00,80.00){\line(1,-1){30}}
\put(40.00,80.00){\line(1,-1){30}}

\put(10.00,60.00){\line(1,-1){30}}
\put(10.00,40.00){\line(1,-1){30}}
\put(20.00,60.00){\line(1,-1){30}}
\put(0.00,60.00){\line(1,-1){30}}

\end{picture}
\caption{The Hasse diagram of ${\rm NC}_{{\rm irr}}(4)\times \mathpzc{M}_{4}$. The vertices of the subposet
$\mathcal{NC}_{4, {\rm irr}}=\bigsqcup_{w\in \mathpzc{M}_{4}}\mathcal{NC}_{{\rm irr}}(w)$ are marked with larger circles.
Each of them corresponds to a non-zero contribution to the free cumulant $r_4$.}
\end{figure}

\begin{Example}
{\rm 
Consider ${\rm NC}_{{\rm irr}}(4)\times \mathpzc{M}_{4}$, the Cartesian product of a 
diamond lattice and a square lattice, together with its subposet $\mathcal{NC}_{4, {\rm irr}}$ 
shown in Fig.~3. Each diamond in this figure corresponds to a
word $w\in \mathpzc{M}_{4}$ and contains five vertices $A,B,C,D,E$, associated with  
the following partitions: 
\[ 
\pi_{A}=\{\{1,2,3,4\}\}, \;\;\;\pi_{B}=\{\{1,3,4\}, \{2\}\}, \;\;\;\pi_{C}=\{\{1,4\}, \{2,3\}\},
\]
\[
\pi_{D}=\{\{1,2,4\}, \{3\}\},\;\;\;\pi_{E}=\{\{1,4\}, \{2\}, \{3\}\}.
\]
In Fig.~3, the vertices belonging to
$\mathcal{NC}_{4,{\rm irr}}$ are indicated by large circles. The definition of this poset yields a decomposition of the diamond lattices into smaller pieces: the diamond corresponding to $12^21$ is entirely contained in $\mathcal{NC}_{4,{\rm irr}}$, those corresponding to $121^2$ and $1^221$ contribute two-element sublattices, and the diamond corresponding to $1^4$ contributes only its maximal element.
We will see later that the definition of $\mathcal{NC}_{4, {\rm irr}}$ corrresponds to the decomposition of $r_4$ 
as a sum of four `partial free cumulants' called Motzkin cumulants. In Fig.~3, Motzkin cumulants correspond to the 
pieces of the diamonds, whereas the decomposition in Proposition 3.8 corresponds to the horizontal layers.
Note that the poset $\mathcal{NC}_{4, {\rm irr}}$
is not a lattice since; for instance, the elements $(\pi_{C}, 12^21)$ and $(\pi_{A}, 1^221)$ do not admit a common lower bound. 
}
\end{Example}

Now, we shall give a {\it constructive} definition of $\mathcal{NC}(w)$ which is both intuitive and 
useful for proofs. We begin by recalling the definition of a nice subset of $\mathcal{NC}(w)$ consisting 
of noncrossing partitions which are monotonically adapted to $w$, denoted by $\mathcal{M}(w)$.
These partitions were introduced in \cite{[L6]}, where it was shown that each $\mathcal{M}(w)$ is a lattice. 
They are characterized by having one-color blocks whose depths coincide with those colors, and thus provide a 
natural decomposition of ${\rm NC}(n)$. 

\begin{Definition}
{\rm 
Let $\pi =(\pi_0,w)\in \mathcal{NC}(w)$, where $w\in \mathpzc{M}$. 
If each block $v_j\in \pi$ is a constant word and satisfies
\[
d(v_j)=h(v_j),
\]
then $\pi$ is called {\it monotonically adapted to} $w$.
We denote by $\mathcal{M}(w)$ the lattice of such partitions, and by
$\mathcal{M}_{{\rm irr}}(w)$ its subset of irreducible elements.
}
\end{Definition}

Each $\mathcal{M}(w)$ 
has a greatest element, namely the coarsest noncrossing partition adapted to $w$ in which there is at most one block of a given color
under any fixed bridge. It also has a least element, denoted by $\hat{0}_{w}$, in which the number of blocks of a given color
under each bridge is maximal subject to the depth condition. Clearly,
$\mathcal{M}(w)\subseteq \mathcal{NC}(w)$ for any $w$. Moreover, $\hat{0}_{w}$ is also the least element in 
$\mathcal{NC}(w)$. 

In contrast, the greatest element of $\mathcal{M}(w)$ does not coincide with the greatest element of $\mathcal{NC}(w)$ unless $w$ is 
constant. For this reason, we reserve the notation $\hat{1}_{w}$ for 
the partition consisting of a single block $\{w\}$, which is the greatest element of both $\mathcal{NC}(w)$ and
$\mathcal{NC}_{{\rm irr}}(w)$ for every $w$.

\begin{Example}
{\rm 
Consider $w=12^31$. From the set $\mathcal{NC}(w)=\mathcal{NC}_{{\rm irr}}(w)$,
consisting of 13 partitions given in Fig.~2, we select those that are monotonically adapted to $w$:

\unitlength=1mm
\special{em.linewidth 0.5pt}
\linethickness{0.5pt}
\begin{picture}(140.00,20.00)(18.00,-18.00)

\put(44.00,-10.00){\line(1,1){4.00}}
\put(48.00,-06.00){\line(1,0){4.00}}
\put(52.00,-06.00){\line(1,0){4.00}}
\put(56.00,-06.00){\line(1,-1){4.00}}
\put(44.00,-10.00){\circle*{1.00}}
\put(48.00,-06.00){\circle*{1.00}}
\put(52.00,-06.00){\circle*{1.00}}
\put(56.00,-06.00){\circle*{1.00}}
\put(60.00,-10.00){\circle*{1.00}}

\put(36.00,-9.00){$\mathcal{M}\big($}
\put(61.00,-9.00){$\big)=$}

%%%%%%%%%%%%%%%%%%%%%%%%%%%%%%%%%%%%%%%%%

\put(70.00,-9.00){$\big\{$}

\put(74.00,-10.00){\circle*{1.00}}
\put(77.00,-10.00){\circle*{1.25}}
\put(80.00,-10.00){\circle*{1.25}}
\put(83.00,-10.00){\circle*{1.25}}
\put(86.00,-10.00){\circle*{1.00}}
\put(74.00,-10.00){\line(0,1){4,5}}
\put(77.00,-10.00){\line(0,1){2,5}}
\put(80.00,-10.00){\line(0,1){2,5}}
\put(83.00,-10.00){\line(0,1){2,5}}
\put(86.00,-10.00){\line(0,1){4,5}}
\put(74.00,-5.50){\line(1,0){12}}
\put(77.00,-7.50){\line(1,0){6}}

\put(88.00,-10.00){,}

\put(92.00,-10.00){\circle*{1.00}}
\put(95.00,-10.00){\circle*{1.25}}
\put(98.00,-10.00){\circle*{1.25}}
\put(101.00,-10.00){\circle*{1.25}}
\put(104.00,-10.00){\circle*{1.00}}
\put(92.00,-10.00){\line(0,1){4,5}}
\put(95.00,-10.00){\line(0,1){2,5}}
\put(98.00,-10.00){\line(0,1){2,5}}
\put(104.00,-10.00){\line(0,1){4,5}}
\put(92.00,-5.50){\line(1,0){12}}
\put(95.00,-7.50){\line(1,0){3}}

\put(106.00,-10.00){,}

\put(110.00,-10.00){\circle*{1.00}}
\put(113.00,-10.00){\circle*{1.25}}
\put(116.00,-10.00){\circle*{1.25}}
\put(119.00,-10.00){\circle*{1.25}}
\put(122.00,-10.00){\circle*{1.00}}
\put(110.00,-10.00){\line(0,1){4,5}}
\put(116.00,-10.00){\line(0,1){2,5}}
\put(119.00,-10.00){\line(0,1){2,5}}
\put(122.00,-10.00){\line(0,1){4,5}}
\put(110.00,-5.50){\line(1,0){12}}
\put(116.00,-7.50){\line(1,0){3}}

\put(124.00,-10.00){,}

\put(128.00,-10.00){\circle*{1.00}}
\put(131.00,-10.00){\circle*{1.25}}
\put(134.00,-10.00){\circle*{1.25}}
\put(137.00,-10.00){\circle*{1.25}}
\put(140.00,-10.00){\circle*{1.00}}
\put(128.00,-10.00){\line(0,1){4,5}}
\put(140.00,-10.00){\line(0,1){4,5}}
\put(128.00,-5.50){\line(1,0){12}}

\put(141.00,-9.00){$\big\}$}

\end{picture}
\\
\noindent
The last partition in this set is the least element $\hat{0}_w$ of $\mathcal{M}(w)$, whereas the first partition
is its greatest element. This example already shows that $\mathcal{NC}(w)$ is, in general, much larger than $\mathcal{M}(w)$, which 
explains why the underlying theory is more involved. }
\end{Example}

\begin{Definition}
{\rm Let $\pi, \pi'\in \mathcal{NC}(w)$, with $\pi\prec \pi'$. We say that $\pi'$ is an {\it admissible coarsening}
of $\pi$ if it is obtained by a finite sequence of operations of two types:
\begin{enumerate}
\item {\it juxtapositions}
\[
(v_1, \ldots, v_j,v_{j+1},\ldots, v_m)\rightarrow (v_1, \ldots, v_jv_{j+1}, \ldots , v_m)
\]
whenever $v_j,v_{j+1}$ are consecutive blocks lying under the same bridge,
\item {\it insertions}
\[
(v_1, \ldots , v_j, v_{j+1}, \ldots , v_m)\rightarrow (v_1, \ldots , v_j'v_{j+1}v_{j}'', \ldots , v_m)
\]
whenever $v_j=v_j'v_j''$ and $v_j$ is the nearest outer block of $v_{j+1}$,
\end{enumerate}
where $1\leq j \leq n-1$. In this case, $\pi$ is called an {\it admissible refinement} of $\pi'$. 
We denote by $\mathcal{AC}(\pi)$ the set of admissible coarsenings of $\pi$, and by 
$\mathcal{AR}(\pi')$ the set of admissible refinements of $\pi'$. 
If $\pi'$ is obtained from $\pi$ using only insertions, then $\pi$ 
is called an {\it ancestor} of $\pi$; the set of all ancestors of $\pi$ is denoted by $\mathcal{A}(\pi)$. 
}
\end{Definition}

\begin{Example}
{\rm Admissible coarsenings can be used to construct $\mathcal{NC}(w)$ from $\hat{0}_{w}$.
Let us construct a partition $\pi \in \mathcal{NC}_{{\rm irr}}(w)$ from $\hat{0}_{w}$, where $w=123^321$,
by applying a juxtaposition followed by two insertions, 

\unitlength=1mm
\special{em.linewidth 0.5pt}
\linethickness{0.5pt}
\begin{picture}(140.00,20.00)(10.00,-17.00)

%%%%%%%%%%%%%%%%%%%%%%%%%%%%%%%%%%%%%%%%%

\put(31.00,-10.00){\circle*{0.75}}
\put(34.00,-10.00){\circle*{1.00}}
\put(37.00,-10.00){\circle*{1.50}}
\put(40.00,-10.00){\circle*{1.50}}
\put(43.00,-10.00){\circle*{1.50}}
\put(46.00,-10.00){\circle*{1.00}}
\put(49.00,-10.00){\circle*{0.75}}
\put(31.00,-10.00){\line(0,1){7}}
\put(34.00,-10.00){\line(0,1){4,5}}
\put(46.00,-10.00){\line(0,1){4,5}}
\put(49.00,-10.00){\line(0,1){7}}
\put(34.00,-5.50){\line(1,0){12}}
\put(31.00,-3.00){\line(1,0){18}}

\put(51.00,-08.00){$\rightarrow$}

\put(39.00,-15.00){$\scriptstyle{\hat{0}_{w}}$}

\put(58.00,-10.00){\circle*{0.75}}
\put(61.00,-10.00){\circle*{1.00}}
\put(64.00,-10.00){\circle*{1.50}}
\put(67.00,-10.00){\circle*{1.50}}
\put(70.00,-10.00){\circle*{1.50}}
\put(73.00,-10.00){\circle*{1.00}}
\put(76.00,-10.00){\circle*{0.75}}
\put(58.00,-10.00){\line(0,1){7}}
\put(61.00,-10.00){\line(0,1){4,5}}
\put(67.00,-10.00){\line(0,1){2,5}}
\put(70.00,-10.00){\line(0,1){2,5}}
\put(73.00,-10.00){\line(0,1){4,5}}
\put(76.00,-10.00){\line(0,1){7}}

\put(58.00,-3.00){\line(1,0){18}}
\put(61.00,-5.50){\line(1,0){12}}
\put(67.00,-7.50){\line(1,0){3}}

\put(66.00,-15.00){$\scriptstyle{\pi_1}$}

\put(78.00,-08.00){$\rightarrow$}

\put(85.00,-10.00){\circle*{0.75}}
\put(88.00,-10.00){\circle*{1.00}}
\put(91.00,-10.00){\circle*{1.50}}
\put(94.00,-10.00){\circle*{1.50}}
\put(97.00,-10.00){\circle*{1.50}}
\put(100.00,-10.00){\circle*{1.00}}
\put(103.00,-10.00){\circle*{0.75}}
\put(85.00,-10.00){\line(0,1){7}}
\put(88.00,-10.00){\line(0,1){4,5}}
\put(94.00,-10.00){\line(0,1){4,5}}
\put(97.00,-10.00){\line(0,1){4,5}}
\put(100.00,-10.00){\line(0,1){4,5}}
\put(103.00,-10.00){\line(0,1){7}}

\put(85.00,-3.00){\line(1,0){18}}
\put(88.00,-5.50){\line(1,0){12}}

\put(105.00,-08.00){$\rightarrow$}
\put(93.00,-15.00){$\scriptstyle{\pi_2}$}

\put(112.00,-10.00){\circle*{0.75}}
\put(115.00,-10.00){\circle*{1.00}}
\put(118.00,-10.00){\circle*{1.50}}
\put(121.00,-10.00){\circle*{1.50}}
\put(124.00,-10.00){\circle*{1.50}}
\put(127.00,-10.00){\circle*{1.00}}
\put(130.00,-10.00){\circle*{0.75}}
\put(112.00,-10.00){\line(0,1){7}}
\put(115.00,-10.00){\line(0,1){7}}
\put(121.00,-10.00){\line(0,1){7}}
\put(124.00,-10.00){\line(0,1){7}}
\put(127.00,-10.00){\line(0,1){7}}
\put(130.00,-10.00){\line(0,1){7}}
\put(112.00,-3.00){\line(1,0){18}}

\put(120.00,-15.00){$\scriptstyle{\pi}$}

\put(136.00,-10.00){,}

\end{picture}
\\
\noindent
which yields a chain $\hat{0}_{w}\prec \pi_1\prec \pi_2 \prec \pi$ in $\mathcal{NC}(w)$.
In terms of diagrams, a juxtaposition corresponds to a {\it horizontal merging} of two neighboring blocks of the same depth
lying under a common bridge, whereas an insertion corresponds to a {\it vertical merging} of a block with its nearest outer block.
}
\end{Example}

\begin{Lemma}
Let $w\in \mathpzc{M}$, and let $\hat{0}_{w}$ be the least  
element in $\mathcal{NC}(w)$. Then
\[
\mathcal{NC}(w)=\mathcal{AC}(\hat{0}_{w}).
\]
\end{Lemma}
{\it Proof.}
If $w$ is a constant word, then the assertions is immediate,
since in that case $\mathcal{NC}(w)\cong {\rm Int}(n)\cong\mathcal{AC}(\hat{0}_{w})$, where
$\hat{0}_{w}$ consists of $n$ singletons and $n=|w|$.  
Thus, assume that $w$ is not constant. 
By definition of admissible coarsenings, we have $\mathcal{AC}({0}_{w})\subseteq \mathcal{NC}(w)$. 
It therefore suffices to show that each 
$\pi\in\mathcal{NC}(w)$ can be obtained from $\hat{0}_{w}$ by a sequence of operations from Definition 3.10. 

Suppose that $w$ is built from letters $1, \ldots, d$. We begin with a sequence of juxtapositions (horizontal mergings). 
Whenever two neighboring letters $d$ belong to the same block of $\pi$, we juxtapose them in $\hat{0}_{w}$. 
Repeating this operation for all such pairs at depth $d$, we obtain 
an admissible coarsening $\pi_1$ of $\hat{0}_{w}$ 
By construction, 
\[
\hat{0}_{w}\preceq \pi_1\preceq \pi.
\] 
We then proceed similarly for blocks of depth and color $d-1$.
Namely, if two neighboring blocks of depth and color $d-1$ lie under the same bridge
and are contained in a single block of $\pi$, we juxtapose them. 
Applying this procedure yields an admissible coarsening $\pi_2$ of $\pi_1$. 
Continuing in this way for all remaining depths, we obtain a sequence
$\pi_1, \ldots \pi_d\in \mathcal{M}(w)$ such that 
\[ 
\hat{0}_{w}\preceq \pi_1\preceq  \ldots \preceq\pi_d\preceq \pi.
\]
By construction, $\pi_d$ is the maximal refinement of $\pi$ that is in $\mathcal{M}(w)$.

We now perform insertions (vertical mergings) starting from $\pi_d$. First consider blocks of depth $d$: whenever an inner block of depth $d$ and its nearest outer block are contained in the same block of $\pi$, we insert the inner block into the outer one. Repeating this operation for all such blocks yields $\pi_{d+1}$. Proceeding similarly for all remaining depths, we obtain a sequence
of partitions $\pi_{d+1}, \ldots, \pi_{2d-1}$, such that
\[
\pi_d\preceq \pi_{d+1}\preceq \cdots \preceq \pi_{2d-1}\preceq \pi.
\]
Since all operations are admissible, we have $\pi_{2d-1}\in \mathcal{NC}(w)$. Moreover, by construction, $\pi_{2d-1}=\pi$, as each step performs the minimal coarsening merging two subblocks of a block of $\pi$ into a single block. This completes the proof. \hfill $\boxvoid$

\begin{Theorem}
For any $w\in \mathpzc{M}_n$, the set $\mathcal{NC}(w)$ is a lattice with the partial order induced from ${\rm NC}(n)$.
\end{Theorem}

{\it Proof.}
If $w=1^n$, then $\mathcal{NC}(w)\cong {\rm Int}(n)$, hence it is a lattice. 
Assume $w\neq 1^n$, and let $\pi,\pi'\in \mathcal{NC}(w)$.
We are going to construct a common coarsening of $\pi$ and $\pi'$ by following the 
admissible merging procedure from the proof of Lemma 3.14, applied 
simultaneously to $\pi$ and $\pi'$.

\medskip

\noindent
{\it Horizontal mergings.}
We construct inductively a sequence
\[
\hat{0}_w = \rho_0 \preceq \rho_1 \preceq \cdots \preceq \rho_d,
\]
where at step $j$ we perform horizontal mergings at depth $d-j+1$.
More precisely, suppose $\rho_{j-1}$ has been constructed. 
Whenever two neighboring blocks of depth and color $d-j+1$ in $\rho_{j-1}$, 
lying under the same bridge, are contained in a single block 
of $\pi$ or of $\pi'$, we juxtapose them. 
This produces $\rho_j$. By the same argument as in Lemma 3.14, each such step 
is an admissible coarsening, hence
\[
\hat{0}_w \preceq \rho_1 \preceq \cdots \preceq \rho_d,
\qquad \rho_j \in \mathcal{NC}(w).
\]
Moreover, $\rho_d$ is the maximal partition which is monotonically adapted 
to $w$ and is dominated by any common coarsening of $\pi$ and $\pi'$.

\medskip

\noindent
{\it Vertical mergings.}
Starting from $\rho_d$, we construct
\[
\rho_d \preceq \rho_{d+1} \preceq \cdots \preceq \rho_{2d-1}.
\]
At step $j=d+1,\ldots,2d-1$, we consider blocks of depth $2d-j$. 
Whenever an inner block and its nearest outer block are contained in a single 
block of $\pi$ or of $\pi'$, we merge them. Again, by the same argument as in Lemma 3.14, each step is an admissible 
coarsening, hence $\rho_j \in \mathcal{NC}(w)$ for all $j$, and in particular 
$\rho_{2d-1}\in \mathcal{NC}(w)$.

By construction, every merging performed in the construction of $\pi$ or $\pi'$ 
is also performed in the construction of $\rho_{2d-1}$. Hence
\[
\pi \preceq \rho_{2d-1}, \qquad \pi' \preceq \rho_{2d-1}.
\]
Conversely, at each step we perform only those mergings which are necessary 
to realize mergings present in $\pi$ or $\pi'$. Therefore $\rho_{2d-1}$ is 
the minimal admissible common coarsening of $\pi$ and $\pi'$, and hence
\[
\rho_{2d-1} = \pi \vee \pi'.
\]
which completes the proof.\hfill $\boxvoid$\\

\section{Motzkin cumulants}

The concept of noncrossing partitions adapted to Motzkin words allows us to introduce the associated cumulants called 
Motzkin cumulants. The construction proceeds in two steps: we first define an intermediate family of cumulants, called $w$-Boolean cumulants, and then use them to introduce the Motzkin cumulants.

We now adapt the notation introduced in Section 2 to the setting of partitions indexed by Motzkin words, as defined in Section 3. In this context, multilinear operator-valued functions will be defined recursively in analogy with Definitions 2.4–2.5 and Theorem 2.6, with ${\mathbb N}$ replaced by $\mathpzc{AM}$. 
The corresponding notation will be
\[
f_{w}(a_1, \ldots, a_n), \quad f_{v}[a_1, \ldots, a_n], \quad \text{and} \quad f_{\pi}[a_{1}, \ldots, a_{n}],
\]
where $w\in \mathpzc{AM}$ and $v\in \pi\in \mathcal{NC}(w)$. Here, $f_{\pi}[a_1, \ldots, a_n]$ is understood as a nested version of the product denoted by the same symbol in Definition 2.4, where the nesting is induced by the block structure of $\pi$. More precisely, each block of $\pi\in \mathcal{NC}(w)$ is viewed as a pair $(V,v)$, where $V$ is a block of $\pi_0\in {\rm NC}(|w|)$ and $v$ is a subword of $w$. In particular, $f_w$ corresponds to the one-block partition $\hat{1}_{w}$. We will also use the notation $f(w)$ and $f(\pi)$ instead of $f_w$ and $f_{\pi}$ when $w$ and $\pi$ are represented diagrammatically.

The constructions below are motivated by the relation between Boolean and free cumulants described in Theorem 2.6. Replacing ${\rm NC}_{{\rm irr}}(n)$ by $\mathcal{NC}_{{\rm irr}}(w)$ leads to a refined combinatorial structure, reflecting the additional information carried by Motzkin paths. In particular, since Motzkin paths carry more structure than natural numbers, the corresponding lattices involve more intricate nesting, in which cumulants associated with nested blocks of $\pi\in \mathcal{NC}_{{\rm irr}}(w)$ depend on the heights of the corresponding subpaths of $w$.

\begin{Definition}
{\rm 
An operator-valued noncommutative probability space is a triple of the form
$({\mathcal C}, {\mathcal B}, E)$, where
${\mathcal C}$ is a unital algebra, ${\mathcal B}$ is a unital subalgebra, and
$E:\mathcal{C}\rightarrow \mathcal{B}$ is a conditional expectation, i.e.
\[
E(b)=b\;\;\;{\rm and}\;\;\; E(b_1ab_2)=b_1E(a)b_2
\]
for any $b,b_1,b_2\in \mathcal{B}$ and any $a\in {\mathcal C}$.
}
\end{Definition}

In order to effectively use the lattices $\mathcal{NC}_{{\rm irr}}(w)$ rather than $\mathcal{NC}(w)$—which is both convenient and simplifies computations—we first introduce an intermediate family of cumulants, which we call operator-valued $w$-Boolean cumulants. These will play a role analogous to Boolean cumulants in the classical setting and will lead to the definition of Motzkin cumulants.
 
Let ${\mathcal Int}(w)$ denote the set of interval partitions of $w$ into Motzkin subwords, namely
\[
{\mathcal Int}(w)=\{(v_1, \ldots , v_k): w=v_1\cdots v_k,\;v_1, \ldots, v_k\in \mathpzc{AM}\},
\]
where $v_1\cdots v_k$ is a concatenation of Motzkin words of height $j$ whenever $w$ is of height $j$. In particular, if $w\in \mathpzc{M}$, then $v_j\in \mathpzc{M}$ for each $j$.

\begin{Definition}
{\rm 
The {\it $w$-Boolean cumulants} are the family $\{B_{w}:w\in \mathpzc{AM}\}$ of multilinear functions 
\[
B_w:\mathcal{C}^{n}\rightarrow {\mathcal B}
\] 
defined by
\[
E(a_1 \cdots a_n)=\sum_{\pi\in {\mathcal Int}(w)} B_{\pi}[a_1, \ldots, a_n],
\] 
where $w=j_1\cdots j_n\in \mathpzc{AM}$, $a_1, \ldots, a_n\in \mathcal{C}$, and 
\[
B_{\pi}=\prod_{v\in \pi}B_{v},
\]
with $B_{v}:=B_{\hat{1}_{v}}$ for any block $v\in \pi$.
}
\end{Definition}

Solving the above relations for the functions $B_{w}$ is straightforward and analogous to the M\"{o}bius inversion formula for standard Boolean cumulants. 

\begin{Proposition}
If $w=j_1\ldots j_n\in \mathpzc{AM}_{n}$, then the following inversion formula holds:
\[
B_{w}(a_1, \ldots , a_n)=\sum_{\pi\in {\mathcal Int}(w)}(-1)^{|\pi|-1}E_{\pi}[a_1, \ldots, a_n],
\]
where $a_1, \ldots, a_{n}\in \mathcal{C}$, and 
\[
E_{\pi}[a_1, \ldots, a_n]=\prod_{v\in \pi}E_v[a_1, \ldots, a_n],
\] 
with $E_{v}[a_1, \ldots , a_n]=E(a_{i_1}\cdots a_{i_k})$ for a block $v=j_{i_1}\cdots j_{i_{k}}$. 
\end{Proposition}

{\it Proof.}
This follows by M\"{o}bius inversion on the lattice ${\mathcal Int}(w)$; see, for example, Anshelevich \cite{[An]} 
for the general moment–cumulant relations obtained via M\"{o}bius inversion on lattices of partitions.
\hfill $\boxvoid$

\begin{Example}
{\rm 
Since ${\mathcal Int}(j^n)\cong {\rm Int}(n)$, it follows that $B_{j^n}$ 
coincides with the operator-valued Boolean cumulant of order $n$ associated with $E$, for any $j$. 
If $w$ is not a constant word, then ${\mathcal Int}(w)$ can be identified 
with a proper subset of ${\rm Int}(n)$, where $n=|w|$. For instance,
\begin{eqnarray*}
B_{12^21}(a_1,a_2,a_3,a_4)&=& E(a_1a_2a_3a_4),\\
B_{121^2}(a_1,a_2,a_3,a_4)&=&E(a_1a_2a_3a_4)-E(a_1a_2a_3)E(a_4),\\
B_{121^3}(a_1,a_2,a_3,a_4,a_5)&=&E(a_1a_2a_3a_4a_5)-E(a_1a_2a_3a_4)E(a_5)\\
&&\!\!\!\!\!\!-\;E(a_1a_2a_3)E(a_4a_5)+E(a_1a_2a_3)E(a_4)E(a_5),
\end{eqnarray*}
and so on. Similar expressions are obtained for analogous words $w\in \mathpzc{AM}$ of height $j$. 
It is easy to see that each $B_{w}(a_1, \ldots, a_n)$ can be expressed as a $w$-Boolean cumulant of products of variables 
associated with shorter words of the same height. For instance,
\[
B_{121^2}(a_1,a_2,a_3,a_4)=B_{11}(a_1a_2a_3,a_4).
\]
Therefore, these quantities can be identified with {\it Boolean block cumulants}, with endpoints of 
blocks determined by the letters $j$. However, we will use the original quantities $B_{w}(a_1, \ldots, a_n)$, since they contain 
all the information about the word $w$ needed in the definition of Motzkin cumulants.
}
\end{Example}

\begin{Definition}
{\rm The {\it Motzkin cumulants} are the family $\{K_w:w\in \mathpzc{AM}\}$
of multilinear functions 
\[
K_{w}:{\mathcal C}^n\rightarrow {\mathcal B}
\]
defined by 
\[
B_{w}(a_1, \ldots, a_n)=\sum_{\pi\in \mathcal{NC}_{{\rm irr}}(w)}K_{\pi}[a_1, \ldots , a_n],
\]
where $w=j_1\cdots j_n\in \mathpzc{AM}$ and $a_1, \ldots, a_n\in \mathcal{C}$, and 
\[
K_{\pi}[a_1, \ldots, a_n]=
\left\{
\begin{array}{ll}
K_{w}(a_1, \ldots, a_n) & {\rm if}\; \pi=\hat{1}_{w},\\
K_{\pi\setminus v}[a_{1}, \ldots, a_{p-1}\,K_v(a_{p},\, \ldots , a_{q}),a_{q+1}, \ldots , a_{n}] & {\rm if} \; \pi \neq \hat{1}_{w},
\end{array}
\right.
\]
where $v\in \pi$ is any block with no inner blocks.
}
\end{Definition}

\begin{Remark}
{\rm
The recursive definition of $K_{\pi}$ is independent of the order in which 
innermost blocks are removed. Indeed, if several innermost blocks have the same 
nearest outer block, each reduction contributes a factor $b_i\in\mathcal B$ 
attached to the same leg of that outer block. Different orders of reduction 
yield only different parenthesizations of the same product
\[
b_1b_2\cdots b_r,
\]
and therefore give the same result by associativity in $\mathcal B$. Moreover,
if two innermost blocks have different nearest outer blocks, the corresponding 
reductions act independently and therefore commute trivially. Hence the recursion 
defining $K_{\pi}$ is well defined.
}
\end{Remark}

\begin{Remark}
{\rm
We make the following remarks on Definitions 4.2 and 4.5.
\begin{enumerate}

\item
The family $\{B_w:w\in\mathpzc{AM}\}$ is well defined, since for any
$w\in\mathpzc{AM}_n$ we can write
\begin{eqnarray*}
E(a_1\cdots a_n)
&=&
B_w(a_1,\ldots,a_n)
+\sum_{\stackrel{\pi\in{\mathcal Int}(w)}
{\scriptscriptstyle \pi\neq \hat1_w}}
B_\pi[a_1,\ldots,a_n].
\end{eqnarray*}
Thus each $B_w(a_1,\ldots,a_n)$ can be expressed recursively in terms of
the moment $E(a_1\cdots a_n)$ and lower-order mixed moments.

\item
Similarly, the family $\{K_w:w\in\mathpzc{AM}\}$ is well defined, since for
any $w\in\mathpzc{AM}_n$,
\begin{eqnarray*}
B_w(a_1,\ldots,a_n)
&=&
K_w(a_1,\ldots,a_n)
+\sum_{\stackrel{\pi\in\mathcal{NC}_{{\rm irr}}(w)}
{\scriptscriptstyle \pi\neq \hat1_w}}
K_\pi[a_1,\ldots,a_n].
\end{eqnarray*}
Hence each $K_w(a_1,\ldots,a_n)$ can be expressed recursively in terms of
$B_w(a_1,\ldots,a_n)$ and lower-order $w$-Boolean cumulants associated
with Motzkin subwords of $w$.

\item
Both definitions apply to arbitrary $a_1,\ldots,a_n\in\mathcal C$.
Therefore, for fixed $a_1,\ldots,a_n$, the cumulants
$B_w(a_1,\ldots,a_n)$ and $K_w(a_1,\ldots,a_n)$ are defined for several
different words $w$ when $n>2$. In particular, we can always compute these
cumulants for $w=j^n$. In that case
\[
K_{j^n}(a_1,\ldots,a_n)=B_{j^n}(a_1,\ldots,a_n),
\]
since
\[
\mathcal{NC}_{{\rm irr}}(j^n)=\{\hat1_{j^n}\}.
\]
Moreover, for fixed $a_1,\ldots,a_n$, the values
$B_w(a_1,\ldots,a_n)$ and $B_{w'}(a_1,\ldots,a_n)$ coincide whenever
$w'$ is the reduced Motzkin word associated with $w$ defined in Section~3. The same holds for
Motzkin cumulants.

\item
Using the bimodule property of the conditional expectation $E$, namely
$E(b)=b$ and
\[
E(bab')=bE(a)b',
\]
for $b,b'\in\mathcal B$ and $a\in\mathcal C$, together with the recursive
definitions of $w$-Boolean and Motzkin cumulants (equivalently, by induction
on $|w|$), we obtain
\[
B_j(b)=b,
\qquad
K_j(b)=b,
\]
for any one-letter word $j\in\mathpzc{AM}$ and $b\in\mathcal B$, and
\begin{eqnarray*}
B_w(b_0a_1b_1,\ldots,a_nb_n)
&=&
b_0B_w(a_1,b_1a_2,\ldots,b_{n-1}a_n)b_n,\\
K_w(b_0a_1b_1,\ldots,a_nb_n)
&=&
b_0K_w(a_1,b_1a_2,\ldots,b_{n-1}a_n)b_n,
\end{eqnarray*}
for any $b_0,b_1,\ldots,b_n\in\mathcal B$,
$a_1,\ldots,a_n\in\mathcal C$, and $w\in\mathpzc{AM}$.

\item
The definition of nested multilinear functions used in the formula for
$K_\pi$ is the same as in \cite{[MS]}. We give an example instead of a
formal definition. If the partition $\pi\in\mathcal{NC}(w)$, where
$w=123^22^2321$, is represented by the diagram

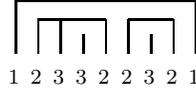
\unitlength=1mm
\special{em:linewidth 1pt}
\linethickness{0.5pt}
\begin{picture}(180.00,22.00)(-30.00,-19.00)

\put(22.00,-15.00){$\scriptstyle{1}$}
\put(25.00,-15.00){$\scriptstyle{2}$}
\put(28.00,-15.00){$\scriptstyle{3}$}
\put(31.00,-15.00){$\scriptstyle{3}$}
\put(34.00,-15.00){$\scriptstyle{2}$}
\put(37.00,-15.00){$\scriptstyle{2}$}
\put(40.00,-15.00){$\scriptstyle{3}$}
\put(43.00,-15.00){$\scriptstyle{2}$}
\put(46.00,-15.00){$\scriptstyle{1}$}

\put(23.00,-11.00){\line(0,1){7}}
\put(26.00,-11.00){\line(0,1){4,5}}
\put(29.00,-11.00){\line(0,1){4,5}}
\put(32.00,-11.00){\line(0,1){2,5}}
\put(35.00,-11.00){\line(0,1){4,5}}
\put(38.00,-11.00){\line(0,1){4,5}}
\put(41.00,-11.00){\line(0,1){2,5}}
\put(44.00,-11.00){\line(0,1){4,5}}
\put(47.00,-11.00){\line(0,1){7}}
\put(23.00,-4.00){\line(1,0){24}}
\put(26.00,-6.50){\line(1,0){9}}
\put(38.00,-6.50){\line(1,0){6}}

\end{picture}
\\
then the corresponding partitioned nested cumulant is
\[
K_\pi[a_1,\ldots,a_9]
=
K_{11}
(a_1
K_{232}(a_2,a_3K_3(a_4),a_5)
K_{22}(a_6K_3(a_7),a_8),
a_9).
\]
One may equivalently adopt the convention of attaching nested cumulants to
the variable immediately to the right of each inner block; by the
well-definedness established above, and associativity of the
$\mathcal B$-valued insertions, the two conventions give the same result.

\end{enumerate}
}
\end{Remark}

\begin{Example}
{\rm 
Clearly, $K_{1}(a_1)=B_1(a_1)$ and $K_{11}(a_1,a_2)=B_{11}(a_1,a_2)$ and the 
difference between $K_w$ and $B_w$ appears only when $|w|\geq 3$. 
The simplest Motzkin cumulant that does not coincide with $B_w$ corresponds to the word $w=121$
since in that case $\mathcal{NC}_{{\rm irr}}(w)=\{\{121\},\{1^2, 2\}\}$. We have
\begin{eqnarray*}
K_{111}(a_1,a_2,a_3)&=&B_{111}(a_1,a_2,a_3),\\
K_{121}(a_1,a_2,a_3)&=&B_{121}(a_1,a_2,a_3)-B_{11}(a_1B_{2}(a_2),a_3).
\end{eqnarray*}
Using Definition 4.5 and the above formulas, we obtain expressions for 
Motzkin cumulants associated with words of length four:
\begin{eqnarray*}
K_{1111}(a_1,a_2,a_3,a_4)&=&B_{1111}(a_1,a_2,a_3,a_4),\\
K_{1211}(a_1,a_2,a_3,a_4)&=&B_{1211}(a_1,a_2,a_3,a_4)-B_{111}(a_1B_{2}(a_2),a_3,a_4),\\
K_{1121}(a_1,a_2,a_3,a_4)&=&B_{1121}(a_1,a_2,a_3,a_4)-B_{111}(a_1,a_2B_{2}(a_3),a_4),\\
K_{1221}(a_1,a_2,a_3,a_4)&=&B_{1221}(a_1,a_2,a_3,a_4)-B_{121}(a_1B_{2}(a_2),a_3,a_4),\\
&&\!\!\!\!\!\!-\;B_{121}(a_1,a_2B_{2}(a_3),a_4)-B_{11}(a_1B_{22}(a_2,a_3),a_4)\\
&&\!\!\!\!\!\!+\; B_{11}(a_1B_2(a_2)B_2(a_3),a_4).
\end{eqnarray*}
These computations illustrate how the combinatorial structure of $\mathcal{NC}_{{\rm irr}}(w)$ governs the correction terms.
Each term on the right-hand side of these formulas corresponds to one of the ten vertices of the poset 
$\mathcal{NC}_{4,{\rm irr}}$ shown in Fig.~3. These formulas jointly resemble 
the expression for $r_4(a_1,a_2,a_3,a_4)$ in Example 2.7, although some partitions $\pi_0\in {\rm NC}(n)$
appear repeatedly on the right-hand side. 
}
\end{Example}

\begin{Proposition}
The definition of Motzkin cumulants is equivalent to the following moment-cumulant formula:
\[
E(a_1 \cdots a_n)=\sum_{\pi\in \mathcal{NC}(w)}K_{\pi}[a_1, \ldots, a_n],
\]
where $w\in \mathpzc{AM}$ and $a_1,\ldots , a_n\in \mathcal{C}$. 
\end{Proposition}

{\it Proof.}
The formula follows directly from the definitions of Motzkin cumulants and $w$-Boolean cumulants. Indeed, each partition $\pi\in \mathcal{NC}(w)$ admits a unique decomposition into irreducible components $\pi_k\in \mathcal{NC}_{{\rm irr}}(w_k)$ corresponding to a factorization $w=w_1\cdots w_p$, where $h(w_k)=h(w)$ for all $1\leq k\leq p$. Combining this decomposition with the defining relations for $B_w$ and $K_w$ yields the desired formula.
\hfill $\boxvoid$

\section{Inversion formula}

We would like to obtain a formula that expresses all Motzkin cumulants 
in terms of $w$-Boolean cumulants. 

\begin{Remark}
{\rm If $w=j^n$, then the corresponding Motzkin cumulants admit a simple expression in terms of mixed moments:
\[
K_{j^n}(a_1, \ldots , a_n)=\sum_{\pi\in {\mathcal Int}(j^n)}(-1)^{|\pi|-1}E_{\pi}[a_1, \ldots, a_n]
\]
by Remark 4.6(3) and Proposition 4.3, where ${\mathcal Int}(j^n)\cong {\rm Int}(n)$. 
If $w$ is not constant, the corresponding inversion formula is more involved.}
\end{Remark}

In order to treat the general case, we need to make some additional assumptions on the values of the 
$w$-Boolean cumulants assigned to blocks.

\begin{Definition}
{\rm
Let  $a_1, \ldots, a_n\in \mathcal{C}$ and let $w=j_1\cdots j_n\in \mathpzc{AM}$. If the following conditions hold:
\begin{enumerate}
\item if $v=j_k\cdots j_pj_{q+1}\cdots j_m,\;v'=j_{p+1}\cdots j_{q}\in \mathpzc{B}(w)$ and $h(v)=h(v')$, then
\[
B_{v}(a_k, \ldots , a_{p}B_{v'}(a_{p+1}, \ldots , a_{q}), a_{q+1}, \ldots , a_m)=0,
\]
\item if $v=j_{k}\cdots j_p, v'=j_{p+1}\cdots j_q\in \mathpzc{B}(w)$ and $h(v)\neq h(v')$, then 
\[
B_{v}(a_{k}, \ldots, a_{p})B_{v'}(a_{p+1}, \ldots, a_{q})=0,
\]
\end{enumerate}
where $\mathpzc{B}(w)$ is given by Definition 3.4,
we say that the family 
$\{B_{v}: v\in \mathpzc{B}(w)\}$ is {\it dedicated} to $(a_1, \ldots, a_n)$.
}
\end{Definition}

The conditions in Definition 5.2 encode the vanishing relations needed for the
recursive expansion of Motzkin cumulants into $w$-Boolean cumulants to preserve
adaptedness to $w$. In the proof of Theorem 5.3, they eliminate precisely the
two types of refinements which are not adapted to $w$: same-height refinements
violating the monotone height structure and mixed-height refinements violating
height separation. These vanishings are not automatic in general, but they hold
in the projection-valued setting considered later, where the relevant cumulants
satisfy the orthogonality relations required in Definition 5.2.

Although formally analogous to Theorem 2.6, the inversion formula below is
not a direct consequence of the usual free-Boolean relation. It is an adapted
M\"{o}bius inversion over $\mathcal{NC}_{{\rm irr}}(w)$, in which the
summation is constrained by the Motzkin structure of the word $w$. To state
it, we use the following notation. For $\pi \in \mathcal{NC}(w)$, we write
$B_{\pi}[a_1,\ldots,a_n]$ in the sense of Section~4, that is, as the
partitioned nested functional obtained by assigning to each block of $\pi$
the corresponding multilinear functional $B_v$, with the insertions
determined by the nesting structure of $\pi$. In particular, when
$\pi\in{\mathcal Int}(w)$, this reduces to the product notation used in
Definition~4.2.

\begin{Theorem}
Let $w=j_1\cdots j_n\in \mathpzc{AM}$ and let $(a_1, \ldots, a_n)\in {\mathcal C}^{n}$. If the family 
$\{B_v:v\in \mathpzc{B}(w)\}$ is dedicated to $(a_1, \ldots, a_n)$ (in the sense of Definition 5.2), then 
\[
K_w(a_1, \ldots , a_n)=\sum_{\pi\in \mathcal{NC}_{{\rm irr}}(w)} (-1)^{|\pi|-1}B_{\pi}[a_1, \ldots, a_n],
\]
where $|\pi|$ denotes the number of blocks of $\pi$.
\end{Theorem}
{\it Proof.}
If $w=j^n$, then $\mathcal{NC}_{{\rm irr}}(w)=\{\hat{1}_{w}\}$ and thus
$
B_w(a_1, \ldots, a_n)=K_w(a_1, \ldots, a_n),
$
therefore there is nothing to prove. Hence, assume that $w$ is not a constant word
and that the assertion holds for all Motzkin words of length smaller than $|w|$.
Then, by Definition 4.5,
\[
B_w(a_1, \ldots, a_n)=K_w(a_1, \ldots, a_n)+
\sum_{\stackrel{\pi \in \mathcal{NC}_{\rm irr}(w)}{\scriptscriptstyle \pi\neq \hat{1}_{w}}}K_{\pi}[a_1, \ldots, a_n],
\]
where each $\pi$ consists of an outermost block $v_1=s_0s_1\cdots s_m$ and 
subpartitions lying under its bridges, which can be graphically represented by the diagram

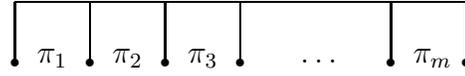
\unitlength=1mm
\special{em:linewidth 1pt}
\linethickness{0.5pt}
\begin{picture}(180.00,20.00)(0.00,15.00)

\put(40.00,28.00){\line(1,0){60.00}}
\put(40.00,20.00){\line(0,1){8.00}}
\put(50.00,20.00){\line(0,1){8.00}}
\put(60.00,20.00){\line(0,1){8.00}}
\put(70.00,20.00){\line(0,1){8.00}}
\put(78.00,20.00){$\ldots$}
\put(90.00,20.00){\line(0,1){8.00}}
\put(100.00,20.00){\line(0,1){8.00}}
\put(43.00,20.00){$\pi_1$}
\put(53.00,20.00){$\pi_2$}
\put(63.00,20.00){$\pi_3$}
\put(93.00,20.00){$\pi_m$}
\put(40.00,20.00){\circle*{1.00}}
\put(50.00,20.00){\circle*{1.00}}
\put(60.00,20.00){\circle*{1.00}}
\put(70.00,20.00){\circle*{1.00}}
\put(90.00,20.00){\circle*{1.00}}
\put(100.00,20.00){\circle*{1.00}}

\end{picture}
\noindent
in which $\pi_1\in \mathcal{NC}(w_1), \ldots, \pi_m\in \mathcal{NC}(w_m)$ and $w_1, \ldots , w_m\in \mathpzc{AM}$, according to the decomposition
\[
w=s_0w_1s_1\cdots w_{m}s_m,
\] 
where at least one subword $w_j$ is non-empty.

Now, we use the inductive assumption inside the blocks of the partitioned cumulants
$K_{\pi}$ which occur on the right-hand side of the above formula. Let
$\hat{1}_{w}\neq \pi\in \mathcal{NC}_{{\rm irr}}(w)$, and let $v$ be a block of
$\pi$. Since $v$ is a Motzkin subword of $w$ of length smaller than $|w|$, the
induction hypothesis gives the identity
\[
K_v
=
\sum_{\rho\in\mathcal{NC}_{{\rm irr}}(v)}
(-1)^{|\rho|-1}B_{\rho}
\]
for the corresponding multilinear functionals. In this sum, the term corresponding
to the one-block partition $\rho=\hat{1}_{v}$ is
\[
(-1)^{|\hat{1}_{v}|-1}B_{\hat{1}_{v}}
=
B_v.
\]
Consequently, when the one-block partition $\hat{1}_{v}$ is chosen in every block
$v$ of $\pi$, the resulting contribution to $K_{\pi}$ is precisely
$B_{\pi}$, while all other choices produce proper refinements of $\pi$.

We claim that, after substituting these expansions into the preceding formula for
$B_w$ and solving for $K_w$, we obtain a formula of the form
\[
K_w(a_1, \ldots , a_n)=B_w(a_1, \ldots, a_n)+
\sum_{\substack{\pi \in \mathcal{NC}_{{\rm irr}}(w)\\ \pi\neq \hat{1}_{w}}}
\mu(\pi)
B_{\pi}[a_1, \ldots, a_n],
\]
where each $\mu(\pi)$ is a constant. Indeed, by the preceding discussion, for each
$\hat{1}_{w}\neq\pi\in\mathcal{NC}_{{\rm irr}}(w)$ we have
\[
K_{\pi}[a_1, \ldots , a_n]
=
B_{\pi}[a_1, \ldots, a_n]
+
\sum_{\sigma \prec \pi}\mu(\pi, \sigma)
B_{\sigma}[a_1, \ldots, a_n]
\]
for some constants $\mu(\pi,\sigma)$. Hence each
$\pi\in \mathcal{NC}_{{\rm irr}}(w)$ may occur in the computations, either directly
from the one-block contributions inside the blocks, or as a refinement produced by
the inductive expansions.

However, we have to justify the fact that no irreducible noncrossing partitions
which are not adapted to $w$ are obtained. This is not {\it a priori} clear
since a refinement of $\pi\in \mathcal{NC}_{{\rm irr}}(w)$ does not have to be
adapted to $w$.

Let us first observe what happens when the induction hypothesis is applied
to the blocks of a fixed partition
$\hat{1}_{w}\neq \pi\in \mathcal{NC}_{{\rm irr}}(w)$. Denote the blocks of
$\pi$ by $v_1,\ldots,v_r$, where $v_1$ is the outermost block. If we express
$K_{v_1}$ in terms of $B_{\rho_1}$, with
$\rho_1\in\mathcal{NC}_{{\rm irr}}(v_1)$, then the blocks of $\rho_1$
satisfy the conditions of Definition 3.1, since $v_1$ is not an inner block
and their heights in $\rho_1$ agree with their heights in the resulting
refinement of $\pi$.

However, if $v_j$ is an inner block of $\pi$ and we express $K_{v_j}$ in
terms of $B_{\rho_j}$, with
$\rho_j\in\mathcal{NC}_{{\rm irr}}(v_j)$, then the heights of inner blocks of
$\rho_j$ may change when $\rho_j$ is inserted into $\pi$. For that reason,
the vertical refinement of $\pi$ induced by $\rho_j$ may have blocks whose
heights are not strictly monotonically increasing with increasing depth, and
therefore may fail to satisfy the conditions of Definition 3.1. Two types of
refinements which should not be allowed have to be analyzed.

{\it Refinement of type 1.} The first type of refinement arises when we have a number of neighboring inner blocks of height $j$ in $\pi$ 
with a bridge $jj$ and a decomposition of their nearest outer block produces a subblock 
which is the nearest outer block of these inner blocks and also has height $j$ (in the simplest case, it is the bridge $jj$). 
Therefore, we obtain a partition $\sigma$ which does not satisfy condition (1) 
of Definition 3.1. Symbolically, we can represent this operation by the diagrams shown below, where squares represent
inner blocks and their nearest outer block is of the form $ujju'$, with the bridge $jj$. The arrows indicate that 
the bridge is separated from the outer block.

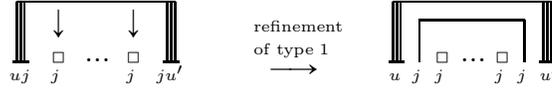
\unitlength=1mm
\special{em:linewidth 1pt}
\linethickness{0.5pt}
\begin{picture}(180.00,20.00)(0.00,13.00)

\put(39.50,28.00){\line(1,0){21.00}}

\put(40.00,20.00){\line(0,1){8.00}}
\put(60.00,20.00){\line(0,1){8.00}}
\put(40.50,20.00){\line(0,1){8.00}}
\put(59.50,20.00){\line(0,1){8.00}}
\put(39.50,20.00){\line(0,1){8.00}}
\put(60.50,20.00){\line(0,1){8.00}}

\put(44.00,19.50){$\square$}
\put(44.00,24.00){$\downarrow$}
\put(48.50,20.00){$...$}
\put(54.00,19.50){$\square$}
\put(54.00,24.00){$\downarrow$}

\put(59.00,20.00){\line(1,0){2.00}}
\put(39.00,20.00){\line(1,0){2.00}}
\put(58.00,17.50){$\scriptscriptstyle{ju'}$}
\put(38.50,17.50){$\scriptscriptstyle{uj}$}
\put(44.00,17.50){$\scriptscriptstyle{j}$}
\put(54.00,17.50){$\scriptscriptstyle{j}$}

\put(71.00,24.00){$\scriptscriptstyle{\rm refinement}$}
\put(73.00,18.00){$\longrightarrow$}
\put(71.00,21.00){$\scriptscriptstyle{\rm of\;type\;1}$}

\put(89.50,28.00){\line(1,0){21.00}}
\put(93.00,25.50){\line(1,0){14.00}}
\put(89.50,20.00){\line(0,1){8.00}}
\put(90.00,20.00){\line(0,1){8.00}}
\put(90.50,20.00){\line(0,1){8.00}}
\put(93.00,20.00){\line(0,1){5.50}}
\put(107.00,20.00){\line(0,1){5.50}}
\put(109.50,20.00){\line(0,1){8.00}}
\put(110.00,20.00){\line(0,1){8.00}}
\put(110.50,20.00){\line(0,1){8.00}}

\put(95.00,19.50){$\square$}
\put(98.50,20.00){$...$}
\put(103.00,19.50){$\square$}

\put(89.00,20.00){\line(1,0){2.00}}
\put(109.00,20.00){\line(1,0){2.00}}

\put(109.00,17.50){$\scriptscriptstyle{u'}$}
\put(89.00,17.50){$\scriptscriptstyle{u}$}
\put(92.00,17.50){$\scriptscriptstyle{j}$}
\put(106.00,17.50){$\scriptscriptstyle{j}$}
\put(95.00,17.50){$\scriptscriptstyle{j}$}
\put(103.00,17.50){$\scriptscriptstyle{j}$}
\end{picture}
\noindent
More generally, a block of the form $j^p$ is separated from the outer block instead of the bridge $jj$.
However, the Motzkin cumulants corresponding to the original inner blocks are associated 
with words of height $j$, say $w_1, \ldots , w_p$, and therefore, when inserted among the arguments of the Motzkin cumulant
corresponding to the new outermost block, they produce
\[
B_{w_0}(a_{i}, \ldots, a_{k}B_{w_1}(\ldots)\cdots B_{w_p}(\ldots), a_{r}, \ldots, a_{m})
\]
for some $a_{i},\ldots, a_k, \ldots, a_r, \ldots , a_{m}$, where $h(w_0)=j$, but this cumulant vanishes by property (1)
of Definition 5.2 (keep $B_{w_1}(\ldots)$ attached to $a_k$ and attach the remaining cumulants to $a_r$).
Therefore, the partitioned cumulant $B_{\sigma}$ does not appear in the decomposition due to property (1) of Definition 5.2.

{\it Refinement of type 2.}
The second type of refinement is obtained when we have a sequence of inner blocks of height $j+1$ in $\pi$ 
and their nearest outer block is decomposed into subblocks in such a way that two or more 
blocks of height $j$ are produced next to these inner blocks at the same depth. 
Symbolically, this type of refinement is shown below in terms of diagrams. The arrows show that the bridge is separated
from the outer block and splits into two subblocks.

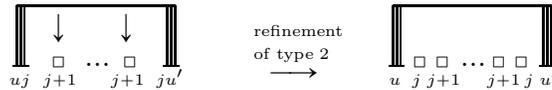
\unitlength=1mm
\special{em:linewidth 1pt}
\linethickness{0.5pt}
\begin{picture}(180.00,20.00)(0.00,13.00)

\put(39.50,28.00){\line(1,0){21.00}}

\put(40.00,20.00){\line(0,1){8.00}}
\put(60.00,20.00){\line(0,1){8.00}}
\put(40.50,20.00){\line(0,1){8.00}}
\put(59.50,20.00){\line(0,1){8.00}}
\put(39.50,20.00){\line(0,1){8.00}}
\put(60.50,20.00){\line(0,1){8.00}}

\put(44.00,19.50){$\square$}
\put(44.00,24.00){$\downarrow$}
\put(48.50,20.00){$...$}
\put(53.00,24.00){$\downarrow$}
\put(53.50,19.50){$\square$}

\put(59.00,20.00){\line(1,0){2.00}}
\put(39.00,20.00){\line(1,0){2.00}}
\put(58.00,17.50){$\scriptscriptstyle{ju'}$}
\put(38.50,17.50){$\scriptscriptstyle{uj}$}
\put(43.00,17.50){$\scriptscriptstyle{j+1}$}
\put(52.00,17.50){$\scriptscriptstyle{j+1}$}

\put(71.00,24.00){$\scriptscriptstyle{\rm refinement}$}
\put(73.00,18.00){$\longrightarrow$}
\put(71.00,21.00){$\scriptscriptstyle{\rm of\;type\;2}$}

\put(89.50,28.00){\line(1,0){21.00}}
\put(89.50,20.00){\line(0,1){8.00}}
\put(90.00,20.00){\line(0,1){8.00}}
\put(90.50,20.00){\line(0,1){8.00}}
\put(109.50,20.00){\line(0,1){8.00}}
\put(110.00,20.00){\line(0,1){8.00}}
\put(110.50,20.00){\line(0,1){8.00}}

\put(95.00,19.50){$\square$}
\put(98.50,20.00){$...$}
\put(102.50,19.50){$\square$}
\put(92.00,19.50){$\square$}
\put(105.50,19.50){$\square$}

\put(89.00,20.00){\line(1,0){2.00}}
\put(109.00,20.00){\line(1,0){2.00}}

\put(109.00,17.50){$\scriptscriptstyle{u'}$}
\put(89.00,17.50){$\scriptscriptstyle{u}$}
\put(92.00,17.50){$\scriptscriptstyle{j}$}
\put(107.00,17.50){$\scriptscriptstyle{j}$}
\put(94.00,17.50){$\scriptscriptstyle{j+1}$}
\put(102.00,17.50){$\scriptscriptstyle{j+1}$}
\end{picture}
\noindent
Then, the cumulants corresponding to the original inner blocks are associated with words of height $j+1$, say $w_{1}, \ldots , w_{p}$, whereas those obtained from the decomposition of the outer block are associated with words, say $w_{0}, w_{p+1}$, of height $j$. 
Since these blocks stand next to each other in the refinement $\sigma$, 
we obtain a cumulant of the form 
\[
B_{uu'}(a_{i}, \ldots, a_{k}B_{w_0}(\ldots )\cdots B_{w_{p+1}}(\ldots ), a_{r}, \ldots, a_{m})
\]
and thus it vanishes by the orthogonality of cumulants associated with words of different heights 
(property (2) of Definition 5.2).
Therefore, the partitioned cumulant corresponding to $B_{\sigma}$ does not appear in the decomposition.
Note that the corresponding partition $\sigma$ is not adapted to $w$ since condition (2) of Definition 3.1 is not 
satisfied. These are the only refinements $\sigma$ of $\pi$ induced by the refinements of some $\pi_j$, $j>1$, that
could produce partitions which are not adapted to $w$. This proves our claim that only partitions adapted to $w$ appear in the decomposition.

We now compute the coefficients $\mu(\pi)$.
We fix $\pi\in \mathcal{NC}_{{\rm irr}}(w)$ in the formula for $K_{w}(a_1, \ldots, a_n)$ 
in order to compute $\mu(\pi)$.
A contribution to $\mu(\pi)$ comes from all partitions $\pi'\succeq \pi$ in $\mathcal{NC}_{{\rm irr}}(w)$ 
which are ancestors of $\pi$, since by induction each $K_{v}$ for any block
$v$ of $\pi'$ is expressed in terms of $B_{\pi_{v}}$, where $\pi_{v}\in \mathcal{NC}_{{\rm irr}}(v)$, and therefore a partition in which a block 
is obtained by a horizontal merging of blocks of $\pi$ (at any depth) will not give $B_{\pi}$ in the formula for 
$K_{\pi'}$ due to the fact that only irreducible partitions are used in the decomposition formulas for $K_{v}$. 
Symbolically, when we restrict our attention to the (locally) deepest inner blocks of $\pi$ and its outer block,
the ancestors of $\pi$ are obtained by lifting up some of these blocks, which is shown in the diagram given below.

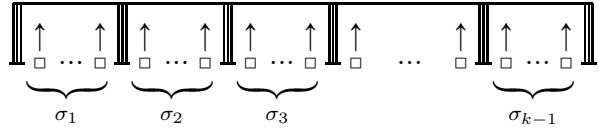
\unitlength=1mm
\special{em:linewidth 1pt}
\linethickness{0.5pt}
\begin{picture}(180.00,30.00)(5.00,7.00)

\put(39.50,28.00){\line(1,0){77.00}}

\put(40.00,20.00){\line(0,1){8.00}}
\put(54.00,20.00){\line(0,1){8.00}}
\put(68.00,20.00){\line(0,1){8.00}}
\put(82.00,20.00){\line(0,1){8.00}}
\put(102.00,20.00){\line(0,1){8.00}}
\put(116.00,20.00){\line(0,1){8.00}}

\put(40.50,20.00){\line(0,1){8.00}}
\put(54.50,20.00){\line(0,1){8.00}}
\put(68.50,20.00){\line(0,1){8.00}}
\put(82.50,20.00){\line(0,1){8.00}}
\put(102.50,20.00){\line(0,1){8.00}}
\put(116.50,20.00){\line(0,1){8.00}}

\put(39.50,20.00){\line(0,1){8.00}}
\put(53.50,20.00){\line(0,1){8.00}}
\put(67.50,20.00){\line(0,1){8.00}}
\put(81.50,20.00){\line(0,1){8.00}}
\put(101.50,20.00){\line(0,1){8.00}}
\put(115.50,20.00){\line(0,1){8.00}}

\put(42.00,19.00){$\square$}
\put(42.00,22.50){$\uparrow$}
\put(45.50,20.00){$...$}
\put(50.00,19.00){$\square$}
\put(50.00,22.50){$\uparrow$}
\put(56.00,19.00){$\square$}
\put(56.00,22.50){$\uparrow$}
\put(59.50,20.00){$...$}
\put(64.00,19.00){$\square$}
\put(64.00,22.50){$\uparrow$}
\put(70.00,19.00){$\square$}
\put(70.00,22.50){$\uparrow$}
\put(73.50,20.00){$...$}
\put(78.00,19.00){$\square$}
\put(78.00,22.50){$\uparrow$}

\put(84.00,19.00){$\square$}
\put(84.00,22.50){$\uparrow$}
\put(90.50,20.00){$...$}
\put(98.00,19.00){$\square$}
\put(98.00,22.50){$\uparrow$}

\put(104.00,19.00){$\square$}
\put(104.00,22.50){$\uparrow$}
\put(107.50,20.00){$...$}
\put(112.00,19.00){$\square$}
\put(112.00,22.50){$\uparrow$}

\put(39.00,20.00){\line(1,0){2.00}}
\put(53.00,20.00){\line(1,0){2.00}}
\put(67.00,20.00){\line(1,0){2.00}}
\put(81.00,20.00){\line(1,0){2.00}}
\put(101.00,20.00){\line(1,0){2.00}}
\put(115.00,20.00){\line(1,0){2.00}}

\put(41.00,18.00){$\underbrace{}_{\sigma_1}$}
\put(55.00,18.00){$\underbrace{}_{\sigma_2}$}
\put(69.00,18.00){$\underbrace{}_{\sigma_3}$}
\put(103.00,18.00){$\underbrace{}_{\sigma_{k-1}}$}

\end{picture}
\noindent
Here, multiple legs correspond to subwords $u_1, \ldots, u_k$ which make up the outer block and 
squares are the inner blocks which make up 
$\sigma_1\in \mathcal{NC}(w_1), \ldots, \sigma_{k-1}\in \mathcal{NC}(w_{k-1})$, respectively. 
Thus, the considered subword of $w$ is of the form
\[
u_1w_1u_2\cdots w_{k-1}u_k,
\]
where each $w_j$ is a juxtaposition of subwords corresponding to the inner blocks lying between $u_i$ and $u_{i+1}$.
Then the ancestors of this part of $\pi$ are obtained by connecting some of the inner blocks with their nearest outer block. 
These operations are denoted by arrows. 
Each arrow indicates that a given block can be connected with its nearest outer block. 
The same lifting procedure can be applied to all blocks and their nearest outer blocks at all depths.
The only restriction is that not all inner blocks can be lifted simultaneously since in that case we obtain 
$\hat{1}_{w}$ which corresponds to $K_w$.
In this fashion, apart from $\pi$ itself, which means no liftings, we obtain the family of all 
ancestors of $\pi$ except $\hat{1}_{w}$. 
The decomposition of each $K_{v}$ in terms of 
$B_{\pi_{v}}$, where $v\in \pi'\in \mathcal{A}(\pi)\setminus \{\hat{1}_{w}\}$ and $\pi_{v}\in \mathcal{NC}_{{\rm irr}}(v)$, 
gives the coefficient $(-1)^{|\pi_{v}|-1}$. 
When using the decomposition of $K_v$, we use the bimodule property of the
cumulants $B_{v'}$, $v'\in\pi_v$, whenever an expression of the form
$B_{v'}(a_kb)$ occurs, namely
$
B_{v'}(a_kb)=B_{v'}(a_k)b
$
for $b\in\mathcal{B}$.

Thus, since the irreducible partitions $\pi_v$, $v\in\pi'$, together form
the partition $\pi$, we have
\[
\sum_{v\in\pi'}|\pi_v|=|\pi|.
\]
Therefore the product over all blocks $v$ of $\pi'$ contributes to $\mu(\pi)$
the coefficient
\[
\prod_{v\in \pi'}(-1)^{|\pi_v|-1}
=
(-1)^{|\pi|-|\pi'|}.
\]
Consequently,
\[
\mu(\pi)
=
-
\sum_{\substack{\pi'\in\mathcal{A}(\pi)\\ \pi'\neq \hat{1}_{w}}}
(-1)^{|\pi|-|\pi'|}.
\]
Now, each $\pi'\in\mathcal{A}(\pi)$ is obtained from $\pi$ by lifting up
$k$ inner blocks of $\pi$, where $k=0,\ldots,i(\pi)-1$ and $i(\pi)$ denotes
the number of inner blocks of $\pi$. The case $k=0$ gives $\pi'=\pi$, whereas
the case $k=i(\pi)$ gives $\pi'=\hat{1}_{w}$ and is excluded. Since
$|\pi|=i(\pi)+1$, we obtain
\[
\mu(\pi)
=
-
\sum_{k=0}^{i(\pi)-1}
(-1)^k {i(\pi)\choose k}
=
(-1)^{i(\pi)}
=
(-1)^{|\pi|-1},
\]
which completes the proof.
\hfill $\boxvoid$\\

\begin{Example}
{\rm 
Let us give simple examples of refinements of the type used in the proof of Theorem 5.3.
In both cases, we consider $K_{\pi'}$, where $\pi'$ is an ancestor of $\pi$, and show that it contributes a scalar multiple of $B_{\pi}$ to the expression for $K_{w}$.
For $w=12321$, a refinement is given by the pair

\unitlength=1mm
\special{em:linewidth 1pt}
\linethickness{0.5pt}
\begin{picture}(180.00,15.00)(20.00,15.00)

\put(62.00,20.50){$\pi'=$}
\put(72.00,20.00){\circle*{1.00}}
\put(75.00,20.00){\circle*{1.25}}
\put(78.00,20.00){\circle*{1.60}}
\put(81.00,20.00){\circle*{1.25}}
\put(84.00,20.00){\circle*{1.00}}
\put(72.00,20.00){\line(0,1){4.5}}
\put(75.00,20.00){\line(0,1){4.5}}
\put(81.00,20.00){\line(0,1){4.5}}
\put(84.00,20.00){\line(0,1){4.5}}
\put(72.00,24.50){\line(1,0){12}}

\put(90.00,20.50){${\rm and}$}

\put(102.00,20.50){$\pi=$}
\put(112.00,20.00){\circle*{1.00}}
\put(115.00,20.00){\circle*{1.25}}
\put(118.00,20.00){\circle*{1.60}}
\put(121.00,20.00){\circle*{1.25}}
\put(124.00,20.00){\circle*{1.00}}
\put(112.00,20.00){\line(0,1){4.5}}
\put(124.00,20.00){\line(0,1){4.5}}
\put(115.00,20.00){\line(0,1){2.5}}
\put(121.00,20.00){\line(0,1){2.5}}
\put(112.00,24.50){\line(1,0){12}}
\put(115.00,22.50){\line(1,0){6}}
\put(126.00,20.00){,}
\end{picture}
\noindent
for which
\begin{eqnarray*}
K_{\pi'}[a_1, \ldots, a_5]&=& K_{1221}(a_1,a_2K_{3}(a_3), a_4, a_5)\\
&=&B_{1221}(a_1,a_2B_{3}(a_3), a_4, a_5)\\
&&\!\!\!\!\!\!-B_{11}(a_1,B_{22}(a_2B_{3}(a_3),a_4), a_5)
\end{eqnarray*}
where we used the expression for $K_{1221}$ from Example 4.8 (which has five terms) and Definition 5.2 (which makes 
three terms vanish).
This shows that $K_{\pi'}$ produces $B_{\pi}$ with the coefficient $\mu(\pi',\pi)=-1$.

For $w=12^31$, a refinement given by the pair

\unitlength=1mm
\special{em:linewidth 1pt}
\linethickness{0.5pt}
\begin{picture}(180.00,15.00)(20.00,15.00)

\put(62.00,20.50){$\pi'=$}
\put(72.00,20.00){\circle*{1.00}}
\put(75.00,20.00){\circle*{1.25}}
\put(78.00,20.00){\circle*{1.25}}
\put(81.00,20.00){\circle*{1.25}}
\put(84.00,20.00){\circle*{1.00}}
\put(72.00,20.00){\line(0,1){4.5}}
\put(75.00,20.00){\line(0,1){4.5}}
\put(84.00,20.00){\line(0,1){4.5}}
\put(72.00,24.50){\line(1,0){12}}

\put(90.00,20.50){${\rm and}$}

\put(102.00,20.50){$\pi=$}
\put(112.00,20.00){\circle*{1.00}}
\put(115.00,20.00){\circle*{1.25}}
\put(118.00,20.00){\circle*{1.25}}
\put(121.00,20.00){\circle*{1.25}}
\put(124.00,20.00){\circle*{1.00}}
\put(112.00,20.00){\line(0,1){4.5}}
\put(124.00,20.00){\line(0,1){4.5}}
\put(112.00,24.50){\line(1,0){12}}

\end{picture}
\noindent
gives 
\begin{eqnarray*}
K_{\pi'}[a_1, \ldots, a_5]&=& K_{121}(a_1,a_2K_{2}(a_3)K_2(a_4), a_5)\\
&=&B_{121}(a_1,a_2B_{2}(a_3)B_2(a_4), a_5)\\
&&\!\!\!\!\!\!-\; B_{11}(a_1B_{2}(a_2)B_{2}(a_3)B_{2}(a_4), a_5),
\end{eqnarray*}
where we used the ${\mathcal B}$-bimodule property $B_2(a_2b)=B_2(a_2)b$ for $b=B_2(a_3)B_2(a_4)$ (see Remark 4.6(5)). 
In particular,
$K_{\pi'}$ produces $B_{\pi}$ with coefficient $\mu(\pi',\pi)=-1$. 
}
\end{Example}

\section{Orthogonal replicas}

Our approach to cumulants has its roots in the decomposition of free random variables 
in terms of orthogonal replicas \cite{[L3],[L6]}. Although it is more general, 
we would like to apply it now to free probability. 

The idea is to realize free random variables as sums of orthogonal replicas indexed by tensor coordinates. This allows us to encode combinatorial information (such as non-crossing structures) via projections acting on tensor positions. The constructions introduced below will later be used to define conditional expectations and cumulants.

Let us recall that for each variable from two given noncommutative probability spaces 
we can construct a countable family of its copies, called its orthogonal replicas, as elements of an 
infinite tensor product of algebras, and define a suitable normalized linear functional on the unital algebra generated 
by these replicas. 

\begin{Definition}
{\rm Let $p$ be an abstract projection, namely $p^2=p$.
For any noncommutative probability space
$({\mathcal A}, \varphi)$, we construct its {\it Boolean extension} (or, $p$-{\it extension})
$(\widetilde{\mathcal A}, \widetilde{\varphi})$, where 
\[
\widetilde{\mathcal A}={\mathcal A}*{\mathbb C}[p]
\] 
is the free product with identified units and $\widetilde{\varphi}$ is the linear extension of 
\[
\widetilde{\varphi}(p^{\alpha}a_1pa_2p\cdots pa_mp^{\gamma})
=
\varphi(a_1)\varphi(a_2)\cdots \varphi(a_m)
\]  
and $\widetilde{\varphi}(p)=1$, where $\alpha, \gamma\in \{0,1\}$ and 
$a_1,\ldots , a_m\in {\mathcal A}$. Clearly, the distributions
of any $a\in {\mathcal A}$ with respect to the functionals $\widetilde{\varphi}$ and $\varphi$ are identical.
}
\end{Definition}

This concept was introduced in \cite{[L1]} and proved to be very useful in
studying independence and various graph products (see, for instance,
\cite{[ALS]}). We now introduce orthogonal replicas, using two projections
$p$ and $q$ to distinguish the two algebras and enforce orthogonality between
tensor components.

\begin{Definition}
{\rm 
Let $({\mathcal A}_{1},\varphi_{1})$ and $({\mathcal A}_{2}, \varphi_{2})$ be noncommutative probability spaces and
let $(\widetilde{\mathcal A}_{1}, \widetilde{\varphi}_{1})$ and $(\widetilde{\mathcal A}_{2}, \widetilde{\varphi}_{2})$ 
be their $p$-extension and $q$-extension, respectively. Define infinite tensor products
\[
\Phi_{\otimes}:=\Phi_1\otimes \Phi_2\;\;\;{\rm and}\;\;\; {\mathcal A}_{\otimes}:=
\widetilde{\mathcal A}_{1}^{\otimes \infty}\otimes \widetilde{\mathcal A}_{2}^{\otimes\infty},
\]
where $\Phi_i=\widetilde{\varphi}_{i}^{\otimes\infty}$ for $i=1,2$. Let $a\in \mathcal{A}_{1}$, $a'\in \mathcal{A}_2$,
and let $1_k$ be the unit in ${\mathcal A}_{k}$ for $k=1,2$. The elements of ${\mathcal A}_{\otimes}$ of the form\\
\begin{eqnarray*}
a(j)&=&(\cdots \otimes 1_1 \otimes \underbrace{a}_{{\rm site}\;j} \otimes 1_1\otimes \cdots )\otimes 
(\cdots \otimes 1_2\otimes \underbrace{q^{\perp}}_{{\rm site}\;j-1} \otimes\; q \otimes \cdots),\\
a'(j)&=&(\cdots \otimes 1_1 \otimes \underbrace{p^{\perp}}_{{\rm site}\;j-1} \otimes p\otimes \cdots )\otimes 
(\cdots \otimes 1_2\otimes \underbrace{a'}_{{\rm site}\;j} \otimes\; 1_2 \otimes \cdots)
\end{eqnarray*}
will be called {\it orthogonal replicas} of 
$a$ and $a'$ of {\it color} $j\in \mathbb{N}$ and {\it labels} $1$ and $2$, respectively. 
By convention, when $j=1$, then the terms involving site $j-1$ are omitted. Thus each replica occupies a single tensor coordinate, while the surrounding projections ensure that replicas of different colors are orthogonal with respect to the product state.
The noncommutative probability space $({\mathcal A}_{{\rm rep}}, \Phi)$, where
${\mathcal A}_{{\rm rep}}$ is the unital subalgebra of ${\mathcal A}_{\otimes}$ generated by orthogonal replicas 
and $\Phi$ is the restriction of $\Phi_{\otimes}$ to ${\mathcal A}_{{\rm rep}}$, will be called the {\it orthogonal replica space}.
}
\end{Definition}

\begin{Remark}
{\rm We now introduce certain sequences of projections that will play an important role in the sequel.
\begin{enumerate}
\item
Partial sums 
\begin{eqnarray*}
e_{1,n}:&=&\sum_{j=1}^{n}1_{1}(j)=1_{1}^{\otimes \infty}\otimes (1_{2}^{\otimes (n-1)}\otimes q^{\otimes \infty})\\
e_{2,n}:&=&\sum_{j=1}^{n}1_{2}(j)=(1_{1}^{\otimes (n-1)}\otimes p^{\otimes \infty})\otimes 1_{2}^{\otimes \infty}
\end{eqnarray*}
are `sequential approximate identities' dedicated to labels $1$ and $2$ in the sense that they 
commute with the associated orthogonal replicas of all colors. 
\item 
The sequence of {\it compressed} projections $(e_{n})_{n\in \mathbb N}$ given by products
\[
e_{n}:=e_{1,n}e_{2,n}=(1_{1}^{\otimes (n-1)}\otimes p^{\otimes \infty})\otimes
(1_{2}^{\otimes (n-1)}\otimes q^{\otimes \infty})
\]
is a `symmetric sequential approximate identity' since for 
any $a,a'\in {\mathcal A}_{{\rm rep}}$ there exists $m$ such that $\Phi(ae_na')=\Phi(aa')$ for any $n\geq m$.
\item
The associated sequence $(p_n)_{n\in \mathbb{N}}$ of orthogonal projections is given by
\[
p_n=e_{n}-e_{n-1},
\]
where $n\geq 1$ and we set $e_{0}=0$. We will say that $p_n$ is an 
{\it orthogonal projection of color $n$}. 
In particular, $p_1=e_1$. 
\item
We also introduce the commutative algebra
\[
{\mathcal B}:={\mathbb C}[p_1,p_2, \ldots ]={\mathbb C}[e_1,e_2, \ldots]\subset {\mathcal A}_{{\rm rep}}.
\]
\item 
We will also need a more general family of projections
\[
\mathcal{F}=\{e(\epsilon,\eta):\epsilon,\eta\in\{0,1\}^{\infty}\},
\]
where
\[
e(\epsilon,\eta)
:=
(p^{\epsilon_1}\otimes p^{\epsilon_2}\otimes \cdots)
\otimes
(q^{\eta_1}\otimes q^{\eta_2}\otimes \cdots).
\]
In particular, if
\[
\epsilon=(\underbrace{0, \ldots, 0}_{n-1\;\text{times}}, 1,1, \ldots )=\eta,
\] 
then $e(\epsilon, \eta)=e_n$.
\item
Let \(S\) denote the set of elementary tensors
\[
a=
\left(\bigotimes_{r=1}^{\infty} x_r\right)
\otimes
\left(\bigotimes_{r=1}^{\infty} y_r\right)
\in
\widetilde{\mathcal A}_{1}^{\otimes\infty}
\otimes
\widetilde{\mathcal A}_{2}^{\otimes\infty},
\]
where all but finitely many \(x_r\)'s and \(y_r\)'s are equal to the corresponding
units, and each non-unit \(x_r\), respectively \(y_r\), is a reduced word in
\[
\widetilde{\mathcal A}_{1}=\mathcal A_1 * \mathbb C[p],
\qquad
\widetilde{\mathcal A}_{2}=\mathcal A_2 * \mathbb C[q],
\]
which neither begins nor ends with \(p\), respectively neither begins nor ends
with \(q\). Thus \(p\)'s and \(q\)'s may occur inside tensor components, but not
at their boundaries. 
Then \({\mathcal A}_{{\rm rep}}\) is spanned by elements of the form
$eaf$, where $e,f\in\mathcal F$ and $a\in S$.
\item
If \(e=e(\epsilon,\eta)\in \mathcal F\) is not the identity, then we set
\[
j(e):=\min\{j:\epsilon_j=1\;\;{\rm or}\;\;\eta_j=1\}.
\]
For the identity projection we put \(j(1)=\infty\), with the convention
\(e_{\infty}=1\). Thus \(j(e)\) is the smallest color occupied by a nontrivial
projection \(p\) or \(q\). In particular,
$j(e_n)=j(e_{1,n})=j(e_{2,n})=n$.
\item
Finally, let us mention a Hilbert space interpretation. 
If $({\mathcal A}_{1}, \varphi_1)$, $({\mathcal A}_{2}, \varphi_2)$ are 
$C^*$-probability spaces and $({\mathcal H}_{1}, \pi_1,\xi_1)$, $({\mathcal H}_{2},\pi_2, \xi_2)$  
are their GNS triples, one can use the free product
${\mathcal H}_{1}\star {\mathcal H}_{2}$ and interpret the replica 
$a_k(j)$ under $\Phi$ as a restriction of 
the left regular representation $\lambda(a_k)$ to the tensor product 
\[
{\mathcal H}_{i}\otimes {\mathcal H}_{i_1}^{\circ}\otimes \ldots {\mathcal H}_{i_{j-1}}^{\circ}
\]
where $i$ is the label of $a_k$, $i\neq i_1\neq \ldots \neq i_{j-1}$ and 
${\mathcal H}_{m}^{\circ}={\mathcal H}_{m}\ominus {\mathbb C}\xi_{m}$, $m=1,2$.
The projections $e_n$ become the projections onto 
\[
\bigoplus_{j=0}^{n-1}\bigoplus_{i_1\neq \cdots \neq i_{j}}
{\mathcal H}_{i_1}^{\circ}\otimes \ldots \otimes {\mathcal H}_{i_{j}}^{\circ}
\]
for $n>1$, with $e_1$ being the projection onto ${\mathbb C}\xi$, where $\xi$ is a unit vector.
In the category of noncommutative probability spaces, the Hilbert space framework 
has to be replaced by a more algebraic one, like that based on tensor products.
Then we can also consider functionals other than $\Phi$.
\end{enumerate}
}
\end{Remark}

\begin{Proposition}
Let $a_k\in \mathcal{A}_{i_k}$ and let $j_k\in {\mathbb N}$, where $k\in \{1,2\}$.
\begin{enumerate}
\item
If $i_1=i_2$ and $j_1\neq j_2$, then for any $j$, 
\[
a_{1}(j_1)p_ja_{2}(j_2)=0.
\]
\item
If $i_1\neq i_2$ and $j_1=j_2=j$, then for any $k>j$,
\[
a_1(j_1)p_{k}a_2(j_2)=0.
\]
\item
If $i_1=i_2$ and $j_1=j_2=j$, 
\[
a_1(j_1)a_2(j_2)=(a_1a_2)(j).
\]
\end{enumerate}
\end{Proposition}
{\it Proof.}
These properties follow directly from the definition of orthogonal replicas 
and reflect the way projections enforce orthogonality between different colors.
\hfill $\boxvoid$\\

The rules of computing mixed moments of orthogonal replicas with respect to $\Phi$ 
were given in Lemma 3.1 of \cite{[L6]} and are stated below for the readers' convenience. 
With Proposition 6.4 treating the case when the neighboring labels are equal, these rules could be viewed 
as a definition of `partial freeness' with respect to $\Phi$.

\begin{Lemma}\cite{[L6]}
Let $a_k\in {\mathcal A}_{i_k}$, where $k=1, \ldots , n$
and $i_1\neq \cdots \neq i_n$. Then we have 
\begin{enumerate}
\item
the freeness property:
\[
\Phi\left(a_1(j_1)\cdots a_{n}(j_n)\right)=0
\]
whenever $a_k\in {\rm Ker}(\varphi_{i_k})$ for $k=1, \ldots n$, and 
$j_1, \ldots , j_n\in \mathbb{N}$,
\item
the reduction property:
\[
\Phi\left(a_1(j_1)\cdots 1_{i_r}(j_r) \cdots a_{n}(j_n)\right)
=
\Phi\left(a_1(j_1)\cdots \check{1}_{i_r}(j_r) \cdots a_{n}(j_n)\right)
\]
whenever $a_k\in {\rm Ker}(\varphi_{i_k})$ for
$k=1, \ldots r-1$, and $(j_1, \ldots , j_r)=(1,\ldots , r)$, 
where the element with $\,\check{}\,$ on top has to be omitted, 
and the moment on the LHS vanishes for the remaining colors.
\end{enumerate}
\end{Lemma}

These rules can be interpreted as a form of `partial freeness', 
with $\Phi$ replacing the free product of states. 
The key feature is the reduction property, in which the role of units is played 
by projections that can be removed only in specific mixed moments. 
This is in contrast to the case of free random variables with respect to the free product 
of functionals, where (identified) units can always be deleted, as observed by Avitzour \cite{[Av]}.

Moreover, it was shown in \cite{[L6]} that if the condition 
$a_k\in {\rm Ker}(\varphi_{i_k})$ for $k=1,\ldots,n$ is dropped, 
then mixed moments of orthogonal replicas vanish unless the associated word 
is a reduced Motzkin word. 

We now extend this result to include mixed moments of orthogonal replicas 
interlaced with elements of $\mathcal{B}$. Let $b_0,\ldots,b_n\in\mathcal{B}$. 
A product of the form
\[
b_{0}a_1(j_1)b_{1}\cdots b_{n-1}a_n(j_n)b_{n}
\]
will be called an {\it alternating product associated with} 
$w=j_1\cdots j_n\in \mathpzc{AM}$. 

As for labels, we associate to this product the word 
$\ell=i_1\cdots i_n\in \{1,2\}^{n}$ and impose additional assumptions when needed.

\begin{Proposition}
If $a_{1}, \ldots, a_n\in \mathcal{A}_{1}\cup \mathcal{A}_{2}$ have labels encoded by $\ell=i_1\cdots i_n$, then
\[
\Phi\left(a_1(j_1)\cdots a_n(j_n)\right)=0
\] 
whenever $w:=j_1\cdots j_n\notin \mathpzc{M}_{n}(\ell)$, where  
\[
\mathpzc{M}_{n}(\ell)=\{w\in \mathpzc{M}_{n}: j_k=j_{k+1} \;if\;i_k=i_{k+1}, \; k\in [n-1]\}.
\]
In particular, if $i_1\neq \ldots \neq i_n$, then the above mixed moment vanishes whenever $w\notin \mathpzc{M}_{n}$.
Moreover, the analogous property holds for alternating products associated with $w$.
\end{Proposition}

{\it Proof.}
The case when $b_j=1$ for $j=0, \ldots, n$ was proved in Proposition 5.1 of \cite{[L6]}.
In the case of arbitrary alternating products, we have 
\[
\Phi\left(b_0a_1(j_1)b_1\cdots b_{n-1}a_{n}(j_n)b_n\right)=
\Phi\left(e_1b_0a_1(j_1)b_1\cdots b_{n-1}a_{n}(j_n)b_ne_1\right).
\]
Recall that the proof in \cite{[L6]} is based on the observation that the moment 
on the right-hand side vanishes when $w\notin \mathpzc{M}_{n}(\ell)$, since 
$p$ meets $p^{\perp}$ or $q$ meets $q^{\perp}$ at some tensor site.
Now, in the case of alternating products
it suffices to consider $b_0, \ldots, b_{n}\in \{1, p_1, p_2, \ldots \}$.
However, the projections from this set cannot change the effect of $p$ meeting $p^{\perp}$
or $q$ meeting $q^{\perp}$. Therefore, the alternating moment must also vanish when 
$w\notin \mathpzc{M}_{n}(\ell)$, which completes the proof. 
\hfill $\boxvoid$\\

Finally, let us remark that the definition of Motzkin cumulants was inspired by 
a formula expressing mixed moments of free random variables in terms of the 
mixed moments of their orthogonal replicas (Corollary 5.1 in \cite{[L6]}).
We recall this formula below, using the unit identification map explicitly.
For unital algebras ${\mathcal A}_{1}$ and ${\mathcal A}_{2}$, 
denote by ${\mathcal A}_1*{\mathcal A}_{2}$ and ${\mathcal A}_{1}\star {\mathcal A}_{2}$ their free products
without and with identification of units, respectively.

\begin{Theorem}\cite{[L6]}
Let $({\mathcal A}_{1}, \varphi_1)$ and $({\mathcal A}_{2}, \varphi_2)$ the noncommutative probability spaces 
and let $({\mathcal A}_{1}\star {\mathcal A}_{2}, \varphi_1\star \varphi_2)$ be their free product.
It holds that
\[
((\varphi_1\star \varphi_2)\circ \tau)(a_1\cdots a_n)=
\sum_{w=j_1\cdots j_n\in \mathpzc{M}_{n}}\Phi(a_{1}(j_1)\cdots a_{n}(j_n))
\]
where $\tau: {\mathcal A}_1*{\mathcal A}_{2}\rightarrow {\mathcal A}_{1}\star {\mathcal A}_{2}$ 
is the unit identification map and $a_1(j_1), \ldots, a_n(j_n)$ are the orthogonal 
replicas of $a_1\in \mathcal{A}_{i_1}, \ldots, a_n\in \mathcal{A}_{i_n}$ of colors $j_1, \ldots, j_n$, respectively, 
where $i_1\neq \cdots \neq i_n$ and $n\in \mathbb{N}$.
\end{Theorem}

This formula shows that mixed moments of free random variables decompose 
into contributions indexed by Motzkin words, which provides the basis 
for the definition of Motzkin cumulants introduced in the next section.

\begin{Remark}
{\rm 
The construction of the replica space can be easily generalized to reproduce the mixed moments of conditionally free random variables. It suffices to take 
\[
\Phi_1=\widetilde{\varphi}_{1}\otimes \widetilde{\psi}_{1}^{\otimes \infty}\;\;\;{\rm and}\;\;\; 
\Phi_2=\widetilde{\varphi}_{2}\otimes \widetilde{\psi}_{2}^{\otimes \infty},
\]
where a pair of normalized linear functionals for each label, say 
$(\varphi_1, \psi_1)$ and $(\varphi_2,\psi_2)$, is used to construct the tensor product functionals $\Phi_1$ and $\Phi_2$. 
Using the language of colors, we associate color $1$ with $\varphi_1$ and $\varphi_2$, 
whereas the remaining colors are associated with $\psi_1$ and $\psi_2$ 
(see \cite{[L1],[L3]}). More generally, we can assume that 
$\Phi_1=\bigotimes_{j=1}^{\infty}\widetilde{\varphi}_{1,j}$ and $\Phi_2=\bigotimes_{j=1}^{\infty}\widetilde{\varphi}_{2,j}$,
which gives the model of `freeness with infinitely many functionals' \cite{[CDI]}.
}
\end{Remark}

\section{Conditional expectation}

\subheading{The \({\mathcal B}\)-valued expectation \(E\)}

Let us define a \({\mathcal B}\)-valued expectation induced by \(\Phi\) on the
orthogonal replica space. Its purpose is to retain, in addition to the scalar
moments computed by \(\Phi\), the information about the heights of the Motzkin
subpaths associated with blocks.

To achieve this, we use the sequence \((e_n)\). The definition of \(E\) is based
on computing moments with respect to \(\Phi\), but requires a suitable treatment
of the projections \(p,q\) appearing at the boundaries of elementary tensors.

We shall use the following convention. Whenever an elementary tensor from the
spanning set described in Remark 6.3 is written in the form
\[
x=eaf,\qquad e,f\in\mathcal F,\quad a\in S,
\]
we always mean the form obtained by moving all boundary projections \(p,q\) of the
tensor components into the exterior factors \(e\) and \(f\). More precisely, left
boundary projections are put into \(e\), and right boundary projections are put
into \(f\). If, at some tensor site, the remaining middle component is the unit and
the whole tensor component is equal to \(p\) or \(q\), then this projection is put
both into the left and into the right exterior factor; thus, locally, we use
\[
p=p\cdot 1\cdot p,\qquad q=q\cdot 1\cdot q.
\]
With this convention the factors \(e\), \(a\), and \(f\) are uniquely determined
for each elementary tensor.

\begin{Definition}
{\rm Let $x=eaf$, where $e,f\in\mathcal F$ and $a\in S$, be an elementary tensor written 
according to the above convention. We set
\[
E(x):=e_{j(e)}\Phi(a)e_{j(f)},
\]
where $j(e), j(f)$ are defined as in Remark 6.3. 
The map \(E:{\mathcal A}_{\rm rep}\to{\mathcal B}\) is obtained by extending this
formula linearly.}
\end{Definition}

Thus \(E\) may be viewed as a \(\mathcal B\)-valued lift of the scalar functional
\(\Phi\): it computes the scalar moment of the middle tensor and records, in
\(\mathcal B\), the boundary levels determined by the projections \(e_n\).

\begin{Proposition}
The mapping \(E\) is a conditional expectation. In particular,
\[
E(e_n)=e_n
\qquad\text{and}\qquad
E(e_nxe_m)=e_nE(x)e_m
\]
for any \(n,m\) and any \(x\in {\mathcal A}_{\rm rep}\).
\end{Proposition}

{\it Proof.}
It is clear that $E(e_n)=e_n$. Indeed, 
we may write $e_n=eaf$, where $a$ is the identity tensor and $e=f=e(\epsilon,\eta)$, with
\[
\epsilon=(\underbrace{0,\ldots,0}_{n-1\;\text{\rm times}},1,1,\ldots)=\eta,
\]
and where the middle tensor is the identity. Hence $E(e_n)=e_n$. In particular,
$E(b)=b$ for any $b\in\mathcal B$.
It remains to verify the bimodule property. By linearity, it is enough to consider
an elementary tensor
\[
x=ea'f,\qquad e,f\in\mathcal F,\quad a'\in S,
\]
written according to the above convention. Then
\[
E(e_nxe_m)
=
e_{j(e)\wedge n}\Phi(a')e_{j(f)\wedge m}
=
e_n\bigl(e_{j(e)}\Phi(a')e_{j(f)}\bigr)e_m
=
e_nE(x)e_m.
\]
Since $\mathcal B$ is generated by the projections $e_n$, this implies
\[
E(bxb')=bE(x)b'
\]
for any $b,b'\in\mathcal B$ and $x\in{\mathcal A}_{\rm rep}$. Thus $E$ is a
conditional expectation.
\hfill $\boxvoid$\\

\begin{Example}
{\rm 
We illustrate the action of $E$ on basic elements.
For the `sequential approximate identities' dedicated to labels, we have
\[
E(e_{1,n})=E(e_{2,n})=e_n
\]
for any $n$. Moreover,  
\[
E(a_1(j))=
\left\{\begin{array}{ll}
\varphi_1(a_1)p_j& {\rm if}\;\;a_1\in {\mathcal A}_{1}\\
\varphi_2(a_1)p_j& {\rm if}\;\;a_1\in {\mathcal A}_{2}
\end{array}
\right.
\]
for any $j$. Next, we obtain
\[
E(a_1(j)a_2(j))=
\left\{\begin{array}{ll}
\varphi_1(a_1a_2)p_j& {\rm if}\;\;a_1\in {\mathcal A}_{1}, a_2\in {\mathcal A}_{1}\\
\varphi_2(a_1a_2)p_j& {\rm if}\;\;a_1\in {\mathcal A}_{2}, a_2\in {\mathcal A}_{2}\\
\varphi_1(a_1)\varphi_2(a_2)p_j &{\rm if} \;\;a_1\in {\mathcal A}_{1}, a_2\in {\mathcal A}_{2}\\
\varphi_2(a_1)\varphi_1(a_2)p_j &{\rm if} \;\;a_1\in {\mathcal A}_{2}, a_2\in {\mathcal A}_{1}
\end{array}
\right.
\]
for any $j$. In turn, $E(a_1(j)a_2(k))=0$ when $j\neq k$. 
In order to obtain these formulas, it suffices to write \(a_1(j)\) and \(a_2(j)\)
as products of the canonical tensor product embeddings of \(a_1\) and \(a_2\)
and projections of the form \(e_{1,j}-e_{1,j-1}\) or \(e_{2,j}-e_{2,j-1}\),
where \(e_{i,0}=0\) for \(i=1,2\).
}
\end{Example}

\subheading{Reduction to scalar functionals \(\Psi_j\)}

There is a simple relation between mixed moments of orthogonal replicas 
associated with Motzkin words of height $j$, computed with respect to $E$, 
and those computed with respect to certain scalar-valued functionals $\Psi_j$. 
This relation allows one to reduce computations at a fixed height to a model 
resembling $\Phi$, while retaining the information encoded by $E$.

\begin{Definition}
{\rm 
For any $j\in \mathbb{N}$, define a normalized linear functional $\Psi_j$ on $\mathcal{A}_{{\rm rep}}$ by 
\[
\Psi_{j}:=\Psi_{1,j}\otimes \Psi_{2,j},
\]
where 
\[
\Psi_{i,j}=\hat{\varphi_{i}}^{\otimes j-1}\otimes \widetilde{\varphi_{i}}^{\otimes \infty}
\]
and $\hat{\varphi_{i}}$ is the trivial extension of $\varphi_{i}$ to $\widetilde{\mathcal A}_{i}$ 
vanishing on all products containing $p$ and $q$, respectively, where $i=1,2$. In particular, $\Psi_1=\Phi$. 
Thus $\Psi_j$ coincides with $\Phi$ on tensor components starting from level $j$, 
while the lower levels are collapsed via $\hat{\varphi_i}$.
}
\end{Definition}

\begin{Proposition}
Let
$a_1, \ldots , a_n\in \mathcal{A}_{1}\cup {\mathcal A}_{2}$ and let $j\in {\mathbb N}$.
\begin{enumerate}
\item
If $w=j_1 \cdots j_n\in \mathpzc{AM}$ and $h(w)=j$, then
\[
E(a_{1}(j_1)\cdots a_{n}(j_n))=\Psi_{j}(a_{1}(j_1)\cdots a_{n}(j_n))p_j.
\]
\item
Under the same assumptions, 
\[
\Psi_{j}(a_{1}(j_1)\cdots a_{n}(j_n))=\Phi(a_{1}(j_1')\cdots a_{n}(j_n')),
\]
where $j_k'=j_k-j+1$ for $k=1, \ldots, n$.
\item
If $(\hat{\beta}_{n})_{n\geq 1}$ are the mixed Boolean cumulants associated with $\Psi_j$, then
\[
\hat{\beta}_{n}(a_1(j), \ldots, a_n(j))=\beta_{n}(a_1, \ldots, a_n),
\]
for any $j$, where $({\beta}_n)_{n\geq 1}$ 
are the mixed Boolean cumulants associated with
the Boolean product of $\varphi_{1}$ and $\varphi_2$.
\end{enumerate}
\end{Proposition}
{\it Proof.}
If $w=j_1 \cdots j_n\in \mathpzc{AM}$ and $h(w)=j$, then 
\begin{eqnarray*}
E(a_{1}(j_1)\cdots a_{n}(j_n))&=&
E((e_{i,j}-e_{i,j-1})a(e_{i',j}-e_{i',j-1}))\\
&=&
(e_{j}-e_{j-1})\Phi(a)(e_j-e_{j-1})
\\
&=&
\Psi_{j}(a_{1}(j_1)\cdots a_{n}(j_n))p_j
\end{eqnarray*}
for some $i,i'\in \{1,2\}$ and $a\in S$, which proves (1). 
Property (2) follows easily from the definition of the sequence $(\Psi_j)_{j\geq 1}$. 
To justify (3), observe that we know from \cite{[L6]} that
the families 
\[
\{a_k(1):i_k=1\}\;\;\;{\rm and}\;\;\;\{a_k(1):i_k=2\}
\]
are Boolean independent with respect to $\Phi$. 
Similarly, the families
\[
\{a_k(j):i_k=1\}\;\;\;{\rm and}\;\;\;\{a_k(j):i_k=2\}
\]
are Boolean independent with respect to $\Psi_j$ for any $j$. 
Moreover, their distributions with respect to $\Psi_j$ 
agree with those of $\{a_k:i_k=1\}$ and $\{a_k:i_k=2\}$ with respect to $\varphi_1$ and
$\varphi_2$, respectively. Therefore, the mixed Boolean cumulants of $a_1(j), \ldots, a_n(j)$ with respect to
$\Psi_j$ agree with the mixed Boolean cumulants of $a_1, \ldots, a_n$ with respect to the Boolean product of $\varphi_1$ and
$\varphi_2$. The index $j$ is suppressed in $(\hat{\beta}_n)_{n\geq 1}$ since the values of these cumulants do not depend on $j$. 
This completes the proof.
\hfill $\boxvoid$\\

This proposition shows that mixed moments with respect to $E$ decompose according 
to the height of the associated Motzkin word, and at each fixed height $j$ they 
reduce to scalar-valued moments with respect to $\Psi_j$, for which replicas of 
color $j$ are Boolean independent.

\begin{Remark}
{\rm 
If we want to consider orthogonal replicas for conditionally free independent random variables, then
the distribution of an orthogonal replica of color $1$ may differ from that of a replica of color $j>1$. The same holds for
mixed Boolean cumulants of orthogonal replicas of color $1$ and of those of color $j>1$. All results of this paper can 
be modified accordingly. 
Let us also observe that we can assume that ${\mathcal A}_1, {\mathcal A}_{2}$ are unital *-algebras and that
normalized linear functionals $\varphi_1,\varphi_2$ are positive.
Then the functionals $\Phi_1,\Phi_2,\Phi,\Psi_j$ are also positive as tensor products of positive functionals if we set
$p^*=p$ and $q^*=q$.
}
\end{Remark}

\subheading{Reduction rules for \(E\)-moments}

There is no simple closed formula for all alternating mixed moments of orthogonal
replicas with respect to \(E\). We therefore record below a collection of reduction
rules sufficient for the moment computations used later. We shall usually omit the
colors of the replicas and write
\[
a:=a_1b_1a_2\cdots b_{n-1}a_n,
\]
where \(a_k=a_k(j_k)\), \(b_1,\ldots,b_{n-1}\in\mathcal B\), and
\(w=j_1\cdots j_n\) is a Motzkin word. The rules first eliminate non-admissible
color words, then remove local maxima by replacing the corresponding replicas by
projection factors, and finally simplify these projections. Iterating this procedure
reduces the computation to Boolean cumulants at constant levels.

\begin{Proposition}
Under the assumptions of Proposition 6.6,
\[
E(a)=0
\] 
whenever $w:=j_1\cdots j_n\notin \mathpzc{AM}_{n}(\ell)$ for any $b_1, \ldots, b_{n-1}
\in \mathcal{B}$, where 
\[
\mathpzc{AM}_{n}(\ell)=\{w\in \mathpzc{AM}_{n}: j_k=j_{k+1} \;\;if\;\;i_k=i_{k+1}, \; k\in [n-1]\}.
\]
In particular, if $i_1\neq \ldots \neq i_n$, then the above equation holds whenever $w\notin \mathpzc{AM}_{n}$.
\end{Proposition}
{\it Proof.}
If $j_1=j_n=1$, then the assertion follows immediately from Proposition 6.6. If $j_1\neq j_n$, then $w\notin \mathpzc{AM}$, but 
in this case the moment with respect to $E$ vanishes since $p_{j_1}\perp p_{j_n}$. If $j_1=j_n\neq 1$, but $w\notin \mathpzc{AM}$, then
either $j_k\geq j_1$ for all $k\in \{2,\ldots, n-1\}$ (then the moment essentially reduces to a moment corresponding to 
$w'=j_1'\cdots j_n'$ and thus it vanishes) or there exists $k\in \{2, \ldots, n-1\}$ such that $j_k<j_1$, but this produces $p^{\perp}$ or $q^{\perp}$ at site 
of color $j_k-1$ and thus the moment with respect to $E$ vanishes since $\Psi_j(a)=\Psi_j(p_jap_j)$ for any alternating product $a$ and $p^{\perp}$ will meet $p$ or $q^{\perp}$ will meet $q$ at this site. This completes the proof.
\hfill $\boxvoid$\\

As concerns `Motzkin moments', we begin with three simple types of mixed moments associated with constant Motzkin paths.

\begin{Lemma}
Let $w=j^{n}$ for some $j\in \mathbb{N}$ and let $(\hat{\beta}_{n})_{n\geq 1}$ be given by Proposition 7.2.
\begin{enumerate}
\item
If $i_1=\cdots =i_n$ and $b_1=\cdots =b_{n-1}=1$, then 
\[
E(a)=\hat{\beta}_{1}(a_1\cdots a_n)p_j.
\]
\item
If $i_1\neq \ldots \neq i_n$ and $b_1=\cdots =b_{n-1}=1$, then 
\[
E(a)= \hat{\beta}_{1}(a_1)\cdots \hat{\beta}_{1}(a_n)p_j.
\]
\item
If $i_{1}= \ldots = i_{n}$ and $b_1=\cdots =b_{n-1}=p_{j+1}$, then
\[
E(a)=\hat{\beta}_{n}(a_{1},\ldots, a_{n})p_j.
\]
\end{enumerate}
\end{Lemma}
{\it Proof.}
The property (1) is obvious.
The factorization in (2) follows from the definition of orthogonal replicas since 
all variables with label $1$ are separated by the projection $p$ at site $(1,j)$
of the infinite tensor product and all variables with label $2$ are separated by the projection $q$
at site $(2,j)$. In effect, this leads to the product of moments of order one which 
are equal to Boolean cumulants when applying $\Psi_j$, which proves (2).
In order to show (3), let us suppose that $a_1, \ldots ,a_n$ are orthogonal replicas of $A_1, \ldots, A_n\in \mathcal{A}_{1}$ respectively (for replicas of label 2 the proof is analogous). 
We will use the decomposition
\[
p_{j+1}=\sum_{(\alpha, \gamma)\in \Theta}(\ldots 1_1 \otimes p^{\alpha} \otimes p \ldots) \otimes 
(\ldots 1_2 \otimes q^{\gamma}\otimes q \ldots),
\]
where $\Theta=\{(0, \perp), (\perp,1)\}$, so that 
$(p^{\alpha}, q^{\gamma})\in \{(1,q^{\perp}), (p^{\perp},q)\}$ and 
only sites of colors $j-1, j, j+1$ of both infinite tensor products are explicitly shown. 
Since $a_1, \ldots , a_n$ contain $q$ at site $(2,j)$,
all terms vanish except the one in which 
$(p^{\alpha}, q^{\gamma})=(p^{\perp}, q)$ in each $p_{j+1}$. Therefore,
\[
a=(\ldots 1_1\otimes (A_1p^{\perp}\cdots p^{\perp}A_n)\otimes p \ldots) 
\otimes 
(\ldots q^{\perp}\otimes q \otimes q \ldots),
\] 
where again only sites of colors $j-1,j,j+1$ are shown.
If we apply $\Psi_j$ to this product, we get
\[
\Psi_j(a)
=
\tilde{\varphi}_{1}(A_1p^{\perp}\cdots p^{\perp}A_n)={\beta}_n(A_1,\ldots , A_n)=\hat{\beta}_{n}(a_1, \ldots, a_n)
\]
by the definition of ${\Psi}_{j}$ and an elementary computation giving the second equality 
(see \cite{[L4]} for the proof). This gives the desired formula for $E(a)$ since $h(w)=j$.
\hfill $\boxvoid$\\

Even if $w$ is a non-constant word of height $j$ and there is a projection $p_j$ between
two replicas in the alternating product $a$, the corresponding moment factorizes.

\begin{Lemma} {\rm (Factorization)}
Let $w=j_1\cdots j_n\in \mathpzc{AM}$ and let $1 \leq k < n$. Then
\[
E(\cdots a_{k}b_kp_ja_{k+1}\cdots )=E(\cdots a_{k}b_k)E(a_{k+1}\cdots )
\]
whenever $j=h(w)$ and the replicas have arbitrary labels.
\end{Lemma}
{\it Proof.}
The projection $p_j$ puts projections $p$ and $q$ at all tensor sites of colors $r\geq j$ (for both labels). 
This leads to the factorization of the mixed moment with respect to $\Psi_j$ for any $j$
and thus it gives the factorization of the moment with respect to $E$ since $p_j^2=p_j$. 
This completes the proof.
\hfill $\boxvoid$\\

\begin{Lemma} {\rm (Local maximum)}
Let \(w=j_1\ldots j_n\in \mathpzc{AM}\). Suppose that, for some
\(1\leq k\leq n\), \(j_{k-1}\leq j_k\geq j_{k+1}\), and
\(i_{k-1}\neq i_k\neq i_{k+1}\). Then
\[
E(\cdots a_{k-1}b_{k-1}a_kb_ka_{k+1}\cdots )
=
\hat{\beta}_{1}(a_k)
E( \cdots a_{k-1}b_{k-1} p_{j_k}b_ka_{k+1}\cdots ).
\]
The cases \(k=1\) and \(k=n\) are included with the natural one-sided
interpretation of the local maximum and neighboring-label conditions, with the missing boundary factors omitted.
\end{Lemma}
{\it Proof.}
To fix attention, let $i_{k-1}=i_{k+1}=1$, $i_k=2$ and take $A_{k-1}, A_{k+1}\in {\mathcal A}_{1}$, $A_k\in {\mathcal A}_{2}$, where
$1<k<n$. First, suppose that $b_{k-1}=b_{k}=1$. Let $(j_{k-1},j_k,j_{k+1})=(j,j+1,j)$. 
We have
\begin{eqnarray*}
a_{k-1}&=&(\ldots 1_1\otimes A_{k-1}\otimes 1_1\ldots )\otimes 
(\ldots q^{\perp}\otimes q \otimes q \ldots),\\
a_{k}&=&(\ldots 1_1\otimes p^{\perp}\otimes p\ldots )\otimes 
(\ldots 1_{2}\otimes 1_{2} \otimes A_{k}\ldots),\\
a_{k+1}&=&(\ldots 1_1\otimes A_{k+1}\otimes 1_1\ldots )\otimes 
(\ldots q^{\perp}\otimes q \otimes q \ldots),
\end{eqnarray*}
where only sites of colors $j-1,j, j+1$ are shown.
Then
\begin{eqnarray*}
a_{k-1}a_{k}a_{k+1}&=&(\ldots 1_1\otimes A_{k-1}p^{\perp}A_{k+1}\otimes p\ldots )\otimes 
(\ldots q^{\perp}\otimes q \otimes qA_{k}q \ldots).
\end{eqnarray*}
Under $E$ this expression can be replaced by 
$\hat{\beta}_{1}(a_{k})a_{k-1}Pa_{k+1}$, where 
\[
P=(\ldots 1_1\otimes p^{\perp}\otimes p\ldots )\otimes 
(\ldots 1_2\otimes 1_2 \otimes q \ldots)
\]
since $qA_{k}q-\varphi_{2}(A_{k})q\in {\rm Ker}\widetilde{\varphi}_{2}$. Next, 
$P$ can be replaced by $p_{j+1}$ since 
\[
p_{j+1}-P=(\ldots 1_1\otimes p\otimes p\ldots )\otimes 
(\ldots 1_2\otimes q^{\perp} \otimes q \ldots)
\]
and $a_{k-1}(p_{j+1}-P)a_{k+1}=0$ since $q$ meets $q^{\perp}$. 
In turn, if 
$(j_{k-1},j_k,j_{k+1})=(j,j+1,j+1)$ or $(j_{k-1},j_k,j_{k+1})=(j+1,j+1,j)$, 
the computation is very similar, the only difference is that $q^{\perp}$ is multiplied by $q$ only on one side.
Finally, if $(j_{k-1},j_k,j_{k+1})=(j,j,j)$, then 
\begin{eqnarray*}
a_{k-1}a_{k}a_{k+1}&=&(\ldots p^{\perp}\otimes A_{k-1}pA_{k+1}\otimes p\ldots )\otimes 
(\ldots q^{\perp}\otimes qA_kq \otimes q \ldots).
\end{eqnarray*}
Under $E$ this expression can be replaced by 
$\hat{\beta}_{1}(a_{k})a_{k-1}e_ja_{k+1}$, where 
\[
e_j=(\ldots 1_1\otimes p\otimes p\ldots )\otimes 
(\ldots 1_2\otimes q \otimes q \ldots)
\]
which can be replaced by $p_{j}$ since 
\[
e_j-p_{j}=(\ldots p\otimes p\otimes p\ldots )\otimes 
(\ldots q\otimes q \otimes q \ldots)
\]
and $a_{k-1}(e_j-p_{j})a_{k+1}=0$ since again $q$ meets $q^{\perp}$. If $k=1$ or $k=n$, the computations are 
similar. If $b_{k-1}\neq 1$ or $b_{k}\neq 1$, the proof is identical 
since \({\mathcal B}\) is spanned by the projections \(e_n\), \(n\in\mathbb N\cup\{\infty\}\)
and each $e_n$ 
has $p$ or $q$ or $1$ at each tensor site and it can be easily seen that 
this does not affect our computations.
This completes the proof.
\hfill $\boxvoid$\\

Now, we would like to treat the case when we have a projection lying 
in between two replicas. The trivial cases were covered by Proposition 6.4.
The non-trivial ones which will be of interest 
are treated in the lemmas given below.

\begin{Lemma} {\rm (Projection reduction)}
If
$i_{k}=i_{k+1}$ and 
$j_{k}=j_{k+1}=j$ for some $1\leq k <  n$, then
\[
E(\cdots a_{k}b_{k}p_{j+1}a_{k+1}\cdots )= E( \cdots a_{k}b_{k}p_{j}^{\perp}a_{k+1}\cdots )
\]
whenever $h(w)\leq j$ and $w=j_1\ldots j_n\in \mathpzc{AM}$, where $p_j^{\perp}=1-p_j$.
\end{Lemma}
{\it Proof.}
To fix attention, 
suppose that $A_{k},A_{k+1} \in \mathcal{A}_{1}$.
As in the proof of Lemma 7.10, it suffices to treat the case when $b_k=1$.
We have
\begin{eqnarray*}
a_{i}&=&(\ldots 1_1\otimes A_i\otimes 1_1\ldots )\otimes 
(\ldots q^{\perp}\otimes q \otimes q \ldots),\\
p_{j+1}&=&(\ldots 1_1\otimes p^{\perp}\otimes p\ldots )\otimes 
(\ldots 1_{2}\otimes 1_{2} \otimes q\ldots)\\
&+&(\ldots 1_1\otimes p\otimes p\ldots )\otimes 
(\ldots 1_{2}\otimes q^{\perp} \otimes q\ldots),\\
a_{k}p_{j+1}a_{k+1}
&=&
(\ldots 1_{1}\otimes A_kp^{\perp}A_{k+1} \otimes p \ldots )\otimes (\ldots q^{\perp} \otimes q \otimes q \ldots),
\end{eqnarray*}
where $i=k,k+1$, and where only sites of colors $j-1,j,j+1$ are shown. Now, if we replace here $p_{j+1}$ by 
\begin{eqnarray*}
p_{j}&=&(\ldots p^{\perp}\otimes p\otimes p\ldots )\otimes 
(\ldots 1_{2}\otimes q \otimes q\ldots)\\
&+&(\ldots p\otimes p\otimes p\ldots )\otimes 
(\ldots q^{\perp}\otimes q \otimes q\ldots),
\end{eqnarray*}
we obtain 
\begin{eqnarray*}
a_{k}p_{j}a_{k+1}
&=&
(\ldots 1_{1}\otimes A_kpA_{k+1} \otimes p \ldots )\otimes (\ldots q^{\perp} \otimes q \otimes q \ldots)
\end{eqnarray*}
since we can replace $1_2$ by $q^{\perp}$ and then use $p+p^{\perp}=1_1$. 
In effect, we obtain
\[
a_{k}p_{j+1}a_{k+1}-a_{k}p_j^{\perp}a_{k+1}
\]
\[
= \left((\ldots 1_{1}\otimes A_kA_{k+1} \otimes p \ldots )
-
(\ldots 1_{1}\otimes A_kA_{k+1} \otimes 1_1 \ldots )\right)
\otimes (\ldots q^{\perp} \otimes q \otimes q \ldots).
\]
However, if we multiply this expression 
from the left by $a_{k-1}$ and from the right by $a_{k+2}$, both
of label $2$ and color from the set $\{j-1,j,j+1\}$, we get zero.
The same holds if $a_{k-1}$ or $a_{k+2}$ is of label $1$ and color $j$ (colors $j-1$ and $j+1$ 
make the whole product vanish by Proposition 6.4, so this case is trivial). Of course, 
we also get zero if $k=1$ or $k=n-1$ since we have $q^{\perp}$ multiplied from the left or from the right 
by $q$ produced by $E$.
Therefore, we can replace $a_{k}p_{j+1}a_{k+1}$ under $E$ by 
$
a_{k}p_{j}^{\perp} a_{k+1}
$
under the assumptions of the Lemma, which completes the proof.
\hfill $\boxvoid$\\

\begin{Lemma} {\rm (Projection absorption)}
If
$i_{k}\neq i_{k+1}$ and 
$(j_{k},j_{k+1})\in \{(j,j+1),(j+1,j)\}$ for some $1\leq k <  n$, then
\[
E(\cdots a_{k}b_{k}p_{j+1}a_{k+1}\cdots )= E( \cdots a_{k}b_{k}a_{k+1}\cdots )
\]
whenever $h(w)\leq j$ and $w=j_1\ldots j_n\in \mathpzc{AM}$.
\end{Lemma}
{\it Proof.}
Without loss of generality, let $i_k=1$, $i_{k+1}=2$ and $(j_k,j_{k+1})=(j,j+1)$. 
As in the preceding lemmas, it suffices to treat the case when $b_k=1$.
We have
\begin{eqnarray*}
a_{k}&=&(\ldots 1_1\otimes A_k\otimes 1_1\ldots )\otimes 
(\ldots q^{\perp}\otimes q \otimes q \ldots),\\
a_{k+1}&=&(\ldots 1_1\otimes p^{\perp}\otimes p\ldots )\otimes 
(\ldots 1_2\otimes 1_2 \otimes A_{k+1} \ldots),\\
p_{j+1}&=&(\ldots 1_1\otimes p^{\perp}\otimes p\ldots )\otimes 
(\ldots 1_{2}\otimes 1_{2} \otimes q\ldots)\\
&+&(\ldots 1_1\otimes p\otimes p\ldots )\otimes 
(\ldots 1_{2}\otimes q^{\perp} \otimes q\ldots),\\
a_{k}p_{j+1}a_{k+1}
&=&
(\ldots 1_{1}\otimes A_kp^{\perp} \otimes p \ldots )\otimes (\ldots q^{\perp} \otimes q \otimes qA_{k+1}\ldots),
\end{eqnarray*}
where only sites of colors $j-1,j,j+1$ are shown. 
Note that the last expression is equal to $a_ka_{k+1}$.
Let us remark that in this case $a_kp_ja_{k+1}=0$ and therefore we can also replace $a_kp_{j+1}a_{k+1}$ by 
$a_kp_{j}^{\perp}a_{k+1}$, which completes the proof.
\hfill $\boxvoid$\\

\begin{Lemma} {\rm (Nested plateau)}
Assume that \(j_k=\cdots=j_m=j\), \(j_{k-1}<j>j_{m+1}\), and
\(i_{k-1}\neq i_k\), \(i_m\neq i_{m+1}\), with the endpoint cases understood in
the natural one-sided sense. Then
\[
E(\cdots a_kb_kp_j\cdots p_ja_m\cdots )
=
E(\cdots E(a_kb_kp_j\cdots p_j a_m)\cdots ).
\]
\end{Lemma}

{\it Proof.}
We prove the assertion in the case \(k>1\) and \(m<n\); the endpoint cases are
obtained by omitting the missing boundary factor. To fix attention, suppose first
that \(A_k,\ldots,A_m\in\mathcal A_1\) and that 
$A_{k-1},A_{m+1}\in \mathcal{A}_{2}$. The other cases are analogous, 
with \(p\) and \(q\) interchanged when
necessary.

As in the preceding lemmas, by multilinearity it is enough to consider the case
in which the elements of \(\mathcal B\) appearing inside the plateau are
projections \(e_r\). Such projections only insert tensor factors \(1,p,q\) into
the internal product over the plateau and do not change the boundary projection
pattern at the levels \(j-1,j,j+1\) used below. Thus we suppress them in the
displayed computation; in the general case the factor
\(A_kp\cdots pA_m\) is replaced by the corresponding product with these inserted
projections.

Using the tensor-product realization and displaying only the sites of colors
\(j-1,j,j+1\), we obtain
\[
a_kp_j a_{k+1}\cdots p_j a_m
=
(\ldots 1_1\otimes (A_kp\cdots pA_m)\otimes p \ldots)
\otimes
(\ldots q^{\perp}\otimes q\otimes q \ldots).
\]
Thus, under \(E\), the plateau contributes its own \(\mathcal B\)-valued moment
\[
E(a_kp_j a_{k+1}\cdots p_j a_m).
\]
It remains to check that the projection left at level \(j\) is the same as the
projection occurring in this \(E\)-moment when the plateau is viewed as a nested
subproduct.

Since \(w\) is a Motzkin word, the inequalities \(j_{k-1}<j\) and \(j_{m+1}<j\)
imply \(j_{k-1}=j_{m+1}=j-1\). The projection
\[
(\ldots 1_1\otimes p\otimes p\ldots )
\otimes
(\ldots q^{\perp}\otimes q\otimes q \ldots)
\]
can therefore be replaced, inside the larger moment, by \(p_j\). Indeed, in the
decomposition
\[
\begin{aligned}
p_j
&=
(\ldots 1_1 \otimes p\otimes p\ldots)
\otimes
(\ldots q^{\perp}\otimes q\otimes q \ldots)  \\
&\quad +
(\ldots p^{\perp}\otimes p\otimes p\ldots)
\otimes
(\ldots q\otimes q\otimes q \ldots),
\end{aligned}
\]
the second term vanishes after multiplication by the neighbouring replicas
\(a_{k-1}\) and \(a_{m+1}\), because their tensor components contain \(q^\perp\)
at the relevant site.
Hence the plateau may be replaced inside the larger \(E\)-moment by its own
\(E\)-moment, which gives the required nested formula.
\hfill \(\boxvoid\)\\

\begin{Example}
{\rm
Let
\[
a=a_1(1)a_2(2)a_3(3)a_4(2)a_5(1),
\]
where \(i_1=i_3=i_5=1\) and \(i_2=i_4=2\). The associated Motzkin word is
\(w=12321\). It has a local maximum at position \(3\), and therefore Lemma 7.10 gives
\[
E(a)
=
\widehat{\beta}_1(a_3)
E(a_1(1)a_2(2)p_3a_4(2)a_5(1)).
\]
Since \(a_2(2)\) and \(a_4(2)\) have the same label and the same color,  Lemma 7.11
gives
\[
E(a_1(1)a_2(2)p_3a_4(2)a_5(1))
=
E(a_1(1)a_2(2)p_2^{\perp}a_4(2)a_5(1)).
\]
Using \(p_2^{\perp}=1-p_2\), we obtain
\[
\begin{aligned}
&E(a_1(1)a_2(2)p_2^{\perp}a_4(2)a_5(1))  \\
&\qquad =
E(a_1(1)a_2(2)a_4(2)a_5(1))
-
E(a_1(1)a_2(2)p_2a_4(2)a_5(1)).
\end{aligned}
\]
We compute the two terms separately. Another application of the local maximum rule,
applied to the product \(a_2(2)a_4(2)\), followed by a projection reduction gives
\[
\begin{aligned}
E(a_1(1)a_2(2)a_4(2)a_5(1))
&=
\widehat{\beta}_1(a_2a_4)\,
E(a_1(1)p_2a_5(1))                                      \\
&=
\widehat{\beta}_1(a_2a_4)\,
E(a_1(1)p_1^{\perp}a_5(1))                                \\
&=
\widehat{\beta}_1(a_2a_4)\,
\widehat{\beta}_2(a_1,a_5)p_1.
\end{aligned}
\]
For the second term, the plateau at level \(2\) can be evaluated inside the larger
moment:
\[
\begin{aligned}
E(a_1(1)a_2(2)p_2a_4(2)a_5(1))
&=
E\bigl(a_1(1)E(a_2(2)p_2a_4(2))a_5(1)\bigr)       \\
&=
\widehat{\beta}_1(a_2)\widehat{\beta}_1(a_4)
E(a_1(1)p_2a_5(1))                                \\
&=
\widehat{\beta}_1(a_2)\widehat{\beta}_1(a_4)
\widehat{\beta}_2(a_1,a_5)p_1.
\end{aligned}
\]
Combining the above identities, and using the definition of the second Boolean
cumulant, we get
\[
\begin{aligned}
E(a)
&=
\widehat{\beta}_1(a_3)
\widehat{\beta}_2(a_2,a_4)
\widehat{\beta}_2(a_1,a_5)p_1.
\end{aligned}
\]
}
\end{Example}

\medskip
\noindent{\it Summary of the reduction procedure.}
The lemmas above, together with the preceding example, illustrate an effective
reduction procedure for the \(E\)-moments associated with Motzkin words: one
successively removes local maxima, simplifies the projections left behind, and
eventually obtains expressions in terms of Boolean cumulants at constant levels.
\medskip

\subheading{Monotone independence from two consecutive colors}

Lemma 7.10 resembles the definition of monotone independence, except that we have a non-strict 
local maximum and the formulation refers to colors rather than labels. Nevertheless, 
the definition of monotone independence can be recovered by combining Lemmas 7.10 and 7.11, as we show below. 

The key mechanism is an asymmetric use of colors: one family is represented by replicas 
of a single color, whereas the other uses replicas of two consecutive colors. This 
asymmetry is responsible for the non-symmetric nature of monotone convolution and 
provides a new perspective on this notion (cf.~\cite{[L6]}).

\begin{Lemma}
Let $(a_k,a_k')$ be a pair of orthogonal replicas of the same variable with label $i_k$ and colors $(j,j+1)$, respectively,
for any $k=1, \ldots, n$, where $j\in \mathbb{N}$ is fixed.
\begin{enumerate}
\item 
If $i_1\neq \ldots \neq i_n$ and $1\leq k \leq n$, then 
\[
E(a_1\cdots a_{k-1}(a_{k}+a_{k}')a_{k+1}\cdots a_n)=E(a_{k})E(a_1\cdots a_{k-1}a_{k+1}\cdots a_n).
\]
\item
If, for $i\in \{1,2\}$, we are given two families of orthogonal replicas
\[
{\mathcal F}_{i}=\{a(j): a\in {\mathcal A}_{i}\} \;\;\;and\;\;\;
{\mathcal G}_{i}:=\{a(j)+a(j+1): a\in \mathcal{A}_{i}\},
\] 
then $({\mathcal F}_{1}, {\mathcal G}_{2})$ and $({\mathcal F}_{2}, {\mathcal G}_{1})$ are monotone 
independent with respect to $E$.
\end{enumerate}
\end{Lemma}
{\it Proof.}
In the proof of Lemma 7.10, let $a_{k-1}=a_{k-1}(j)$, $a_{k}=a_{k}(j)+a_k(j+1)$, $a_{k+1}=a_{k+1}(j)$, where
$a_{k-1},a_{k+1}\in {\mathcal A}_{1}$ and $a_k\in \mathcal{A}_{2}$. Then, using Lemma 7.10 and Lemma 7.11, we obtain
\begin{eqnarray*}
E(a_1\cdots a_{k-1}(a_{k}+a_k')a_{k+1}\cdots a_n)&=&E(a_{k})E(a_1\cdots a_{k-1}(p_j+p_{j}^{\perp})a_{k+1}\cdots a_n)\\
&=&E(a_{k})E(a_1\cdots a_{k-1}a_{k+1}\cdots a_n)
\end{eqnarray*}
for any $1\leq k \leq n$, which proves (1). 
Since $E(a_k')$ is supported on $p_{j+1}$, whereas 
$E(a_1\cdots a_{k-1}a_{k+1}\cdots a_n)$ is supported on $p_j$,
their product vanishes. Thus we obtain 
\[
E(a_1\cdots a_{k-1}(a_{k}+a_k')a_{k+1}\cdots a_n)=E(a_{k}+a_k')E(a_1\cdots a_{k-1}a_{k+1}\cdots a_n),
\]
which gives monotone independence of $({\mathcal F}_1, {\mathcal G}_{2})$ with respect to $E$ 
since we allow $k=1$ and $k=n$. The proof for $({\mathcal F}_{2}, {\mathcal G}_{1})$ is similar.
This completes the proof of (2).
\hfill $\boxvoid$\\

In consequence, in the case of two algebras, it is possible to describe mixed 
moments of monotone independent random variables using Motzkin cumulants rather than 
monotone cumulants of Hasebe and Saigo \cite{[HS]}.

\section{$w$-Boolean cumulants of orthogonal replicas}

We now study the $w$-Boolean cumulants associated with the
$\mathcal B$-valued conditional expectation $E$ introduced in Section~7.
These results will later be used to analyze Motzkin cumulants of orthogonal
replicas.

Although $w$-Boolean cumulants were defined for more general arguments, we
now specialize to the case when these arguments are orthogonal replicas.
Moreover, we assume that the colors of these replicas are encoded by the word
$w$. In this section, we study $B_w(a_1,\ldots,a_n)$ when
$a_1,\ldots,a_n$ are orthogonal replicas with labels and colors encoded by
$\ell=i_1\cdots i_n$ and $w=j_1\cdots j_n\in\mathpzc{AM}$, respectively.
We would like to study $w$-Boolean cumulants and Motzkin cumulants 
associated with the sequence of conditional expectations of Section 7. 

Although they were defined for more general arguments, we 
would like to study them now in the case when 
these arguments are orthogonal replicas. 
Moreover, we will assume that the colors of these replicas are encoded by the
word $w$.
In this Section, we study $B_{w}(a_1, \ldots, a_n)$
when $a_1, \ldots, a_n$ are orthogonal replicas 
with labels and colors encoded by $\ell=i_1 \cdots i_n$ and
$w=j_1\cdots j_n\in \mathpzc{AM}$, respectively.

As we demonstrated in \cite{[L6]}, in order to compute such $B_w(a_1, \ldots, a_n)$ in terms of mixed Boolean cumulants, 
one can use the set of noncrossing partitions which are monotonically adapted to colors and labels.

\begin{Definition}
{\rm 
We will say that the partition $\pi=(\pi_0,w)\in \mathcal{NC}(w)$ 
is {\it adapted} to the labeling $\ell=i_1\cdots i_n$ if $i_k=i_m$ whenever $k,m\in V$ 
for some $V\in \pi$. The set of such partitions will be denoted by $\mathcal{NC}(w,\ell)$. 
The partition $\pi=(\pi_0,w)\in \mathcal{M}(w)\cap \mathcal{NC}(w, \ell)$  
will be called {\it monotonically adapted to $\ell$} if 
the labels assigned to blocks of $\pi$ alternate in saturated chains of $\pi$. 
By a {\it saturated chain} in $\pi\in \mathcal{M}(w)$
we understand a sequence of its blocks $(v_{j_1}, \ldots , v_{j_p})$ such that 
\[
v_{j_1}=o(v_{j_2}), \ldots , v_{j_{p-1}}=o(v_{j_{p}}),
\]
where $v_{j_1}, \ldots , v_{j_p}$ are constant words. 
The set of such partitions will be denoted by $\mathcal{M}(w, \ell)$.
We will also use the self-explanatory notations $\mathcal{NC}_{{\rm irr}}(w, \ell)$ and $\mathcal{M}_{{\rm irr}}(w,\ell)$ 
for the subfamilies of those defined above consisting of irreducible partitions.
}
\end{Definition}

\begin{Example}
{\rm
Let
\[
w=1232321,\qquad \ell=1212121.
\]
Consider the irreducible noncrossing partition
\[
\pi_0=\big\{\{1,7\},\{2,4,6\},\{3\},\{5\}\big\}.
\]
Following the convention introduced earlier, we identify the corresponding
$\pi=(\pi_0,w)$ with the set of pairs
\[
\pi=
\{(\{1,7\},11),
(\{2,4,6\},222),
(\{3\},3),
(\{5\},3)\}.
\]
It can be represented by the diagram

\unitlength=1mm
\special{em:linewidth 1pt}
\linethickness{0.5pt}
\begin{picture}(180.00,20.00)(5.00,15.00)

\put(62.00,20.50){$\pi=$}

\put(72.00,20.00){\circle*{1.00}}
\put(75.00,20.00){\circle*{1.25}}
\put(78.00,20.00){\circle*{1.60}}
\put(81.00,20.00){\circle*{1.25}}
\put(84.00,20.00){\circle*{1.60}}
\put(87.00,20.00){\circle*{1.25}}
\put(90.00,20.00){\circle*{1.00}}

% outer block {1,7}
\put(72.00,20.00){\line(0,1){7.0}}
\put(90.00,20.00){\line(0,1){7.0}}
\put(72.00,27.00){\line(1,0){18.0}}

% middle block {2,4,6}
\put(75.00,20.00){\line(0,1){4.5}}
\put(81.00,20.00){\line(0,1){4.5}}
\put(87.00,20.00){\line(0,1){4.5}}
\put(75.00,24.50){\line(1,0){12.0}}

% singletons {3} and {5}
%\put(78.00,20.00){\line(0,1){2.0}}
%\put(84.00,20.00){\line(0,1){2.0}}

% word labels
\put(71.20,17.00){$\scriptscriptstyle 1$}
\put(74.20,17.00){$\scriptscriptstyle 2$}
\put(77.20,17.00){$\scriptscriptstyle 3$}
\put(80.20,17.00){$\scriptscriptstyle 2$}
\put(83.20,17.00){$\scriptscriptstyle 3$}
\put(86.20,17.00){$\scriptscriptstyle 2$}
\put(89.20,17.00){$\scriptscriptstyle 1$}

%\put(92.00,20.00){,}

\end{picture}

\noindent
and it belongs to $\mathcal{M}_{{\rm irr}}(w,\ell)$. Indeed, it is
irreducible and adapted both to the coloring $w$ and to the labeling $\ell$:
the block $\{1,7\}$ has label $1$, the block $\{2,4,6\}$ has label $2$,
and the singleton blocks $\{3\}$ and $\{5\}$ have label $1$. Moreover,
the nearest outer block of the block with color-word $222$ has color-word
$11$, while each of the two singleton blocks with color-word $3$ has nearest
outer block with color-word $222$. Thus the two saturated chains have
color-words
\[
11,\ 222,\ 3,
\]
and the corresponding labels are $1,2,1$. Hence the labels alternate along
these chains.
On the other hand, if we keep the same word $w$ and the same partition $\pi$,
but take
\[
\ell'=1222121,
\]
then $\pi$ is no longer in $\mathcal{M}_{{\rm irr}}(w,\ell')$. It is still
adapted to the labeling $\ell'$, since each block carries one label, but it is
not monotonically adapted to $\ell'$. Namely, along the saturated chain ending
at the singleton block $\{3\}$, the corresponding labels are $1,2,2$, which do
not alternate. Hence $\pi\notin\mathcal{M}_{{\rm irr}}(w,\ell')$.
}
\end{Example}

\begin{Definition}
{\rm 
For any $w=j_1\cdots j_n\in \mathpzc{AM}$ and any $\ell=i_1\cdots i_n$,
define a family
$
\{\hat{\beta}_{\pi}:\pi \in \mathcal{M}(w,\ell)\}
$
of multilinear functionals on ${\mathcal A}_{{\rm rep}}^{n}$ by
\[
\hat{\beta}_{\pi}[a_1, \ldots, a_n]
=
\prod_{v\in \pi}\hat{\beta}_v[a_1, \ldots, a_n],
\]
where the product is taken over all blocks $v$ of $\pi$. More precisely, if
$V=\{k_1<\cdots<k_m\}$ is the corresponding block of $\pi_0$ and
$v=j_{k_1}\cdots j_{k_m}$ is the associated subword of $w$, then we set
\[
\hat{\beta}_{v}[a_1,\ldots,a_n]
=
\hat{\beta}_{m}(a_{k_1},\ldots,a_{k_m}),
\]
where $\hat{\beta}_{m}$ is the Boolean cumulant associated with the functional
$\Psi_{h(v)}$, where the height $h(v)$ is
suppressed in the notation.
}
\end{Definition}

\begin{Theorem}
If $a_1, \ldots, a_n$ are orthogonal replicas with labels and colors encoded by $\ell=i_1 \cdots i_n$ and
$w=j_1\cdots j_n\in \mathpzc{AM}$, respectively, where $h(w)=j$, then 
\[
B_w(a_1, \ldots , a_n)=\left(\sum_{\pi\in \mathcal{M}_{{\rm irr}}(w, \ell)}\hat{\beta}_{\pi}[a_1,\ldots, a_n] \right)p_j.
\]
In particular, if $w=j^n$, then $B_w(a_1, \ldots, a_n)=\hat{\beta}_{n}(a_1, \ldots, a_n)p_j$.
If  the labels are identical and the colors are not all equal, then $B_w(a_1, \ldots, a_n)=0$.
\end{Theorem}
{\it Proof.}
The general formula is an adaptation of Theorem 8.2 in \cite{[L6]}, where an
analogous formula was proved for scalar-valued Boolean cumulants of orthogonal
replicas associated with the functional $\Phi$. In the present setting we use
the $\mathcal B$-valued expectation $E$. By Proposition 7.5, when the word
$w=j_1\cdots j_n$ has height $h(w)=j$, the corresponding $E$-moment is obtained
from the scalar moment with respect to $\Psi_j$ and is then multiplied by the
projection $p_j$. Thus the proof of Theorem 8.2 in \cite{[L6]}, with $\Phi$
replaced by $\Psi_j$, gives the same sum over partitions monotonically adapted
to the colors and labels, with the additional final factor $p_j$.
Indeed, the Boolean cumulants associated with $\Psi_j$, when multiplied by
$p_j$, agree with $B_w(a_1,\ldots,a_n)$, since $\Psi_j(a_k\cdots a_m)=0$
whenever $j_k\neq j$ or $j_m\neq j$. Therefore ${\rm Int}(n)$ reduces to
${\mathcal Int}(w)$ for $w\in\mathpzc{AM}$. This proves the general formula
for $B_w(a_1,\ldots,a_n)$.

If $w=j^n$ and the labels are identical, then
$\mathcal{M}_{{\rm irr}}(w,\ell)=\{\hat{1}_w\}$, and hence
\[
B_w(a_1,\ldots,a_n)=\hat{\beta}_{n}(a_1,\ldots,a_n)p_j.
\]
If $w=j^n$ and the labels are not identical, then
$\mathcal{M}_{{\rm irr}}(w,\ell)=\emptyset$, hence
$B_w(a_1,\ldots,a_n)=0$, while also
$\hat{\beta}_{n}(a_1,\ldots,a_n)=0$ by Proposition 7.5(2).
Finally, if the labels are identical but the colors are not all equal,
then there exist neighboring replicas with different colors. In the Boolean
moment-cumulant expansion, these replicas either occur in the same moment or in
two neighboring moments. In the first case the moment vanishes by orthogonality;
in the second case one of the neighboring moments starts or ends at a color
different from $j$, and therefore vanishes under $\Psi_j$. Thus, using
Proposition 4.3, we get
\[
B_w(a_1,\ldots,a_n)=0.
\]
This completes the proof. \hfill $\boxvoid$\\

\begin{Example}
{\rm Let us illustrate Theorem 8.4 for orthogonal replicas with non-identical labels and colors such that $\mathcal{M}_{{\rm irr}}(w, \ell)\neq \emptyset$. We have
\begin{eqnarray*}
B_{121}(x_{k_1}, y_{k_2}, x_{k_3})&=&\hat{\beta}_{2}(x_{k_1}, x_{k_3})\hat{\beta}_{1}(y_{k_2})p_1,\\
B_{1221}(x_{k_1},y_{k_2},y_{k_3},x_{k_4})&=&\hat{\beta}_{2}(x_{k_1},x_{k_4})\hat{\beta}_{2}(y_{k_2},y_{k_3})p_1\\
&&\!\!\!\!\!+\;\hat{\beta}_{2}(x_{k_1},x_{k_4})\hat{\beta}_{1}(y_{k_2})\hat{\beta}_{1}(y_{k_3})p_1,\\
B_{1121}(x_{k_1},x_{k_2},y_{k_3},x_{k_4})&=&\hat{\beta}_{3}(x_{k_1},x_{k_2},x_{k_4})\hat{\beta}_{1}(y_{k_3})p_1,\\
B_{12321}(x_{k_1},y_{k_2},x_{k_3},y_{k_4},x_{k_5})&=&\hat{\beta}_2(x_{k_1},x_{k_5})\hat{\beta}_2(y_{k_2},y_{k_4})\hat{\beta}_1(x_{k_3})p_1
\end{eqnarray*}
where $\{x_k:k\in I_1\}$ and $\{y_k:k\in I_2\}$ are families of orthogonal replicas with labels $1$ and $2$, respectively, and colors are encoded by $w$.
}
\end{Example}

For comparison, recall that the formula for mixed Boolean cumulants of free
random variables derived in \cite{[FMNS]} can be recovered from the formula for
orthogonal replicas obtained in \cite{[L6]} by summing over the corresponding
Motzkin refinements. Thus the replica approach reflects the combinatorial
structure of freeness at a finer level. Keeping the notation
$(\hat{\beta}_n)_{n\geq 1}$ for the Boolean cumulants associated with
$\Psi_j$, let $(\beta_n)_{n\geq 1}$ denote the Boolean cumulants associated
with $\varphi$.

\begin{Corollary} \cite{[L6]}
Let ${\mathcal A}_{1}$ and ${\mathcal A}_{2}$ be unital free subalgebras of $({\mathcal A}, \varphi)$.
Then
\[
\beta_n(a_{1}, \ldots , a_{n})=
\sum_{\pi_{0}\in {\rm NC}_{{\rm irr}}(n, \ell)}\beta_{\pi_0}[a_1,\ldots, a_n],
\]
where $a_1, \ldots, a_n\in {\mathcal A}_1\cup {\mathcal A}_2$, $\ell=i_1 \cdots i_n$, and ${\rm NC}_{{\rm irr}}(n, \ell)$ consists of partitions in ${\rm NC}_{{\rm irr}}(n)$ whose labels alternate along saturated chains of blocks.
\end{Corollary}

Nevertheless, our goal in this paper is to derive formulas for {\it Motzkin cumulants} of orthogonal replicas since we would like 
to check whether they can be viewed as `partial free cumulants'. For that purpose, let us prove some lemmas on 
the cumulants $B_{w}$.

\begin{Lemma}
Let $a_1, \ldots , a_n$ be orthogonal replicas with labels $\ell=i_1\cdots i_n$ and colors $w=j_1\cdots j_n\in \mathpzc{AM}$.
\begin{enumerate}
\item
The family $\{B_v:v\in \mathpzc{B}(w)\}$ is dedicated to $(a_1, \ldots, a_n)$.
\item
If $h(w)=j$, then 
\[
B_w(a_{1},\ldots ,a_{k}p_{j},a_{k+1}, \ldots, a_{n})=0
\]
for any $1\leq k <n$.
\item
If $h(w)=j_k=j_{k+1}=j$ and $i_k\neq i_{k+1}$, then 
\[
B_{w}(a_1,\ldots, a_kp_{j+1}, a_{k+1}, \ldots, a_n)=0.
\]
\end{enumerate}
These properties remain valid when $a_i$ is replaced by $a_ib_i$, $b_i\in \mathcal{B}$, 
for $1\leq i<n$.
\end{Lemma}

{\it Proof.}
To prove (1), we verify both conditions of Definition 5.2.
The second condition is clearly satisfied since the values of $B_{w}$ are proportional to $p_j$ when 
$h(w)=j$ and $p_j\perp p_{j'}$ when $j\neq j'$. 
For the first condition, it suffices to show that
\[
B_{w}(a_{1},\ldots, a_kp_j, a_{k+1}, \ldots ,a_n)=0
\]
for any $a_1=a_1(j_1), \ldots, a_n=a_n(j_n)$ 
whenever $h(w)=j$ and $1\leq k<n$. 
If $n=2$, we clearly have 
\[
B_{jj}(a_1p_j, a_2)=E(a_1p_ja_2)-E(a_1)E(a_2)=0
\] 
since
\[
E(a_1p_j)=E(a_1)p_j=E(a_1).
\]
Now, let us recall the factorization property of Lemma 7.9:
\[
E(a_1 \cdots a_{k}p_ja_{k+1}\cdots a_n)=E(a_1 \cdots a_{k})E(a_{k+1} \cdots a_n).
\]
Suppose now that the required vanishing of $B_{w'}$ holds for $|w'|<n$. 
If  $w_1=j_1\cdots j_{k}$ and $w_2=j_{k+1}\cdots j_{n}$, where 
$j_k=j_{k+1}=j$, then, by Definition 4.2, we have
\[
E(a_1 \cdots a_{k}p_ja_{k+1}\cdots a_n)=B_{w}(a_1, \ldots, a_{k}p_j, a_{k+1}, \ldots, a_n)
\]
\[
+\sum_{\pi_1\in {\mathcal Int}(w_1)}\sum_{\pi_2\in {\mathcal Int}(w_2)}
B_{\pi_1}(a_1, \ldots, a_{k})B_{\pi_2}(a_{k+1}, \ldots, a_n),
\]
since the remaining $\pi\in {\mathcal Int}(w)$ give zero contribution by the inductive assumption. Thus,
we obtain that $B_w(a_1, \ldots, a_{k}p_j, a_{k+1}, \ldots, a_n)=0$. Clearly,
if $j_k\neq j$ or $j_{k+1}\neq j$, we arrive at the same conclusion since either $E(a_1\cdots a_k)=0$ or 
$E(a_{k+1}\cdots a_n)=0$ and, moreover, the double sum in the above equation vanishes.
This proves that $\{B_{v}:v\in \mathpzc{B}(w)\}$ is dedicated to $(a_1, \ldots, a_n)$.
In particular, we obtain (2). To prove (3), it suffices to observe that
if $i_k\neq i_{k+1}$ and $j_k=j_{k+1}=j$, then
\[
E(a_1\cdots a_kp_{j+1}a_{k+1}\cdots a_n)=0
\]
by Proposition 6.4(3) and therefore, with the use of Definition 4.2 we obtain (3) by induction.
If we replace $a_i$ by $a_ib_i$ for $1\leq i <n$, the arguments used above remain valid.
Thus our proof is complete.
\hfill $\boxvoid$\\

\begin{Lemma}
Let the assumptions of Lemma 8.7 be satisfied and let $h(w)=j$.
\begin{enumerate}
\item
If the replicas have identical labels and $(j_k,j_{k+1})=(j,j)$, then 
\[
B_{w}(a_1,\ldots, a_kp_{j+1}, a_{k+1}, \ldots, a_n)=B_{w}(a_1,\ldots, a_n)
\]
for any $1\leq k <n$.
\item
If $i_k\neq i_{k+1}$ and $(j_k,j_{k+1})\in \{(j,j+1), (j+1,j)\}$, then
\[
B_{w}(a_1,\ldots, a_kp_{j+1}, a_{k+1}, \ldots, a_n)=B_{w}(a_1,\ldots, a_n)
\]
for any $1\leq k <n$.
\end{enumerate}
Both equations remain valid when $a_i$ is replaced by $a_ib_i$, $b_i\in \mathcal{B}$, for
$1\leq i <n$.
\end{Lemma}
{\it Proof.}
If the labels are identical and the colors are not, then
$
B_{w}(a_1, \ldots, a_n)=0
$
by Theorem 8.4.
Moreover, if we replace any $a_i$ by $a_ib_i$, this cumulant still vanishes since 
the property that $p$ meets $p^{\perp}$ or $q$ meets $q^{\perp}$ remains true if $p$ or $q$, respectively, 
is placed between them. Therefore, the cumulant on the LHS in (1) also vanishes.
Thus, let us assume that the labels are identical and all colors are equal to $j$ (if $w\neq j^n$, then
both sides vanish by the definition of $E$). 
If $w=j^2$ and $i_1=i_2$, then 
\[
B_{j^2}(a_1p_{j+1},a_2)=E(a_1p_{j+1}a_2)=E(a_1a_2)-E(a_1)E(a_2)=B_{j^2}(a_1,a_2)
\]
since $E(a_1p_{j+1})=0$.
Let us assume now that (1) holds for any cumulant of order smaller than $n$, where $n>1$. 
We will show that it holds for any cumulant of order $n$. Using Definition 4.2 and the fact that $B_{v}(c_1, \ldots, c_kp_{j+1})=0$
when $v\in \mathpzc{AM}$ and $h(v)=j$, we obtain
\begin{eqnarray*}
B_{w}(a_1,\ldots,a_kp_{j+1}, a_{k+1}, \ldots, a_n)&=&E(a_1\cdots a_kp_{j+1}a_{k+1}\cdots a_n)\\
&&\!\!\!\!\!\!\!-
\sum_{\pi\in {\mathcal Int}_{k}(w)}
B_{\pi}[a_1, \ldots,a_kp_{j+1}, a_{k+1}, \ldots, a_n]
\end{eqnarray*}
where ${\mathcal Int}_{k}(w)$ is the subset of ${\mathcal Int}(w)$ 
consisting of partitions that have at least two blocks and $k$ is not the endpoint of any block.
Now, if $\pi\in {\mathcal Int}_{k}(w)$, then $B_{\pi}[a_1, \ldots,a_kp_{j+1}, a_{k+1}, \ldots, a_n]$ is a product over blocks of
$w$-Boolean cumulants, one of which contains $p_{j+1}$ in the middle and thus
$p_{j+1}$ can be deleted by the inductive assumption. Therefore,
\begin{eqnarray*}
B_{\pi}[a_1, \ldots,a_kp_{j+1}, a_{k+1}, \ldots, a_n]&=&
B_{\pi}[a_1, \ldots, a_n]
\end{eqnarray*}
On the other hand, applying Lemma 7.11 and then Lemma 7.9, we obtain
\begin{eqnarray*}
E(a_1\cdots a_kp_{j+1}a_{k+1}\cdots a_n)&=&E(a_1\cdots a_kp_{j}^{\perp}a_{k+1}\cdots a_n)\\
&=&E(a_1\cdots a_n)-E(a_1\cdots a_k)E(a_{k+1}\cdots a_n)\\
&=&E(a_1\cdots a_n)
-
\sum_{\pi'\in {\mathcal Int}(w')}\sum_{\pi''\in {\mathcal Int}(w'')}\\
&&B_{\pi'}[a_1, \ldots, a_k]B_{\pi''}[a_{k+1}, \ldots, a_n],
\end{eqnarray*}
where $w'=j^k$ and $w''=j^{n-k}$.
Summarizing, we arrive at
\begin{eqnarray*}
B_{w}(a_1,\ldots,a_kp_{j+1}, a_{k+1}, \ldots, a_n)&=&E(a_1\cdots a_n)-
\sum_{\pi\in {\mathcal Int}(w)\setminus \hat{1}_{w}}
B_{\pi}[a_1,\ldots, a_n]\\
&=& B_{w}(a_1,\ldots,a_n)
\end{eqnarray*}
for identical labels and $w=j^n$, which completes the proof of (1).
If $i_k\neq i_{k+1}$ and $(j_k,j_{k+1})\in \{(j,j+1),(j+1,j)\}$, 
it suffices to use Definition 4.2, Lemma 7.12 and induction to obtain (2).
Let us also note that if we replace $a_i$ by $a_ib_i$ for any $1\leq i<n$
and any $b_i\in \mathcal{B}$, the arguments used above remain valid.
Thus our proof is complete.
\hfill $\boxvoid$

\section{Motzkin cumulants of orthogonal replicas}

We are ready to investigate the mixed Motzkin cumulants
$K_{w}(a_1, \ldots, a_n)$, where $a_1, \ldots, a_n$ are orthogonal replicas 
with labels and colors encoded by $\ell=i_1 \cdots i_n$ and
$w=j_1\cdots j_n\in \mathpzc{AM}$, respectively.

\begin{Lemma}
Let the assumptions of Lemma 8.7 be satisfied and let $h(w)=j$.
\begin{enumerate}
\item
If the labels are arbitrary, then
\begin{eqnarray*}
K_w(a_1, \ldots,a_{k}p_j,a_{k+1}, \ldots, a_n)\;\;\;&=&0
\end{eqnarray*}
for any $1\leq k <n$.
\item
If the labels are identical and $(j_{k},j_{k+1})\in \{(j,j),(j,j+1),(j+1,j)\}$, then
\begin{eqnarray*}
K_w(a_1, \ldots, a_kp_{j+1}, a_{k+1}, \ldots, a_n)&=&K_{w}(a_1, \ldots, a_n)
\end{eqnarray*}
for any $1\leq k <n$. 
\end{enumerate}
Both equations remain valid when $a_i$ is replaced by $a_ib_i$, $b_i\in \mathcal{B}$, where
$1\leq i <n$.
\end{Lemma}
{\it Proof.}
By Lemma 8.7(2), we have $B_w(a_1, \ldots, a_kp_{j}, a_{k+1}, \ldots, a_n)=0$, hence
\begin{eqnarray*}
K_w(a_1, \ldots, a_kp_j, a_{k+1}, \ldots, a_n)&=&
-\sum_{\stackrel{\pi \in \mathcal{NC}_{{\rm irr}}(w)}{\scriptscriptstyle \pi\neq \hat{1}_{w}}}
K_{\pi}(a_1, \ldots, a_kp_{j}, a_{k+1}, \ldots, a_n).
\end{eqnarray*}
If $w=j^n$, then
$\mathcal{NC}_{{\rm irr}}(w)=\{\hat{1}_{w}\}$ and the right-hand side vanishes. 
Let $w\neq j^n$ and assume that Motzkin cumulants 
of the considered type (with $a_kp_j$ in the middle) of orders smaller than $n$ vanish. 
Then, the right-hand side also vanishes since no block of $\pi\in \mathcal{NC}_{{\rm irr}}(w)$ that gives a nonzero $K_{\pi}$ on the RHS
can end at $k$ since in that case it would be an inner block contributing a cumulant proportional to $p_j$, which must vanish
since the family $\{B_{v}:v\in \mathpzc{B}(w)\}$ is dedicated to $(a_1, \ldots, a_n)$ by Lemma 8.7(1).
This proves (1). 

In order to prove (2), note that if $h(w)=j$ and $(j_{k},j_{k+1})\in \{(j,j),(j,j+1),(j+1,j)\}$, 
then we obtain
\begin{eqnarray*}
K_{w}(a_1, \ldots, a_kp_{j+1},a_{k+1}, \ldots, a_n)&=&
B_{w}(a_1, \ldots, a_{k}p_{j+1}, a_{k+1}, \ldots, a_n)\\
& -&\!\!\!\!\!\sum_{\stackrel{\pi\in \mathcal{NC}_{{\rm irr}}(w)}{\scriptscriptstyle \pi\neq \hat{1}_{w}}}
K_{\pi}[a_1, \ldots, a_kp_{j+1},a_{k+1}, \ldots, a_n]\\
&=&
B_{w}(a_1, \ldots, a_n)\\
& -&\!\!\!\!\!\sum_{\stackrel{\pi\in \mathcal{NC}_{{\rm irr}}(w)}{\scriptscriptstyle \pi\neq \hat{1}_{w}}}
K_{\pi}[a_1, \ldots, a_k,a_{k+1}, \ldots, a_n]
\\
&=&K_{w}(a_1, \ldots, a_n)
\end{eqnarray*}
by Lemma 8.8(1) since no block of $\pi$ can end at $k$ and thus $a_kp_{j+1}$ is 
in the middle of the arguments of one cumulant $K_{v}$ contributing to $K_{\pi}$ and $h(v)=j$. Therefore, by induction,
we can delete this $p_{j+1}$. This completes the proof.
\hfill $\boxvoid$\\

\begin{Theorem}
If $a_1, \ldots, a_n$ are orthogonal replicas with identical labels and colors encoded by $w\in \mathpzc{AM}$, then 
\[
K_w(a_1, \ldots , a_n)=\left(\sum_{\pi\in \mathcal{M}_{{\rm irr}}(w)}(-1)^{|\pi|-1}\hat{\beta}_{\pi}[a_1,\ldots, a_n] \right)p_j
\]
where $h(w)=j$. If their labels are not identical, then $K_{w}(a_1, \ldots, a_n)=0$.
\end{Theorem}
{\it Proof.}
First, assume that the replicas $a_1,\ldots,a_n$ have identical labels.
By Lemma 8.7(1), the family
$\{B_v:v\in \mathpzc{B}(w)\}$ is dedicated to $(a_1,\ldots,a_n)$.
Hence Theorem 5.3 can be applied. Moreover, by Theorem 8.4, the cumulants
$B_v(a_{k_1},\ldots,a_{k_p})$ vanish whenever the labels of
$a_{k_1},\ldots,a_{k_p}$ are identical and the corresponding colors are not
all equal. Let us note that this vanishing remains valid if
$a_{k_1},\ldots,a_{k_{p-1}}$ are multiplied on the right by elements of
$\mathcal{B}$, by the $\mathcal B$-multilinearity properties of the cumulants
and the same orthogonality argument. Therefore, the summation over
$\mathcal{NC}_{{\rm irr}}(w)$ in Theorem 5.3 reduces to the summation over
$\mathcal{M}_{{\rm irr}}(w)$. Finally, using Lemma 8.8(1) repeatedly, we get
\[
B_{\pi}[a_1,\ldots,a_n]
=
\hat{\beta}_{\pi}[a_1,\ldots,a_n]p_j
\]
for every $\pi\in\mathcal{M}_{{\rm irr}}(w)$. This proves the assertion in
the case of identical labels.

Now assume that the labels are not identical. 
We distinguish three cases. 

{\it Case 1.} Let $w=j^n$ for some $j\in \mathbb{N}$. Then
\[
K_w(a_1,\ldots,a_n)
=
B_w(a_1,\ldots,a_n)
=
\hat{\beta}_{n}(a_1,\ldots,a_n)p_j
=
0,
\]
since mixed Boolean cumulants of replicas with non-identical labels vanish.

{\it Case 2.} Let $w\neq j^n$ and suppose that the labels of the replicas
$a_1,\ldots,a_n$ alternate when their colors increase, by which we mean that
any two consecutive colors appearing in $w$ correspond to replicas with different
labels. We argue by induction on the length of $w$, using Case 1 for constant
subwords.

The first nonconstant case is $w=j(j+1)j$. By Definition 4.5,
\[
B_{j(j+1)j}(a_1,a_2,a_3)
=
K_{j(j+1)j}(a_1,a_2,a_3)
+
K_{j^2}(a_1K_{j+1}(a_2),a_3).
\]
On the other hand, by Theorem 8.4, the only irreducible monotonically adapted
partition which may give a non-zero contribution to
$B_{j(j+1)j}(a_1,a_2,a_3)$ is the partition with outer block $\{1,3\}$ and
inner singleton $\{2\}$. The one-block partition gives zero, since its block
contains colors which are not all equal. Hence, using Lemma 9.1(2), we obtain
\[
B_{j(j+1)j}(a_1,a_2,a_3)
=
K_{j^2}(a_1K_{j+1}(a_2),a_3)
\]
and hence
\[
K_{j(j+1)j}(a_1,a_2,a_3)=0.
\]

Assume now that the assertion has already been proved for all Motzkin words of
length smaller than $n$, constant words being covered by Case 1, and let
$|w|=n$. By Theorem 8.4 and repeated use of Lemma 9.1(2), we have
\[
B_w(a_1,\ldots,a_n)
=
\sum_{\pi\in \mathcal{M}_{{\rm irr}}(w,\ell)}
\hat{\beta}_{\pi}[a_1,\ldots,a_n]p_j
=
\sum_{\pi\in \mathcal{M}_{{\rm irr}}(w,\ell)}
K_{\pi}[a_1,\ldots,a_n].
\]
On the other hand, by Definition 4.5,
\[
B_w(a_1,\ldots,a_n)
=
K_w(a_1,\ldots,a_n)
+
\sum_{\substack{\pi\in \mathcal{NC}_{{\rm irr}}(w,\ell)\\
\pi\neq \hat{1}_{w}}}
K_{\pi}[a_1,\ldots,a_n].
\]
Using the inductive hypothesis for the proper blocks which occur in
$K_{\pi}$, we obtain
\[
K_w(a_1,\ldots,a_n)
=
-
\sum_{\pi\in
\mathcal{NC}_{{\rm irr}}(w,\ell)\setminus
\mathcal{M}_{{\rm irr}}(w,\ell)}
K_{\pi}[a_1,\ldots,a_n].
\]
We now show that the set on the right-hand side is empty. Let
$\pi\in\mathcal{NC}_{{\rm irr}}(w,\ell)$. If $\pi$ were not monotonically
adapted to $\ell$, then there would exist a saturated chain of blocks of
$\pi$ along which the corresponding labels do not alternate. This contradicts
the assumption that the labels of replicas alternate whenever their colors
increase. Therefore every partition in
$\mathcal{NC}_{{\rm irr}}(w,\ell)$ is monotonically adapted, and hence
\[
\mathcal{NC}_{{\rm irr}}(w,\ell)
=
\mathcal{M}_{{\rm irr}}(w,\ell).
\]
Consequently,
\[
\mathcal{NC}_{{\rm irr}}(w,\ell)\setminus
\mathcal{M}_{{\rm irr}}(w,\ell)=\emptyset,
\]
and thus $K_w(a_1,\ldots,a_n)=0$.

{\it Case 3.}
Let $w\neq j^n$ and assume that the labels of the replicas
$a_1,\ldots,a_n$ are not identical and do not alternate when the colors
increase. Thus, there exists a label, say $i_0$, such that the subword of
$w$ corresponding to the replicas with label $i_0$ contains the two
consecutive colors $l$ and $l+1$. Let
\[
S=\{r_1<\cdots<r_m\}\subset \{1,\ldots,n\}
\]
be the set of positions at which the label is $i_0$, and put
\[
w_0=j_{r_1}\cdots j_{r_m}.
\]
Then $w_0$ is not a constant word, since both colors $l$ and $l+1$ occur
in $w_0$.

In this situation no irreducible partition can be monotonically adapted to both
\(w\) and \(\ell\); equivalently,
\[
\mathcal M_{\rm irr}(w,\ell)=\emptyset.
\]
Indeed, in an irreducible partition belonging to \(\mathcal M(w,\ell)\),
there is one outer block of the lowest color, and every block of color
\(r+1\) lies below a block of color \(r\). Since the labels alternate along
chains and there are only two labels, the label is then forced to change along
each chain whenever the color increases by one. Hence a fixed label cannot occur at two
consecutive colors.

Thus Theorem 8.4 gives \(B_w(a_1,\ldots,a_n)=0\).
Hence, by the definition of Motzkin cumulants,
\[
0
=
K_w(a_1,\ldots,a_n)
+
\sum_{\substack{\pi\in \mathcal{NC}_{{\rm irr}}(w)\\
\pi\neq \hat{1}_{w}}}
K_{\pi}[a_1,\ldots,a_n].
\]
Therefore,
\[
K_w(a_1,\ldots,a_n)
=
-
\sum_{\substack{\pi\in \mathcal{NC}_{{\rm irr}}(w)\\
\pi\neq \hat{1}_{w}}}
K_{\pi}[a_1,\ldots,a_n].
\]
We now apply the induction hypothesis. If a block of $\pi$ contains replicas
with non-identical labels, then the corresponding Motzkin cumulant has order
smaller than $n$ and vanishes. Hence only those terms may survive for which
every block of $\pi$ contains replicas with one fixed label.

We group the remaining partitions $\pi$ according to their complementary
part, by which we mean their restriction to the complement of $S$, together
with the nesting of this restriction relative to the positions from $S$.
Let $\sigma$ denote such a fixed complementary part. For this fixed
$\sigma$, the cumulants corresponding to the complementary blocks are fixed.
If they are nonzero, they may be factored out, while their nesting may leave
boundary projections attached to the arguments whose positions belong to
$S$. 

Thus the total contribution of all partitions with this fixed
complementary part $\sigma$ can be grouped according to the induced
partition of the positions in $S$. 
If the color word induced on these positions is not admissible, then no
nonzero nested cumulant contribution is obtained from this induced part, and
the corresponding grouped contribution is zero.
Otherwise, this induced word,
denoted by $v_\sigma$, is an admissible nonconstant word. Moreover, all
arguments corresponding to this word have the same label $i_0$, while the
nesting of the complementary part may only leave boundary projections attached
to these arguments. Hence the grouped contribution is, up to fixed factors
from $\mathcal B$, of the form
\[
\sum_{\rho\in \mathcal{NC}_{{\rm irr}}(v_{\sigma})}
K_{\rho}[a_{r_1}b_1^{(\sigma)},\ldots,a_{r_m}b_m^{(\sigma)}],
\]
where $b_1^{(\sigma)},\ldots,b_m^{(\sigma)}\in \mathcal B$ are the boundary
projections produced by the nesting.
By the defining relation between $v_\sigma$-Boolean cumulants and Motzkin
cumulants, this sum is precisely
\[
B_{v_\sigma}
(a_{r_1}b_1^{(\sigma)},\ldots,a_{r_m}b_m^{(\sigma)}).
\]
Since all arguments of this $v_\sigma$-Boolean cumulant have the same label
$i_0$, whereas $v_\sigma$ is not constant, this expression vanishes by
Theorem 8.4 and $\mathcal B$-multilinearity.

Thus the contribution corresponding to the fixed complementary part
$\sigma$ is zero. Since this holds for every $\sigma$, the whole remaining
contribution vanishes. Consequently, all terms on the right-hand side of the
above formula for $K_w(a_1,\ldots,a_n)$ vanish, and therefore
\[
K_w(a_1,\ldots,a_n)=0.
\]
This completes the proof of Case 3, and hence the induction. \hfill $\boxvoid$\\

\begin{Example}
{\rm
Let $w=12^21$ and suppose we have replicas with identical labels,
where we use the notation of Example 3.9 for partitions and that of Example 8.2 for replicas.
This simple example is instructive because it shows how the summation over 
\[
{\mathcal NC}_{{\rm irr}}(w)=\{(\pi_A,w), (\pi_B,w), (\pi_C,w), (\pi_D,w), (\pi_E, w)\},
\]
which has its origin in Theorem 5.3, reduces to that over 
\[
\mathcal{M}_{{\rm irr}}(w)=\{(\pi_{C},w), (\pi_E,w)\}
\] 
in Theorem 9.2, with the use of Lemma 9.1. First, observe that $B_{121}(x,x',x'')=0$ 
for any replicas $x,x',x''$ of label $1$ and colors encoded by $121$, hence
\[
K_{121}(x,x',x'')=-\hat{\beta}_1(x')K_{11}(xp_2,x'')=-\hat{\beta}_1(x')\hat{\beta}_{2}(x,x'')p_1.
\]
Similarly,
\begin{eqnarray*}
K_w(x_{k_1},x_{k_2},x_{k_3},x_{k_4})
&=&
-\hat{\beta}_1(x_{k_2})K_{121}(x_{k_1}p_2,x_{k_3},x_{k_4})-
\hat{\beta}_1(x_{k_3})K_{121}(x_{k_1},x_{k_2}p_2,x_{k_4})\\
&& \!\!\!\!\! -\; \hat{\beta}_2(x_{k_2},x_{k_3})K_{11}(x_{k_1}p_2,x_{k_4})
-
\hat{\beta}_1(x_{k_2})\hat{\beta}_1(x_{k_3})K_{11}(x_{k_1}p_2, x_{k_4})\\
&=&
\hat{\beta}_2(x_{k_1},x_{k_4})(-\hat{\beta}_2(x_{k_2},x_{k_3})+\hat{\beta}_{1}(x_{k_2})\hat{\beta}_{1}(x_{k_3}))p_1\\
&=&
-\hat{\beta}_{\pi_{C}}[x_{k_1}, x_{k_2}, x_{k_3}, x_{k_4}]p_1+\hat{\beta}_{\pi_{E}}[x_{k_1}, x_{k_2}, x_{k_3}, x_{k_4}]p_1,
\end{eqnarray*}
by Definition 4.2, where we again use Theorem 8.4 and then Lemma 9.1(2).
Of course, this result agrees with the formula of
Theorem 9.2 for replicas with identical labels, proved with the use of Lemma 8.8.
}
\end{Example}

\begin{Example}
{\rm 
Let $w=1^22^21$ and consider the corresponding Motzkin cumulants for two
choices of mixed labels: one in which the labels of replicas alternate when
their colors increase and one in which they do not alternate. Keeping the
notation of Example 8.2 for replicas, we have
\begin{eqnarray*}
K_w(x_{k_1},x_{k_2},y_{k_3},y_{k_4},x_{k_5})&=&B_w(x_{k_1},x_{k_2},y_{k_3},y_{k_4},x_{k_5})\\
&&\!\!\!\!\!-\;\hat{\beta}_{2}(y_{k_3},y_{k_4})K_{111}(x_{k_1},x_{k_2}p_2,x_{k_5})\\
&&\!\!\!\!\! -\; \hat{\beta}_1(y_{k_3})\hat{\beta}_{1}(y_{k_4})K_{111}(x_{k_1},x_{k_2}p_2,x_{k_5}),\\
K_w(x_{k_1},x_{k_2},y_{k_3},x_{k_4},x_{k_5})&=&B_w(x_{k_1},x_{k_2},y_{k_3},x_{k_4},x_{k_5})\\
&&\!\!\!\!\!-\;\hat{\beta}_1(y_{k_3})K_{1121}(x_{k_1},x_{k_2}p_2,x_{k_4},x_{k_5})\\
&&\!\!\!\!\! -\; \hat{\beta}_1(y_{k_3})\hat{\beta}_{1}(x_{k_4})K_{111}(x_{k_1},x_{k_2}p_2,x_{k_5}),
\end{eqnarray*}
where we use the vanishing of Motzkin cumulants with mixed labels of orders
smaller than $5$. By Theorem 8.4, we have
\begin{eqnarray*}
B_w(x_{k_1},x_{k_2},y_{k_3},y_{k_4},x_{k_5})&=&
\hat{\beta}_{3}(x_{k_1},x_{k_2},x_{k_5})\hat{\beta}_{2}(y_{k_3},y_{k_4})p_1\\
&&\!\!\!\!\!\!+\;\hat{\beta}_{3}(x_{k_1}, x_{k_2}, x_{k_5})\hat{\beta}_{1}(y_{k_3})\hat{\beta}_{1}(y_{k_4})p_1,\\
B_w(x_{k_1},x_{k_2},y_{k_3},x_{k_4},x_{k_5})&=&0.
\end{eqnarray*}
On the other hand, by Lemma 9.1(2) and Theorem 9.2, we have
\begin{eqnarray*}
K_{111}(x_{k_1},x_{k_2}p_2,x_{k_5})\;\;\;\;\;\;\;\,&=&\;\;\;
\hat{\beta}_{3}(x_{k_1},x_{k_2},x_{k_5})p_1,\\
K_{1121}(x_{k_1},x_{k_2}p_2,x_{k_4},x_{k_5})&=&
-\hat{\beta}_{1}(x_{k_4})\hat{\beta}_{3}(x_{k_1},x_{k_2}, x_{k_5})p_1.
\end{eqnarray*}
In the second computation,
the last two terms combine into
\[
-\hat{\beta}_1(y_{k_3})
B_{1121}(x_{k_1},x_{k_2}p_2,x_{k_4},x_{k_5}),
\]
which vanishes by Theorem 8.4 and $\mathcal B$-multilinearity.
Therefore, both considered mixed cumulants vanish. 
}
\end{Example}

\begin{Corollary}
Let $I$ be a set of indices and let $\{x_{k}\}_{k\in I}$ and $\{y_{k}\}_{k\in I}$ be families of orthogonal replicas with
labels $1$ and $2$, respectively. Then
\[
K_w(x_{k_1}+y_{k_1}, \ldots , x_{k_n}+y_{k_n})=K_w(x_{k_1}, \ldots, x_{k_n})+ K_w(y_{k_1}, \ldots, y_{k_n})
\]
for any $k_1, \ldots, k_n\in I$ whenever $w=j_1\cdots j_n\in \mathpzc{AM}$ encodes the colors of 
$x_{k_1}, \ldots, x_{k_n}$ and $y_{k_1}, \ldots, y_{k_n}$, respectively.
\end{Corollary}
{\it Proof.}
This property follows from the vanishing of Motzkin cumulants with different labels shown in 
Theorem 9.2.
\hfill $\boxvoid$\\

It is worth collecting the already proven results on $w$-Boolean cumulants and Motzkin cumulants and 
observe an interplay between labels and colors.

\begin{Corollary}
Let $a_1, \ldots, a_n$ be orthogonal replicas 
with labels and colors encoded by 
$\ell=i_1 \cdots i_n$ and $w=j_1\cdots j_n\in \mathpzc{AM}$, respectively. 
\begin{enumerate}
\item
If $i_1, \ldots, i_n$ are not identical, then 
\[
K_{w}(a_1, \ldots, a_n)=0.
\]
\item
If $i_1=\cdots = i_n$ and $j_1, \ldots, j_n$ are not identical, then 
\[
B_w(a_1, \ldots, a_n)=0.
\]
\item
If $i_1=\cdots = i_n$ and $j_1=\cdots = j_n$, then 
\[
K_{w}(a_1, \ldots, a_n)=B_{w}(a_1, \ldots , a_n).
\]
\end{enumerate}
\end{Corollary}
{\it Proof.}
These properties follow from Theorem 9.2, Theorem 8.4 and Remark 4.6(3), respectively.
\hfill $\boxvoid$

\begin{Example}
{\rm
If $i_1=\cdots = i_n$ and $j_1, \ldots, j_n$ are not identical, then it is typical that
\[
E(a_1\cdots a_n)=0\;\;\;{\rm and}\;\;\;K_w(a_1, \ldots, a_n)\neq 0
\] 
for $h(w)=j$, which is not the case for other types of cumulants. 
The simplest example is when $w=121$ and $x_{k_1},x_{k_2},x_{k_3}$ are orthogonal replicas of label $1$
and colors $1,2,1$, respectively. Then $E_1(x_{k_1}x_{k_2}x_{k_3})=0$,
whereas
\[
K_{w}(x_{k_1},x_{k_2},x_{k_3})=-\hat{\beta}_{2}(x_{k_1},x_{k_3})\hat{\beta}_{1}(x_{k_2})p_1
\]
by Theorem 9.2. In fact, this cumulant corresponds to the `missing' expression in the formula for the free 
cumulant $r_3(x_{k_1}, x_{k_2}, x_{k_3})$ in the third equation of Example 2.7.
}
\end{Example}

Finally, let us look at cumulants in which some of the arguments are orthogonal projections $p_j$. Such cumulants have not 
appeared in our computations yet since Motzkin cumulants are nested, but they are closely related to the well-known property saying that
free cumulants with some arguments equal to the unit vanish.

\begin{Lemma}
Let the assumptions of Lemma 8.7 be satisfied.
\begin{enumerate}
\item
For any $j,j_1$ it holds that
$B_{j_1}(p_j)=K_{j_1}(p_j)=p_j$.
\item
For any $n\geq 1$, $0\leq k\leq n$ and any $j$, it holds that
\[
B_{\hat w}(a_1,\ldots,a_k,p_j,a_{k+1},\ldots,a_n)=0
\]
whenever $h(w)=j$ and
$
\hat w=j_1\cdots j_k\,j\,j_{k+1}\cdots j_n\in \mathpzc{AM}.
$
\item
For any $n\geq 1$, $0\leq k\leq n$ and any $j$, it holds that
\[
K_{\hat w}(a_1,\ldots,a_k,p_j,a_{k+1},\ldots,a_n)=0
\]
whenever
$
\hat w=j_1\cdots j_k\,j\,j_{k+1}\cdots j_n\in \mathpzc{AM}.
$
\end{enumerate}
Both equations remain valid when $a_i$ is replaced by
$a_i b_i$, where $b_i\in\mathcal B$ and $1\leq i<n$.
\end{Lemma}

{\it Proof.}
By Definitions 4.2, 4.5 and 7.1,
\[
K_{j_1}(p_j)=B_{j_1}(p_j)=E(p_j)=p_j,
\]
which proves (1). 
Next, for $\hat w=j_1j$, we have
\[
B_{\hat w}(a_1,p_j)=E(a_1p_j)-E(a_1)E(p_j)=0
\]
and similarly
\[
B_{\hat w}(p_j,a_1)=0
\]
for $\hat w=jj_1$. By induction over interval partitions, we then obtain
\[
B_{\hat w}(p_j,a_1,\ldots,a_n)=
B_{\hat w}(a_1,\ldots,a_n,p_j)=0
\]
whenever $\hat w=jj_1\cdots j_n\in\mathpzc{AM}$ or
$\hat w=j_1\cdots j_nj\in\mathpzc{AM}$, respectively. Using Definition 4.5
and induction over irreducible noncrossing partitions, we also obtain
\[
K_{\hat w}(p_j,a_1,\ldots,a_n)=
K_{\hat w}(a_1,\ldots,a_n,p_j)=0.
\tag{9.1}
\]

We shall also use the following consequence of the recursive definition of
$w$-Boolean cumulants. If $1\leq k<n$ and $\hat w=j_1\cdots j_k\,j\,j_{k+1}\cdots j_n\in\mathpzc{AM}$,
then
\[
B_{\hat w}(a_1,\ldots,a_k,p_j,a_{k+1},\ldots,a_n)
=
B_w(a_1,\ldots,a_kp_j,a_{k+1},\ldots,a_n).
\tag{9.2}
\]
Indeed, this follows by induction on $|w|$. After subtracting the two defining
moment-cumulant formulas, the moment terms cancel. The proper interval terms
cancel pairwise: if the block containing $p_j$ is non-singleton, this follows
from the induction hypothesis or from the endpoint case (9.1), while if
$\{k+1\}$ is a singleton block, then it contributes $B_j(p_j)=p_j$, which is
absorbed into the argument immediately to its left.

We now prove (2). The endpoint cases follow from (9.1). If $p_j$ is in the
middle, then (9.2) and Lemma 8.7 give
\[
B_{\hat w}(a_1,\ldots,a_k,p_j,a_{k+1},\ldots,a_n)=0,
\]
since $h(w)=j$. Thus (2) is proved. The assumption $h(w)=j$ is essential in
this part.

It remains to prove (3). The endpoint cases follow from (9.1), so assume that
$p_j$ is in the middle. We argue by induction on $|\hat w|$. By Definition 4.5,
\[
\begin{aligned}
K_{\hat w}(a_1,\ldots,a_k,p_j,a_{k+1},\ldots,a_n)
&=
B_{\hat w}(a_1,\ldots,a_k,p_j,a_{k+1},\ldots,a_n)  \\
&\quad -
\sum_{\substack{\hat\pi\in\mathcal{NC}_{\rm irr}(\hat w)\\
\hat\pi\neq \hat 1_{\hat w}}}
K_{\hat\pi}[a_1,\ldots,a_k,p_j,a_{k+1},\ldots,a_n].
\end{aligned}
\]
By the induction hypothesis, together with the endpoint cases (9.1), all terms
in the sum vanish except those for which $\{k+1\}$ is a singleton block.
Removing this singleton and absorbing
$K_j(p_j)=p_j$
into the argument immediately to its left gives a partition
$\pi\in\mathcal{NC}_{\rm irr}(w)$, and every such $\pi$ is obtained in this
way. Therefore the surviving part of the sum is
\[
\sum_{\pi\in\mathcal{NC}_{\rm irr}(w)}
K_\pi[a_1,\ldots,a_kp_j,a_{k+1},\ldots,a_n]
=
B_w(a_1,\ldots,a_kp_j,a_{k+1},\ldots,a_n).
\]
Substituting this into the defining formula above, we obtain exactly the
difference between the two sides of (9.2), hence zero. This proves (3).
Finally, the version with $a_i$ replaced by $a_i b_i$ follows from the
$\mathcal B$-bimodule property of $E$ and from the corresponding bimodule
property of the cumulants proved in Remark 4.6.
\hfill $\boxvoid$

\section{Decomposition of free cumulants}

Relations between Motzkin cumulants and $w$-Boolean cumulants resemble those between
free cumulants and Boolean cumulants. This indicates that there should 
be some relation between free cumulants and Motzkin cumulants especially since each  
$\mathcal{NC}_{{\rm irr}}(w)$ is isomorphic to a sublattice of ${\rm NC}_{{\rm irr}}(n)$. We will show that
free cumulants of random variables
can be decomposed in terms of `scalar-valued Motzkin cumulants'. 

Let us first write the decomposition formula for the mixed moments of free random variables of Theorem 6.7 in an operator-valued version.

\begin{Theorem} 
Under the assumptions of Theorem 6.7, the mixed moments of free random variables have the decomposition
\[
\left(\left(\varphi_1\star \varphi_2\right)\circ \tau\right)(a_1\cdots a_n)=
\sum_{w=j_1\cdots j_n\in \mathpzc{M}_{n}}\zeta\left(E(a_{1}(j_1)\cdots a_{n}(j_n))\right),
\]
where $a_1\in {\mathcal A}_{i_1}, \ldots, a_n\in {\mathcal A}_{i_n}$ and $i_1\neq \cdots \neq i_n$, $n\in \mathbb{N}$, and where
$\zeta:{\mathcal B}\rightarrow {\mathbb C}$ is the normalized linear functional defined by $\zeta(p_k)=\delta_{k,1}$, 
where $k\in \mathbb{N}$.
\end{Theorem}
{\it Proof.}
This formula follows immediately from Theorem 6.7, Proposition 7.5 and the fact that $\zeta(p_1)=1$.
\hfill $\boxvoid$\\

We would like to find the corresponding decomposition formula for free cumulants.
It can already be seen from Theorem 9.2 that Motzkin cumulants of orthogonal replicas can be associated with fragments of
free cumulants of free random variables. Those with mixed labels vanish and those with identical labels reproduce 
fragments of the second formula of Theorem 2.6 multiplied by a projection.

\begin{Definition}
{\rm
Let $({\mathcal A}_{1}, \varphi_1)$ and $({\mathcal A}_{2}, \varphi_2)$ be noncommutative probability spaces
and let ${\mathcal A}^{\circ}$ be the vector space sum of ${\mathcal A}_{1}$ and ${\mathcal A}_{2}$.
The family of multilinear functionals $\{k_w:w\in \mathpzc{M}\}$ on ${\mathcal A}^{\circ}$ defined as 
the multilinear extensions of
\[
k_w(a_1, \ldots, a_n)=\zeta (K_{w}(a_{1}(j_1), \ldots, a_{n}(j_n))),
\] 
where $a_1, \ldots, a_n\in \mathcal{A}_{1}\cup \mathcal{A}_{2}$, $w=j_1\cdots j_n$, and the Motzkin cumulants
of the orthogonal replicas on the RHS of the above equation are given by Theorem 9.2, will be 
called the family of {\it scalar-valued Motzkin cumulants}.
}
\end{Definition}

\begin{Example}
{\rm In contrast to free cumulants, scalar-valued Motzkin cumulants distinguish the units $1_1$ and $1_2$. 
For instance, using Theorem 9.2, we obtain
\[
k_{w}(x_1,1_1,x_2)=\left\{
\begin{array}{rr}\beta_{2}(x_1,x_2)&{\rm if}\;w=111,\\
-\beta_{2}(x_1,x_2)&{\rm if}\; w=121,
\end{array}
\right.
\]
whereas $k_{w}(x_1,1_2,x_2)=0$ for $w\in \{111, 121\}$ and $x_1,x_2\in \mathcal{A}_{1}$. 
In particular, this implies that 
\[
k_{111}(x_1,1_i,x_2)+k_{121}(x_1,1_i,x_2)=0
\]
for $i\in \{1,2\}$, and since ${\rm LHS}=r_3(x_1,1,x_2)$ by the theorem given below, we obtain the well-known property of free cumulants
$r_3(x_1,1,x_2)=0$, where $1$ is the unit in $\mathcal{A}_{1}\star \mathcal{A}_{2}$.
}
\end{Example}

\begin{Theorem}
Let $a_1, \ldots, a_n\in \mathcal{A}_{1}\cup {\mathcal A}_{2}$ and let $\tau(a_1), \ldots, \tau(a_n)$ 
be their canonical embeddings in ${\mathcal A}_{1}\star {\mathcal A}_{2}$.
Their mixed cumulants associated with $\varphi_1\star \varphi_2$ have the decomposition
\[
r_{n}(\tau(a_1), \ldots, \tau(a_n))=\sum_{w\in \mathpzc{M}_{n}}k_{w}(a_1, \ldots, a_n),
\]
where
\[
k_w(a_1, \ldots, a_n)=\sum_{\pi\in \mathcal{M}_{{\rm irr}}(w)}(-1)^{|\pi|-1}\beta_{\pi}[a_1,\ldots, a_n]
\] 
if all $a_{i}\in \mathcal{A}_{1}$ or all $a_i\in \mathcal{A}_{2}$ and otherwise
we have $k_w(a_1, \ldots, a_n)=0$, where $n\in \mathbb{N}$.
\end{Theorem}
{\it Proof.}
If the variables have identical labels, we have
\begin{eqnarray*}
k_w(a_1, \ldots ,a_n)
&=&\sum_{\pi\in \mathcal{M}_{{\rm irr}}(w)}(-1)^{|\pi|-1}\hat{\beta}_{\pi}[a_1(j_1),\ldots, a_n(j_n)]\\
&=&\sum_{\pi\in \mathcal{M}_{{\rm irr}}(w)}(-1)^{|\pi|-1}\beta_{\pi}[a_1,\ldots, a_n]
\end{eqnarray*}
by Theorem 9.2
and Proposition 7.2.
Moreover, we have a natural bijection
\[
\eta: {\rm NC}_{{\rm irr}}(n)\rightarrow \bigsqcup_{w\in \mathpzc{M}_{n}}{\mathcal M}_{{\rm irr}}(w)
\]
given by
\[
\eta(\{V_1, \ldots, V_p\})=\{(V_1, v_1), \ldots , (V_p, v_p)\},
\] 
where $v_k=d_k^{|V_k|}$ and $d_k$ is the depth of $V_k$ (see Proposition 7.1 in \cite{[L6]}).
Therefore, after summing over $w\in \mathpzc{M}_{n}$ and using the
Boolean-cumulant formula for free cumulants from Theorem 2.6, we obtain
\[
r_n(\tau(a_1),\ldots,\tau(a_n))
=
\sum_{w\in \mathpzc{M}_{n}} k_w(a_1,\ldots,a_n)
\]
whenever the variables have identical labels. Moreover, 
\[
k_w(a_1, \ldots, a_n)=\zeta(K_w(a_1(j_1), \ldots, a_n(j_n)))=0
\] 
if not all variables (and thus not all replicas) have identical labels. Therefore, we obtain the desired decomposition
for any $a_1, \ldots, a_n\in \mathcal{A}_{1}\cup \mathcal{A}_{2}$, which completes the proof.
\hfill $\boxvoid$\\

\begin{Remark}
{\rm 
This result shows that one can construct free cumulants of variables living in the vector space
$\tau({\mathcal A}^{\circ})\subset {\mathcal A}_{1}\star {\mathcal A}_{2}$ by adding Motzkin cumulants of 
their orthogonal replicas.
These free cumulants can be uniquely extended to the free product ${\mathcal A}_1\star {\mathcal A}_{2}$
by the property called c-multiplicativity \cite{[Sp]}.
In consequence, one can construct the free product of functionals $\varphi_1\star\varphi_2$ 
from Motzkin cumulants of orthogonal replicas. The fact that the arguments of these cumulants 
have non-identified units in contrast to the arguments of free cumulants is natural since
the summands include the Boolean cumulants.
}
\end{Remark}

\begin{Example}
{\rm In the univariate case, we set $a_1=\ldots =a_n=a$ in 
Theorem 10.4 and obtain the decomposition of free cumulants $(r_n)$ of the variable $a$ of the form
\[
r_n=\sum_{w\in \mathpzc{M}_{n}}k_w,
\]
where
\[
k_w=
\sum_{\pi\in {\mathcal M}_{{\rm irr}}(w)}
(-1)^{|\pi|-1}\beta_{\pi},
\]
which is a refinement of the formula of Theorem 2.2.
For example, the cumulants given in Example 2.3 can be decomposed as 
\begin{eqnarray*}
r_1&=&k_1\;\;=\;\;\beta_1,\\
r_2&=&k_{11}\;=\;\;\beta_2,\\
r_3&=&k_{111}+k_{121} \;\;=\;\;\beta_3+(-\beta_2\beta_1),\\
r_4&=&k_{1111}+k_{1121}+k_{1211}+k_{1221}\\
&=&\beta_4+(-\beta_3\beta_1)+(-\beta_3\beta_1)+ (-\beta_2^{2}+\beta_2\beta_1^2),
\end{eqnarray*}
where $\beta_n$ is the Boolean cumulant of $a$ of order $n$, and parentheses
are used only to indicate which Boolean polynomial corresponds to a given
$k_w$.
The next free cumulant is $r_5$, for which there are 9 Motzkin paths contributing
the corresponding cumulants $k_w$. These are 
expressed in terms of $\beta_{\pi}\equiv \beta(\pi)$ (from which `Boolean polynomials' can be easily
obtained) shown in Fig.~4. Similar formulas hold for cumulants of any orthogonal
replicas $a_1, \ldots, a_n$, provided they have identical labels.
}
\end{Example}

\begin{figure}
\unitlength=1mm
\special{em:linewidth 1pt}
\linethickness{0.5pt}
\begin{picture}(180.00,140.00)(-20.00,-45.00)

%%%%%%%%%%%%%%%%%%%%%%%%%%%%%%%%%%%%%%%%%%%%%%%%%%%
\put(-05.00,82.00){$w$}
\put(16.00,82.00){$Y_w$}
\put(44.00,82.00){$k_w$}
%%%%%%%%%%%%%%%%%%% 1 %%%%%%%%%%%%%%%%%%%%%%%%%%%%%
%%%%%%%%%%%%% 11111 %%%%%%%%%%%%%%%%%%%%%%%%%%%%%%

\put(-11.00,70.00){\line(1,0){4.00}}
\put(-07.00,70.00){\line(1,0){4.00}}
\put(-02.00,70.00){\line(1,0){4.00}}
\put(01.00,70.00){\line(1,0){4.00}}
\put(-11.00,70.00){\circle*{1.00}}
\put(-07.00,70.00){\circle*{1.00}}
\put(-03.00,70.00){\circle*{1.00}}
\put(01.00,70.00){\circle*{1.00}}
\put(05.00,70.00){\circle*{1.00}}

%%%%%%%%%%%%%%%%%% Y %%%%%%%%%%%%%%%%%%%%

\put(14.00,70.00){\line(1,0){16}}
\put(14.00,74.00){\line(1,0){16}}
\put(14.00,70.00){\line(0,1){4}}
\put(18.00,70.00){\line(0,1){4}}
\put(22.00,70.00){\line(0,1){4}}
\put(26.00,70.00){\line(0,1){4}}
\put(30.00,70.00){\line(0,1){4}}

\put(15.00,71.00){\small{1}}
\put(19.00,71.00){\small{2}}
\put(23.00,71.00){\small{3}}
\put(27.00,71.00){\small{4}}

%%%%%%%%%%%%%%%%%%%%%%%%%%%%%%%%%%%%%%

\put(36.00,71.00){$\beta$}
\put(39.00,71.00){$\big($}
\put(42.00,70.00){\circle*{1.00}}
\put(45.00,70.00){\circle*{1.00}}
\put(48.00,70.00){\circle*{1.00}}
\put(51.00,70.00){\circle*{1.00}}
\put(54.00,70.00){\circle*{1.00}}
\put(42.00,70.00){\line(0,1){4,5}}
\put(45.00,70.00){\line(0,1){4,5}}
\put(48.00,70.00){\line(0,1){4,5}}
\put(51.00,70.00){\line(0,1){4,5}}
\put(54.00,70.00){\line(0,1){4,5}}
\put(42.00,74.50){\line(1,0){12}}
\put(55.00,71.00){$\big)$}

%%%%%%%%%%%%%%%%%%%%%%%%%% 2 %%%%%%%%%%%%%%%%%%%%%%%%%%%%%%

%%%%%%%%%%%%%%%%%%%%%%%%%% 11121 %%%%%%%%%%%%%%%%%%%%%%%%%%

\put(-11.00,57.00){\line(1,0){4.00}}
\put(-07.00,57.00){\line(1,0){4.00}}
\put(-03.00,57.00){\line(1,1){4.00}}
\put(01.00,61.00){\line(1,-1){4.00}}
\put(-11.00,57.00){\circle*{1.00}}
\put(-07.00,57.00){\circle*{1.00}}
\put(-03.00,57.00){\circle*{1.00}}
\put(01.00,61.00){\circle*{1.00}}
\put(05.00,57.00){\circle*{1.00}}

%%%%%%%%%%%%%%%%%%% Y %%%%%%%%%%%%%%%%%%%

\put(14.00,61.00){\line(1,0){12}}
\put(14.00,65.00){\line(1,0){12}}
\put(14.00,57.00){\line(1,0){4}}
\put(14.00,57.00){\line(0,1){8}}
\put(18.00,57.00){\line(0,1){8}}
\put(22.00,61.00){\line(0,1){4}}
\put(26.00,61.00){\line(0,1){4}}

\put(15.00,62.00){\small{1}}
\put(19.00,62.00){\small{2}}
\put(23.00,62.00){\small{3}}
\put(15.00,58.00){\small{4}}
%%%%%%%%%%%%%%%%%%%%%%%%%%%%%%%%%%%%%%

\put(33.00,57.50){$-$}
\put(36.00,57.50){$\beta$}
\put(39.00,58.00){$\big($}
\put(42.00,57.00){\circle*{1.00}}
\put(45.00,57.00){\circle*{1.00}}
\put(48.00,57.00){\circle*{1.00}}
\put(51.00,57.00){\circle*{1.25}}
\put(54.00,57.00){\circle*{1.00}}
\put(42.00,57.00){\line(0,1){4,5}}
\put(45.00,57.00){\line(0,1){4,5}}
\put(48.00,57.00){\line(0,1){4,5}}
\put(54.00,57.00){\line(0,1){4,5}}
\put(42.00,61.50){\line(1,0){12}}
\put(55.00,58.00){$\big)$}

%%%%%%%%%%%%%%%%%%%%%%%%%% 3 %%%%%%%%%%%%%%%%%%%%%%%%%%%%

%%%%%%%%%%%%%%% 11211 %%%%%%%%%%%%%%%%%%%%%%%%%%%%%%%%%

\put(-11.00,45.00){\line(1,0){4.00}}
\put(-07.00,45.00){\line(1,1){4.00}}
\put(-03.00,49.00){\line(1,-1){4.00}}
\put(01.00,45.00){\line(1,0){4.00}}
\put(-11.00,45.00){\circle*{1.00}}
\put(-07.00,45.00){\circle*{1.00}}
\put(-03.00,49.00){\circle*{1.00}}
\put(01.00,45.00){\circle*{1.00}}
\put(05.00,45.00){\circle*{1.00}}

%%%%%%%%%%%%%%%% Y %%%%%%%%%%%%%%%%%%%%%%

\put(14.00,49.00){\line(1,0){12}}
\put(14.00,53.00){\line(1,0){12}}
\put(14.00,45.00){\line(1,0){4}}
\put(14.00,45.00){\line(0,1){8}}
\put(18.00,45.00){\line(0,1){8}}
\put(22.00,49.00){\line(0,1){4}}
\put(26.00,49.00){\line(0,1){4}}

\put(15.00,50.00){\small{1}}
\put(19.00,50.00){\small{2}}
\put(23.00,50.00){\small{4}}
\put(15.00,46.00){\small{3}}
%%%%%%%%%%%%%%%%%%%%%%%%%%%%%%%%%%%%%%

\put(33.00,45.50){$-$}
\put(36.00,45.50){$\beta$}
\put(39.00,46.00){$\big($}
\put(42.00,45.00){\circle*{1.00}}
\put(45.00,45.00){\circle*{1.00}}
\put(48.00,45.00){\circle*{1.25}}
\put(51.00,45.00){\circle*{1.00}}
\put(54.00,45.00){\circle*{1.00}}
\put(42.00,45.00){\line(0,1){4,5}}
\put(45.00,45.00){\line(0,1){4,5}}
\put(51.00,45.00){\line(0,1){4,5}}
\put(54.00,45.00){\line(0,1){4,5}}
\put(42.00,49.50){\line(1,0){12}}
\put(55.00,46.00){$\big)$}

%%%%%%%%%%%%%%%%%%%%%%%%% 4 %%%%%%%%%%%%%%%%%%%%%%%%%%%%
 
%%%%%%%%%%%%%% 12111 %%%%%%%%%%%%%%%%%%%%%

\put(-11.00,33.00){\line(1,1){4.00}}
\put(-07.00,37.00){\line(1,-1){4.00}}
\put(-03.00,33.00){\line(1,0){4.00}}
\put(01.00,33.00){\line(1,0){4.00}}
\put(-11.00,33.00){\circle*{1.00}}
\put(-07.00,37.00){\circle*{1.00}}
\put(-03.00,33.00){\circle*{1.00}}
\put(01.00,33.00){\circle*{1.00}}
\put(05.00,33.00){\circle*{1.00}}

%%%%%%%%%%%%%%%%%%% Y %%%%%%%%%%%%%%%%%%%

\put(14.00,37.00){\line(1,0){12}}
\put(14.00,41.00){\line(1,0){12}}
\put(14.00,33.00){\line(1,0){4}}
\put(14.00,33.00){\line(0,1){8}}
\put(18.00,33.00){\line(0,1){8}}
\put(22.00,37.00){\line(0,1){4}}
\put(26.00,37.00){\line(0,1){4}}

\put(15.00,38.00){\small{1}}
\put(19.00,38.00){\small{3}}
\put(23.00,38.00){\small{4}}
\put(15.00,34.00){\small{2}}
%%%%%%%%%%%%%%%%%%%%%%%%%%%%%%%%%%%%%%

\put(33.00,33.50){$-$}
\put(36.00,33.50){$\beta$}
\put(39.00,34.00){$\big($}
\put(42.00,33.00){\circle*{1.00}}
\put(45.00,33.00){\circle*{1.25}}
\put(48.00,33.00){\circle*{1.00}}
\put(51.00,33.00){\circle*{1.00}}
\put(54.00,33.00){\circle*{1.00}}
\put(42.00,33.00){\line(0,1){4,5}}
\put(48.00,33.00){\line(0,1){4,5}}
\put(51.00,33.00){\line(0,1){4,5}}
\put(54.00,33.00){\line(0,1){4,5}}
\put(42.00,37.50){\line(1,0){12}}
\put(55.00,34.00){$\big)$}

%%%%%%%%%%%%%%%%%%%%% 5 %%%%%%%%%%%%%%%%%%%%%%%%%%%%%%%%%

%%%%%%%%%%%% 12121 %%%%%%%%%%%%%%%%%%%%%%%%%%%%%%%%%%

\put(-11.00,21.00){\line(1,1){4.00}}
\put(-07.00,25.00){\line(1,-1){4.00}}
\put(-03.00,21.00){\line(1,1){4.00}}
\put(01.00,25.00){\line(1,-1){4.00}}
\put(-11.00,21.00){\circle*{1.00}}
\put(-07.00,25.00){\circle*{1.00}}
\put(-03.00,21.00){\circle*{1.00}}
\put(01.00,25.00){\circle*{1.00}}
\put(05.00,21.00){\circle*{1.00}}

%%%%%%%%%%%%%%%% Y %%%%%%%%%%%%%%%%%%%%%%

\put(14.00,29.00){\line(1,0){8}}
\put(14.00,25.00){\line(1,0){8}}
\put(14.00,21.00){\line(1,0){8}}
\put(14.00,21.00){\line(0,1){8}}
\put(18.00,21.00){\line(0,1){8}}
\put(22.00,21.00){\line(0,1){8}}

\put(15.00,26.00){\small{1}}
\put(19.00,26.00){\small{3}}
\put(15.00,22.00){\small{2}}
\put(19.00,22.00){\small{4}}

%%%%%%%%%%%%%%%%%%%%%%%%%%%%%%%%%%%%%%

\put(36.00,21.50){$\beta$}
\put(39.00,22.00){$\big($}
\put(42.00,21.00){\circle*{1.00}}
\put(45.00,21.00){\circle*{1.25}}
\put(48.00,21.00){\circle*{1.00}}
\put(51.00,21.00){\circle*{1.25}}
\put(54.00,21.00){\circle*{1.00}}
\put(42.00,21.00){\line(0,1){4,5}}
\put(48.00,21.00){\line(0,1){4,5}}
\put(54.00,21.00){\line(0,1){4,5}}
\put(42.00,25.50){\line(1,0){12}}
\put(55.00,22.00){$\big)$}

%%%%%%%%%%%%%%%%%%%%%% 6 %%%%%%%%%%%%%%%%%%%%%%%%%

%%%%%%%%%%%%%%%%%% 11221 %%%%%%%%%%%%%%%%%%%%%%

\put(-11.00,5.00){\line(1,0){4.00}}
\put(-07.00,5.00){\line(1,1){4.00}}
\put(-03.00,9.00){\line(1,0){4.00}}
\put(01.00,9.00){\line(1,-1){4.00}}
\put(-11.00,5.00){\circle*{1.00}}
\put(-07.00,5.00){\circle*{1.00}}
\put(-03.00,9.00){\circle*{1.00}}
\put(01.00,9.00){\circle*{1.00}}
\put(05.00,5.00){\circle*{1.00}}

%%%%%%%%%%%%%%%%%%%%%%%%%%%%%%%%%%%%%%

\put(14.00,17.00){\line(1,0){8}}
\put(14.00,13.00){\line(1,0){8}}
\put(14.00,09.00){\line(1,0){4}}
\put(14.00,05.00){\line(1,0){4}}
\put(14.00,05.00){\line(0,1){12}}
\put(18.00,05.00){\line(0,1){12}}
\put(22.00,13.00){\line(0,1){4}}

\put(15.00,14.00){\small{1}}
\put(19.00,14.00){\small{2}}
\put(15.00,10.00){\small{3}}
\put(15.00,6.00){\small{4}}
%%%%%%%%%%%%%%%%%%%%%%%%%%%%%%%%%%%%%%

\put(33.00,5.50){$-$}
\put(36.00,5.50){$\beta$}
\put(39.00,6.00){$\big($}
\put(42.00,5.00){\circle*{1.00}}
\put(45.00,5.00){\circle*{1.00}}
\put(48.00,5.00){\circle*{1.25}}
\put(51.00,5.00){\circle*{1.25}}
\put(54.00,5.00){\circle*{1.00}}
\put(42.00,5.00){\line(0,1){4,5}}
\put(45.00,5.00){\line(0,1){4,5}}
\put(48.00,5.00){\line(0,1){2,5}}
\put(51.00,5.00){\line(0,1){2,5}}
\put(54.00,5.00){\line(0,1){4,5}}
\put(42.00,9.50){\line(1,0){12}}
\put(48.00,7.50){\line(1,0){3}}
\put(55.00,6.00){$\big)$}

\put(57.00,5.50){$+$}

\put(60.00,5.50){$\beta$}
\put(63.00,6.00){$\big($}
\put(66.00,5.00){\circle*{1.00}}
\put(69.00,5.00){\circle*{1.00}}
\put(72.00,5.00){\circle*{1.25}}
\put(75.00,5.00){\circle*{1.25}}
\put(78.00,5.00){\circle*{1.00}}
\put(66.00,5.00){\line(0,1){4,5}}
\put(69.00,5.00){\line(0,1){4,5}}
\put(78.00,5.00){\line(0,1){4,5}}
\put(66.00,9.50){\line(1,0){12}}
\put(79.00,6.00){$\big)$}

%%%%%%%%%%%%%%  7 %%%%%%%%%%%%%%%%%%%%%%%%%%%%%%%%%%%%
%%%%%%%%%%%%%% 12211 %%%%%%%%%%%%%%%%%%%

\put(-11.00,-11.00){\line(1,1){4.00}}
\put(-07.00,-07.00){\line(1,0){4.00}}
\put(-03.00,-07.00){\line(1,-1){4.00}}
\put(01.00,-11.00){\line(1,0){4.00}}
\put(-11.00,-11.00){\circle*{1.00}}
\put(-07.00,-07.00){\circle*{1.00}}
\put(-03.00,-07.00){\circle*{1.00}}
\put(01.00,-11.00){\circle*{1.00}}
\put(05.00,-11.00){\circle*{1.00}}

%%%%%%%%%%%%%%%%%%%%%%%%%%%%%%%%%%%%%%

\put(14.00,01.00){\line(1,0){8}}
\put(14.00,-03.00){\line(1,0){8}}
\put(14.00,-07.00){\line(1,0){4}}
\put(14.00,-11.00){\line(1,0){4}}
\put(14.00,-11.00){\line(0,1){12}}
\put(18.00,-11.00){\line(0,1){12}}
\put(22.00,-03.00){\line(0,1){4}}

\put(15.00,-2.00){\small{1}}
\put(19.00,-2.00){\small{4}}
\put(15.00,-6.00){\small{2}}
\put(15.00,-10.00){\small{3}}

%%%%%%%%%%%%%%%%%%%%%%%%%%%%%%%%%%%%%%

\put(33.00,-10.50){$-$}
\put(36.00,-10.50){$\beta$}
\put(39.00,-10.00){$\big($}
\put(42.00,-11.00){\circle*{1.00}}
\put(45.00,-11.00){\circle*{1.25}}
\put(48.00,-11.00){\circle*{1.25}}
\put(51.00,-11.00){\circle*{1.00}}
\put(54.00,-11.00){\circle*{1.00}}
\put(42.00,-11.00){\line(0,1){4,5}}
\put(45.00,-11.00){\line(0,1){2,5}}
\put(48.00,-11.00){\line(0,1){2,5}}
\put(51.00,-11.00){\line(0,1){4,5}}
\put(54.00,-11.00){\line(0,1){4,5}}
\put(42.00,-6.50){\line(1,0){12}}
\put(45.00,-8.50){\line(1,0){3}}
\put(55.00,-10.00){$\big)$}

\put(57.00,-10.50){$+$}

\put(60.00,-10.50){$\beta$}
\put(63.00,-10.00){$\big($}
\put(66.00,-11.00){\circle*{1.00}}
\put(69.00,-11.00){\circle*{1.25}}
\put(72.00,-11.00){\circle*{1.25}}
\put(75.00,-11.00){\circle*{1.00}}
\put(78.00,-11.00){\circle*{1.00}}
\put(66.00,-11.00){\line(0,1){4,5}}
\put(75.00,-11.00){\line(0,1){4,5}}
\put(78.00,-11.00){\line(0,1){4,5}}
\put(66.00,-6.50){\line(1,0){12}}
\put(79.00,-10.00){$\big)$}

%%%%%%%%%%%%%%  8  %%%%%%%%%%%%%%%%%%%%%%%%%%%%%%%%%%%%%%%%%%
%%%%%%%%%%%%%%%%%12221%%%%%%%%%%%%%%%%%%

\put(-11.00,-27.00){\line(1,1){4.00}}
\put(-07.00,-23.00){\line(1,0){4.00}}
\put(-03.00,-23.00){\line(1,0){4.00}}
\put(01.00,-23.00){\line(1,-1){4.00}}
\put(-11.00,-27.00){\circle*{1.00}}
\put(-07.00,-23.00){\circle*{1.00}}
\put(-03.00,-23.00){\circle*{1.00}}
\put(01.00,-23.00){\circle*{1.00}}
\put(05.00,-27.00){\circle*{1.00}}

%%%%%%%%%%%%%%%%%%%%%%%%%%%%%%%%%%%%%%

\put(14.00,-15.00){\line(1,0){8}}
\put(14.00,-19.00){\line(1,0){8}}
\put(14.00,-23.00){\line(1,0){4}}
\put(14.00,-27.00){\line(1,0){4}}
\put(14.00,-27.00){\line(0,1){12}}
\put(18.00,-27.00){\line(0,1){12}}
\put(22.00,-19.00){\line(0,1){4}}

\put(15.00,-18.00){\small{1}}
\put(19.00,-18.00){\small{3}}
\put(15.00,-22.00){\small{2}}
\put(15.00,-26.00){\small{4}}

%%%%%%%%%%%%%%%%%%%%%%%%%%%%%%%%%%%%%%

\put(33.00,-26.50){$-$}
\put(36.00,-26.50){$\beta$}
\put(39.00,-26.00){$\big($}
\put(42.00,-27.00){\circle*{1.00}}
\put(45.00,-27.00){\circle*{1.25}}
\put(48.00,-27.00){\circle*{1.25}}
\put(51.00,-27.00){\circle*{1.25}}
\put(54.00,-27.00){\circle*{1.00}}
\put(42.00,-27.00){\line(0,1){4,5}}
\put(45.00,-27.00){\line(0,1){2,5}}
\put(48.00,-27.00){\line(0,1){2,5}}
\put(51.00,-27.00){\line(0,1){2,5}}
\put(54.00,-27.00){\line(0,1){4,5}}
\put(42.00,-22.50){\line(1,0){12}}
\put(45.00,-24.50){\line(1,0){6}}
\put(55.00,-26.00){$\big)$}

\put(57.00,-26.50){$+$}

\put(61.00,-26.50){$\beta$}
\put(64.00,-26.00){$\big($}
\put(67.00,-27.00){\circle*{1.00}}
\put(70.00,-27.00){\circle*{1.25}}
\put(73.00,-27.00){\circle*{1.25}}
\put(76.00,-27.00){\circle*{1.25}}
\put(79.00,-27.00){\circle*{1.00}}
\put(67.00,-27.00){\line(0,1){4,5}}
\put(70.00,-27.00){\line(0,1){2,5}}
\put(73.00,-27.00){\line(0,1){2,5}}
\put(79.00,-27.00){\line(0,1){4,5}}
\put(67.00,-22.50){\line(1,0){12}}
\put(70.00,-24.50){\line(1,0){3}}
\put(80.00,-26.00){$\big)$}

\put(82.00,-26.50){$+$}

\put(86.00,-26.50){$\beta$}
\put(89.00,-26.00){$\big($}
\put(92.00,-27.00){\circle*{1.00}}
\put(95.00,-27.00){\circle*{1.25}}
\put(98.00,-27.00){\circle*{1.25}}
\put(101.00,-27.00){\circle*{1.25}}
\put(104.00,-27.00){\circle*{1.00}}
\put(92.00,-27.00){\line(0,1){4,5}}
\put(98.00,-27.00){\line(0,1){2,5}}
\put(101.00,-27.00){\line(0,1){2,5}}
\put(104.00,-27.00){\line(0,1){4,5}}
\put(92.00,-22.50){\line(1,0){12}}
\put(98.00,-24.50){\line(1,0){3}}
\put(105.00,-26.00){$\big)$}

\put(107.00,-26.50){$-$}

\put(111.00,-26.50){$\beta$}
\put(114.00,-26.00){$\big($}
\put(117.00,-27.00){\circle*{1.00}}
\put(120.00,-27.00){\circle*{1.25}}
\put(123.00,-27.00){\circle*{1.25}}
\put(126.00,-27.00){\circle*{1.25}}
\put(129.00,-27.00){\circle*{1.00}}
\put(117.00,-27.00){\line(0,1){4,5}}
\put(129.00,-27.00){\line(0,1){4,5}}
\put(117.00,-22.50){\line(1,0){12}}
\put(130.00,-26.00){$\big)$}

%%%%%%%%%%%%%%%%%%%%%% 9 %%%%%%%%%%%%%%%%%%%%%%%%%%
%%%%%%%%%% 12321 %%%%%%%%%%%%%%%%%%%%%%%%%%%%%%%%%

\put(-11.00,-39.00){\line(1,1){4.00}}
\put(-07.00,-35.00){\line(1,1){4.00}}
\put(-03.00,-31.00){\line(1,-1){4.00}}
\put(01.00,-35.00){\line(1,-1){4.00}}
\put(-11.00,-39.00){\circle*{1.00}}
\put(-07.00,-35.00){\circle*{1.00}}
\put(-03.00,-31.00){\circle*{1.00}}
\put(01.00,-35.00){\circle*{1.00}}
\put(05.00,-39.00){\circle*{1.00}}

%%%%%%%%%%%%%%%%% Y %%%%%%%%%%%%%%%%%%%%%

\put(14.00,-31.00){\line(1,0){8}}
\put(14.00,-35.00){\line(1,0){8}}
\put(14.00,-39.00){\line(1,0){8}}
\put(14.00,-39.00){\line(0,1){8}}
\put(18.00,-39.00){\line(0,1){8}}
\put(22.00,-39.00){\line(0,1){8}}

\put(15.00,-34.00){\small{1}}
\put(19.00,-34.00){\small{2}}
\put(15.00,-38.00){\small{3}}
\put(19.00,-38.00){\small{4}}

%%%%%%%%%%%%%%%%%%%%%%%%%%%%%%%%%%%%%%

\put(36.00,-38.50){$\beta$}
\put(39.00,-38.00){$\big($}
\put(42.00,-39.00){\circle*{1.00}}
\put(45.00,-39.00){\circle*{1.25}}
\put(48.00,-39.00){\circle*{1.60}}
\put(51.00,-39.00){\circle*{1.25}}
\put(54.00,-39.00){\circle*{1.00}}
\put(42.00,-39.00){\line(0,1){4,5}}
\put(45.00,-39.00){\line(0,1){2,5}}
\put(51.00,-39.00){\line(0,1){2,5}}
\put(54.00,-39.00){\line(0,1){4,5}}
\put(42.00,-34.50){\line(1,0){12}}
\put(45.00,-36.50){\line(1,0){6}}
\put(55.00,-38.00){$\big)$}

\end{picture}
\caption{Motzkin paths $w\in \mathpzc{M}_{5}$, the corresponding Standard Young Tableaux $Y_w$ and cumulants 
$k_w$ expressed in terms of $\beta(\pi)$.
}
\end{figure}
\noindent

Theorem 10.4 gives the decomposition of free cumulants in terms of scalar-valued Motzkin cumulants. Therefore, if there exists
at least one $k$ such that $\tau(a_k)=1$, then both sides of the equation in Theorem 10.4 vanish for $n\geq 2$ (as in 
Example 10.3). However, Lemma 9.8 suggests that it should be possible to obtain
this property in a more direct way, using the sequence $(p_j)$ as the arguments of the $K_w$, since
\[
e_n=p_1+\ldots + p_n
\]
for any $n$ and $(e_n)$ is the `sequential symmetric approximate identity'. 
This suggests that the vanishing of free cumulants of order greater than one when 
$1$ is among its arguments and the vanishing of Motzkin cumulants of order greater than one 
when $p_j$ is among its arguments (a special case of this was proved in Lemma 9.8) are closely related to each other.
Indeed, each $p_j$ plays the role of a `partial identity' (in particular, 
$p_1$ is the `Boolean identity'). As concerns the $k_w$, projections $p_j$ do not appear as their 
arguments since the information about colors is hidden in the word $w$.

\begin{Corollary}
Let $a_1, \ldots, a_n\in \mathcal{A}_{1}\cup {\mathcal A}_{2}\cup {\mathbb C}1$ and let $\tau(a_1), \ldots, \tau(a_n)$ 
be their canonical embeddings in ${\mathcal A}_{1}\star {\mathcal A}_{2}$, where $n\geq 2$. Then 
\[
r_{n}(\tau(a_1), \ldots ,\tau(a_{n}))=0
\]
whenever $a_k=1$ for at least one $k$, where $1\leq k \leq n$.
\end{Corollary}
{\it Proof.}
Let us extend Definition 10.2 by allowing $a_1, \ldots, a_n\in \mathcal{A}_{1}\cup \mathcal{A}_{2}\cup \mathbb{C}1$ and 
setting $1(j)=p_j$ for any $j$. Then we get $k_{j}(1)=\zeta(K_{j}(p_j))=\delta_{j,1}$ for 
any $j$ and thus $r_1(1)=k_1(1)=1$. 
Moreover, if $w=j_1\cdots j_n$ and $j_k=j$, then
\[
\begin{aligned}
&k_w(a_1,\ldots,a_{k-1},1,a_{k+1},\ldots,a_n)  \\
&\qquad =
\zeta\bigl(K_w(A_1,\ldots,A_{k-1},p_j,A_{k+1},\ldots,A_n)\bigr)=0.
\end{aligned}
\]
by Lemma 9.8, where, for simplicity, we put $A_i=a_i(j_i)$ for $i\neq k$.
Therefore, if $a_k=1$, then 
\[
r_{n}(\tau(a_1), \ldots, \tau(a_{n}))
=
\sum_{w=j_1\cdots j_n\in \mathpzc{M}_{n}}
k_{w}(a_1, \ldots, a_{k-1}, 1, a_{k+1}, \ldots, a_n)=0
\]
for $n>1$, by Theorem 10.4 and Lemma 9.8, which completes the proof.
\hfill $\boxvoid$

\begin{Remark}
{\rm 
One can decompose free cumulants of {\it any} variables 
in terms of their scalar-valued Motzkin cumulants. It suffices 
to use orthogonal replicas of these variables to which we assign only one label. In that case, 
all information about the distributions of variables is contained in the mixed moments in the given state (even if
this state is a free product of states) or in the associated mixed Boolean cumulants, but
there are no cumulants with mixed labels and therefore the associated Motzkin cumulants 
cannot be used to linearize addition of `independent' variables.
For instance,
\[
k_{w}(a_1,a_2,a_3)=\left\{
\begin{array}{rr}\beta_{3}(a_1,a_2,a_3)&{\rm if}\;w=111,\\
-\beta_1(a_2)\beta_{2}(a_1,a_3)&{\rm if}\; w=121,
\end{array}
\right.
\]
for any $a_1,a_2,a_3\in \mathcal{A}$ since $B_{121}(a_1,a_2,a_3)=0$. Thus,
\[
r_3(a_1,a_2,a_3)=k_{111}(a_1,a_2,a_3)+k_{121}(a_1,a_2,a_3),
\]
regardless of the relation between $a_1,a_2,a_3$. Of course, here 
$\beta_1,\beta_2, \beta_3$ are Boolean cumulants associated with $\varphi$, and not with $\varphi_1$ and $\varphi_2$.
}
\end{Remark}

\begin{Remark}
{\rm
Let us also remark that there are no `repetitions' of identical 
$w$-Boolean cumulants when adding $k_{w}(a_1, \ldots, a_n)$ for $w\in \mathpzc{M}_{n}$ 
simply because these cumulants are defined in terms of $K_{w}(a_1(j_1), \ldots, a_n(j_n))$.
This is natural since computing various cumulants of variables $a_1, \dots, a_n$ is based 
on different ways of convolving them and each $w$ indicates how this should be done. As concerns Motzkin cumulants,
they serve as $w$-dependent tools to convolve $a_1, \ldots, a_n$. Of course, 
it may happen, as in the case of one variable considered in Example 10.6, that 
$B_w(a,a,a,a)=0$ for non-constant $w$ and thus all such `repetitions' vanish. Nevertheless, 
`repetitions' are needed when we consider orthogonal replicas with non-identical labels since 
otherwise the Motzkin cumulants and thus also scalar-valued Motzkin cumulants with mixed 
labels would not vanish.
}
\end{Remark}
\begin{Remark}
{\rm 
For each $n\in \mathbb{N}$, there is a bijective correspondence $\mathpzc{M}_n\cong {\mathpzc T}_{n}^{(3)}$, where ${\mathpzc T}_{n}^{(3)}$
is the set of standard Young tableaux with $n$ cells and at most three rows \cite{[EFHH]}. Therefore,
the set $\mathpzc{M}_{n}$ can be replaced by ${\mathpzc T}_{n}^{(3)}$ in the decomposition 
of free cumulants of Theorem 10.4. In Fig.~4, we give the explicit 
form of the Standard Young Tableaux $Y_{w}$ with $5$ cells corresponding to $w\in \mathpzc{M}_{5}$, 
using the bijective algorithm constructed in \cite{[EFHH]}, where a more general result has been proved: 
$\mathpzc{M}_{n}^{(d)}\cong {\mathpzc T}_{n}^{(2d+1)}$, where 
$\mathpzc{M}_{n}^{(d)}$ is the set of $d$-colored Motzkin paths and,
for each $q$, ${\mathpzc T}_{n}^{(q)}$ denotes the set of standard Young
tableaux with $n$ cells and at most $q$ rows.
Although this correspondence is used here only as a bijection between
indexing sets, it gives another combinatorial interpretation of the
decomposition in Theorem 10.4. Namely, the scalar-valued Motzkin cumulants
appearing in the decomposition of $r_n$ may equivalently be indexed by
standard Young tableaux with $n$ cells and at most three rows.

}
\end{Remark}

\section{Convolutions}

Motzkin paths are useful in decomposing the free additive convolution 
of distributions. 
In fact, one of our original motivations was to represent the moments of the free additive convolution as a sum of partial free convolutions expressible in terms of partial free cumulants.

\begin{Remark}
{\rm Our notation is similar to that of Nica \cite{[N1]}.
\begin{enumerate}
\item
For an index set $I$, let ${\mathbb C}\langle \{X_i\}_{i\in I}\rangle$ 
be the unital algebra of polynomials in noncommutative indeterminates $(X_{i})_{i\in I}$. 
Let
\[
\Sigma_{I}=\{ \mu: {\mathbb C}\langle \{X_{i}\}_{i\in I}\rangle \rightarrow {\mathbb C}\;|\; \mu \;{\rm linear},\;\mu(1)=1\}.
\]
be the associated family of distributions.
\item
Decompose ${\mathbb C}\langle \{X_{i}\}_{i\in I}\rangle$  into a vector space direct 
sum of {\it homogeneous components}
\[
{\mathbb C}\langle \{X_{i}\}_{i\in I}\rangle =\bigoplus_{n=0}^{\infty} H_{n},
\]
where $H_n$ is spanned by all monomials $X_{k_1}\cdots X_{k_n}$ of degree $n$ and 
$H_{0}={\mathbb C}1$.
\item 
An arbitrary polynomial $f$ can then be uniquely decomposed as a sum of homogeneous polynomials
$f_k$, namely
\[
f=\sum_{k=0}^{\infty}f_k,
\]
where $f_k\in H_k$ for all $k$ and only a finite number of terms is nonzero.
\item
If $\mu_1,\mu_2\in \Sigma_{I}$, then their free additive convolution is defined by
\[
(\mu_1\boxplus \mu_2)(f)=(\mu_1\star \mu_2)(f((X_{k,1}+X_{k,2})_{k\in I})),
\] 
where $\mu_1\star \mu_2$ is the free product of the functionals $\mu_1$ and $\mu_2$  
defined on the free product of unital algebras 
${\mathbb C}\langle (X_{k,1})_{k\in I}\rangle\star {\mathbb C}\langle (X_{k,2})_{k\in I}\rangle
\cong {\mathbb C}\langle (X_{h})_{h\in I\times \{1,2\}}\rangle$ with identification of units.
\item
If we are given unital free subalgebras ${\mathcal A}_{1}, {\mathcal A}_{2}$ of a noncommutative probability space $({\mathcal A}, \varphi)$ and families $(a_{k,1})_{k\in I}$ and $(a_{k,2})_{k\in I}$ of elements in ${\mathcal A}_{1}$ and ${\mathcal A}_{2}$, respectively, then
\[
\left(\mu_1\boxplus \mu_2\right) (X_{k_1}\cdots X_{k_n})=
\varphi \left(\left(a_{k_1,1}+a_{k_1,2}\right)\cdots \left(a_{k_n,1}+a_{k_n,2}\right)\right),
\]
where $\mu_1,\mu_2\in \Sigma_{I}$ are the joint distributions of these families, see \cite{[N1]}.
\item
We can derive the formula for the difference between the moments of the free convolution and of the Boolean convolution. Denote
\[
\Delta(X_{k_1}\cdots X_{k_n}):=(\mu_1\boxplus \mu_2)(X_{k_1}\cdots X_{k_n})-(\mu_1\uplus \mu_2)(X_{k_1}\cdots X_{k_n}),
\]
where the moments $(\mu_1\uplus \mu_2)(X_{k_1}\cdots X_{k_n})$ are defined in terms of 
the Boolean product of functionals by the formula
\[
(\mu_1\uplus \mu_2)(f)=(\mu_1\convolution \mu_2)(f((X_{k,1}+X_{k,2})_{k\in I})),
\] 
where $\mu_1\convolution\mu_2$ is the Boolean product of the functionals $\mu_1$ and $\mu_2$  
defined on the free product of ${\mathbb C}\langle (X_{k,1})_{k\in I}\rangle$ 
and ${\mathbb C}\langle (X_{k,2})_{k\in I}\rangle$ without identification of units.
\item
In general, \(\Delta(X_{k_1}\cdots X_{k_n})\) is a sum of mixed products of Boolean cumulants, where by mixed products we mean products of Boolean cumulants involving variables from different algebras.
\item
Using the framework of Motzkin cumulants, we will show 
that the mixed products of Boolean cumulants can be written as linear combinations of
products of `linearized expressions' (they will be related to convolutions indexed by nonconstant
Motzkin paths). Denote the scalar-valued Motzkin cumulants
\[
k_{w}(a_{k_1,i_1}, \ldots, a_{k_n,i_n}):=\zeta (K_{w}(a_{k_1,i_1}(j_1), \ldots, a_{k_n,i_n}(j_n)))
\]
for $w=j_1\ldots j_n$, where $\zeta$ is defined in Theorem 10.1.
\end{enumerate}
}
\end{Remark}

\begin{Example}
{\rm
Before we treat the general case, let us present a simple example of a monomial of order three. 
For simplicity, denote $x_k=a_{k,1}, y_k=a_{k,2}$.
Using Theorem 2.6 and Example 4.8, we obtain
\begin{eqnarray*}
\Delta(X_{k_1}X_{k_2}X_{k_3})&=&\beta_1(x_{k_2})\beta_2(y_{k_1},y_{k_3})+\beta_1(y_{k_2})\beta_2(x_{k_1},x_{k_3})\\
&=&(\beta_1(x_{k_2})+\beta_1(y_{k_2}))(\beta_2(x_{k_1}, x_{k_3})+\beta_2(y_{k_1},y_{k_3}))\\
&&\!\!\!\!\!-\;\beta_1(x_{k_2})\beta_2(x_{k_1},x_{k_3})-\beta_1(y_{k_2})\beta_2(y_{k_1},y_{k_3})\\
&=&
(k_2(x_{k_2})+k_2(y_{k_2}))(k_{11}(x_{k_1}, x_{k_3})+k_{11}(y_{k_1},y_{k_3}))\\
&&\!\!\!\!\!+\;k_{121}(x_{k_1}, x_{k_2}, x_{k_3})+
k_{121}(y_{k_1},y_{k_2},y_{k_3})
\end{eqnarray*}
since 
\begin{eqnarray*}
k_{121}(x_{k_1},x_{k_2},x_{k_3})&=&-\beta_{2}(x_{k_1},x_{k_3})\beta_1(x_{k_2}),\\
k_{121}(y_{k_1},y_{k_2},y_{k_3})&=&-\beta_{2}(y_{k_1},y_{k_3})\beta_1(y_{k_2}).
\end{eqnarray*}
Thus, the mixed products of Boolean cumulants can be expressed in terms of products of 
linearized scalar-valued Motzkin cumulants.
}
\end{Example}

In order to write $\Delta$ as a sum of convolutions, we will define a family of `Motzkin homogeneous parts' of
the free additive convolution associated with reduced Motzkin words. 
In the convolution context, the situation is quite easy since we only have to 
assign values to monomials of the form $X_{k_1}\cdots X_{k_n}$. Therefore, our definition will parallel 
that in (5) of Remark 11.1. For that purpose, we will use the noncommutative probability space
$({\mathcal A}_{{\rm rep}}, \Phi)$ and the families of orthogonal replicas $(a_{k,1}(j))_{k\in I, j\in \mathbb{N}}$
and $(a_{k,2}(j))_{k\in I, j\in \mathbb{N}}$ with labels 1 and 2, respectively. As compared with the free case, we have 
one additional index $j$ that refers to colors of replicas. 

\begin{Definition}
{\rm
If $\mu_1, \mu_2\in \Sigma_{I}$, then by a {\it Motzkin homogeneous part of $\mu_1 \boxplus \mu_2$ associated with} 
$w\in \mathpzc{M}$ we understand the linear functional on $H_n$ given by the linear extension of
\[
(\mu_1\boxplus_{\,w}\mu_2)(X_{k_1}\cdots X_{k_n}):=\Phi((a_{k_1,1}(j_1)+a_{k_1,2}(j_1))\cdots (a_{k_n,1}(j_n)+a_{k_n,2}(j_n)))
\]
when $|w|=n$, and we set $(\mu_1\boxplus_{\,\varnothing}\mu_2)(c\cdot 1)=c$, for any 
$c\in{\mathbb C}$.
The family 
\[
\mathpzc{H}(\mu_1 , \mu_2):=\{\mu_1\boxplus_{\,w}\mu_2: w\in \mathpzc{M}\}
\] 
will be called the {\it family of Motzkin homogeneous parts} of $\mu_1\boxplus \mu_2$.}
\end{Definition}

Motzkin homogeneous parts of the free additive convolution 
encode finer information about the free additive convolution.
We will derive two formulas for the mixed moments of $\mu_{1}\boxplus_{\,w}\mu_2$: the first one is phrased in terms of Boolean cumulants and uses partitions 
which are monotonically adapted to $w$ and $\ell$ and the second one is a kind of `polarization formula' which expresses these moments in terms of linearized Motzkin cumulants and uses all noncrossing partitions adapted to $w$ and $\ell$.

\begin{Proposition}
Under the assumptions of Theorem 8.4 it holds that
\begin{eqnarray*}
(\mu_1 \boxplus_{\,w} \mu_2)(X_{k_1} \cdots  X_{k_n})
&=&
\sum_{\pi\in \mathcal{M}(w)}\sum_{\ell\in \mathcal{L}_0(\pi)}\beta_{\pi}[a_{k_1,i_1}, \ldots , a_{k_n,i_n}]
\end{eqnarray*}
for any $k_1, \ldots , k_n\in I$, 
where \({\mathcal L}_0(\pi)\) denotes the set of labelings
\(\ell=i_1\cdots i_n\in\{1,2\}^n\)
such that \(\pi\) is monotonically adapted to \(w\) and \(\ell\)
in the sense of Definition 8.1.
\end{Proposition}
{\it Proof.}
By Theorem 8.4, we have, whenever $w=j_1\cdots j_n$,
\begin{eqnarray*}
(\mu_1 \boxplus_{\,w} \mu_2)(X_{k_1}\cdots X_{k_n})
&=& \sum_{\ell\in \{1,2\}^{n}}\zeta(E(a_{k_1,i_1}(j_1) \cdots a_{k_n,i_n}(j_n)))\\
&=& \sum_{\ell\in \{1,2\}^{n}}\sum_{\pi\in {\mathcal Int}(w)}\zeta(B_{\pi}[a_{k_1,i_1}(j_1), \ldots, a_{k_n,i_n}(j_n)])\\
&=& \sum_{\ell\in \{1,2\}^{n}}\sum_{\pi\in \mathcal{M}(w,\ell)}
\beta_{\pi}[a_{k_1,i_1}, \ldots , a_{k_n,i_n}]\\
&=&
\sum_{\pi\in \mathcal{M}(w)}\sum_{\ell\in {\mathcal L}_0(\pi)}
\beta_{\pi}[a_{k_1,i_1}, \ldots ,a_{k_n,i_n}]\\
\end{eqnarray*}
which completes the proof. \hfill $\boxvoid$\\

\begin{Example}
{\rm Let us compute Motzkin convolutions for $f=X_{k_1}X_{k_2}X_{k_3}X_{k_4}$
for all $w\in \mathpzc{M}_{4}$. For simplicity,
we omit the explicit form of $f$ and denote $x_k=a_{k,1}$, $y_k=a_{k,2}$. We obtain
\begin{eqnarray*}
(\mu_1 \boxplus_{\,w_1} \mu_2)(f) &=&(\mu_1\uplus \mu_2) (f), \\
(\mu_1 \boxplus_{\,w_2} \mu_2)(f) &=&\beta_3(x_{k_1},x_{k_2},x_{k_4})\beta_1(y_{k_3})+ \beta_1(x_{k_1})\beta_2(x_{k_2},x_{k_4})\beta_1(y_{k_3})\\
&+& \beta_1(y_{k_1})\beta_2(x_{k_2},x_{k_4})\beta_1(y_{k_3}) +x\rightleftarrows y,\\
(\mu_1 \boxplus_{\,w_3} \mu_2)(f)&=&\beta_3(x_{k_1},x_{k_3},x_{k_4})\beta_1(y_{k_2})+\beta_2(x_{k_1},x_{k_3})\beta_1(y_{k_2})\beta_1(x_{k_4})\\
&+&\beta_2(x_{k_1},x_{k_3})\beta_1(y_{k_2})\beta_1(y_{k_4})+ x\rightleftarrows y,\\
(\mu_1 \boxplus_{\,w_4} \mu_2)(f)&=&\beta_2(x_{k_1},x_{k_4})(\beta_2(y_{k_2},y_{k_3})+\beta_1(y_{k_2})\beta_1(y_{k_3}))+ x\rightleftarrows y,
\end{eqnarray*}
where $w_1=1^4$, $w_2=1^221$, $w_3=121^2$, $w_4=12^21$ and $\mu_1\uplus \mu_2$ is the 
Boolean additive convolution expressed in terms of partitioned Boolean 
cumulants assigned to all interval partitions. The symbol $x\rightleftarrows y$ stands for products obtained by 
exchanging $x$ and $y$. Observe that on the right-hand side the labels alternate in all saturated chains of blocks.
}
\end{Example}

\begin{Example}
{\rm 
In the above example, one can express the RHS of each equation in terms of scalar-valued Motzkin cumulants. 
It is not hard to see that we then obtain the linearized form
\begin{eqnarray*}
(\mu_1\boxplus_{\,w_{s}}\mu_2)(f)&=&k_{w_{s}}(x_{k_1}, x_{k_2}, x_{k_{3}}, x_{k_{4}})+
k_{w_{s}}(y_{k_1}, y_{k_2}, y_{k_{3}}, y_{k_{4}})\\
&&
+\;\;\text{terms involving products of linearized cumulants}.
\end{eqnarray*}
for $s=1,2,3,4$, where the remaining terms are products of linearized expressions corresponding
to blocks of all partitions $\pi\in \mathcal{NC}(w_s)$, respectively. For instance, the sum of such products
for $\mu_1\boxplus_{\,w_4}\mu_2$ takes the form:
\begin{eqnarray*}
&&(k_{121}(x_{k_1},x_{k_2},x_{k_4})+k_{121}(y_{k_1},y_{k_2},y_{k_4}))(k_{2}(x_{k_3})+k_{2}(y_{k_3}))\\
&&+\;(k_{121}(x_{k_1},x_{k_3},x_{k_4})+k_{121}(y_{k_1},y_{k_3},y_{k_4}))(k_{2}(x_{k_2})+k_{2}(y_{k_2}))\\
%\!\!\!\!\!\!\!\!\!\!\!\!\!\!\!
&&+\;(k_{11}(x_{k_1},x_{k_4})+k_{11}(y_{k_1},y_{k_4}))(k_{22}(x_{k_2},x_{k_3})+k_{22}(y_{k_2},y_{k_3}))\\
&&+\;(k_{11}(x_{k_1},x_{k_4})+k_{11}(y_{k_1},y_{k_4}))(k_{2}(x_{k_2})+k_2(y_{k_2}))(k_{2}(x_{k_3})+k_{2}(y_{k_3})),
\end{eqnarray*}
where the partitions can be recovered from the way the indices are written. Note that in contrast to the formulas in Example 11.5 we obtain here all partitions from $\mathcal{NC}(w_4)$.
It is then easy to verify that by adding all these expressions we obtain the corresponding moment 
of the free additive convolution. 
}
\end{Example}

We were able to use the products of linearized scalar-valued Motzkin cumulants above since the length
of all Motzkin words was equal to $4$. Starting from length $5$, operator-valued Motzkin cumulants have to be used
to get nested products of linearized cumulants in all expressions. Therefore, the general formula for mixed 
moments of $\mu_1\boxplus_{\,w}\mu_2$ given below involves $K_{\pi}$ rather than $k_{\pi}$. Linearization
appears in an implicit form in the definition of the set $\mathcal{L}(\pi)$ of
all labelings to which $\pi$ is adapted in the sense of Definition 8.1, which means that the 
labels are identical within blocks.

\begin{Proposition}
For any $k_1, \ldots , k_n\in I$ and any $w=j_1\cdots j_n\in \mathpzc{M}_{n}$, it holds that
\begin{eqnarray*}
(\mu_1 \boxplus_{\,w} \mu_2)(X_{k_1} \cdots  X_{k_n})
&=&
\sum_{\pi\in \mathcal{NC}(w)}
\sum_{\ell\in \mathcal{L}(\pi)}
\zeta(K_{\pi}[a_{k_1,i_1}(j_1), \ldots, a_{k_n,i_n}(j_n)]),
\end{eqnarray*}
where $\mathcal{L}(\pi)$ is the set of labelings $\ell=i_1\cdots i_n$ to which $\pi$ is adapted.
\end{Proposition}
{\it Proof.}
We have
\begin{eqnarray*}
(\mu_1 \boxplus_{\,w} \mu_2)(X_{k_1} \cdots  X_{k_n})
&=& \sum_{\ell\in \{1,2\}^{n}}\zeta(E(a_{k_1,i_1}(j_1) \cdots a_{k_n,i_n}(j_n)))\\
&=& \sum_{\ell\in \{1,2\}^{n}}\sum_{\pi\in {\mathcal NC}(w,\ell)}\zeta(K_{\pi}[a_{k_1,i_1}(j_1), \ldots, a_{k_n,i_n}(j_n)])\\
&=& \sum_{\pi\in \mathcal{NC}(w)}
\sum_{\ell\in \mathcal{L}(\pi)}
\zeta(K_{\pi}[a_{k_1,i_1}(j_1), \ldots, a_{k_n,i_n}(j_n)]),
\end{eqnarray*}
which completes the proof.
\hfill $\boxvoid$\\

\begin{Example}
{\rm In order to see the linearization effect in the formula of Proposition 11.7, let $w=12^31$. In the computation of 
$(\mu_1 \boxplus_{\,w} \mu_2)(f)$, where $f=X_{k_1}X_{k_2}X_{k_3}X_{k_4}X_{k_5}$, the contribution associated with the
partition $\pi=(\pi_0,w)$, where $\pi_0=\{\{1,2,4,5\}, \{3\}\}$, is of the form
\begin{eqnarray*}
&&\zeta(K_{12^21}(x_{k_1},x_{k_2}(K_{2}(x_{k_3})+K_{2}(y_{k_3})),x_{k_4},x_{k_5}))\\
&+&
\zeta(K_{12^21}(y_{k_1},y_{k_2}(K_{2}(x_{k_3})+K_{2}(y_{k_3})),y_{k_4},y_{k_5}))
\end{eqnarray*}
which, due to the fact that $K_{2}(x_{k_3})+K_{2}(y_{k_3})$ is proportional to $p_2$, cannot be
written in the usual product form.
}
\end{Example}

Our last result gives the decomposition of the free additive convolution in terms of
its Motzkin homogeneous parts, which is a natural consequence of Theorem 6.7.
It is also natural to consider the Boolean additive convolution in this context.

\begin{Proposition}
Let \(f=\sum_{m=0}^{\infty} f_m\), where \(f_m\in H_m\) and only finitely many terms are nonzero. Then
Boolean and free additive convolutions have the decompositions
\[
(\mu_1\uplus\mu_2)(f)
=
(\mu_1\boxplus_{\varnothing}\mu_2)(f_0)
+
\sum_{m\geq 1}(\mu_1\boxplus_{\,1^m}\mu_2)(f_m),
\]
and
\[
(\mu_1\boxplus\mu_2)(f)
=
(\mu_1\boxplus_{\varnothing}\mu_2)(f_0)
+
\sum_{m\geq 1}\sum_{w\in\mathpzc M_m}
(\mu_1\boxplus_{\,w}\mu_2)(f_m).
\]
\end{Proposition}
{\it Proof.}
The equation for the Boolean convolution is a consequence of the fact that the families
\((a_{k,1}(1))_{k\in I}\) and \((a_{k,2}(1))_{k\in I}\)
are Boolean independent with respect to \(\Phi\). Therefore, for any \(f_m\in H_m\), \(m\geq 1\),
\[
(\mu_1\uplus \mu_2 )(f_m)
=
(\mu_1 \boxplus_{\,1^m}\mu_2) (f_m).
\]
Together with
\[
(\mu_1\uplus\mu_2)(f_0)
=
(\mu_1\boxplus_{\varnothing}\mu_2)(f_0),
\]
this gives the formula for \((\mu_1\uplus\mu_2)(f)\).

As concerns the free convolution, for any \(k_1, \ldots, k_n\in I\), we obtain, by Theorem 6.7,
\begin{eqnarray*}
(\mu_1 \boxplus \mu_2)(X_{k_1}\cdots X_{k_n})
&=&
\sum_{i_1, \ldots, i_n\in \{1,2\}}
\sum_{j_1\cdots j_n\in \mathpzc{M}_{n}}
\Phi(a_{k_1,i_1}(j_1)\cdots a_{k_n,i_n}(j_n))\\
&=&
\sum_{w\in \mathpzc{M}_{n}}
(\mu_{1}\boxplus_{\,w}\mu_2 )(X_{k_1}\cdots X_{k_n}).
\end{eqnarray*}
By linearity, the same formula holds for all \(f_n\in H_n\), \(n\geq 1\).
Together with
\[
(\mu_1\boxplus\mu_2)(f_0)
=
(\mu_1\boxplus_{\varnothing}\mu_2)(f_0),
\]
we obtain the asserted decomposition for every polynomial
\(f=\sum_{m=0}^{\infty}f_m\).
\hfill $\boxvoid$\\

\section*{Acknowledgments}

The author would like to thank the anonymous referee for valuable comments and suggestions, which were helpful in preparing the revised version of the paper.

\end{document}